\documentclass[11pt,a4paper]{article}
\usepackage{makeidx}
\usepackage{amssymb}
\usepackage[all]{xy}
\usepackage[latin1]{inputenc}
\usepackage{amsfonts}
\usepackage{amsmath}
\usepackage{amssymb}
\usepackage{url}
\usepackage{hyperref}
\usepackage{graphicx}
\usepackage{RR}
\usepackage{RRthemes}

\providecommand{\U}[1]{\protect\rule{.1in}{.1in}}
\newtheorem{theo}{Theorem}[section]
\newtheorem{prop}[theo]{Proposition}
\newtheorem{lem}[theo]{Lemma}
\newtheorem{cor}[theo]{Corollary}

\newtheorem{rem}[theo]{Remark}
\newtheorem{defi}[theo]{Definition}

\numberwithin{equation}{section}

\newcommand{\LL}{\mathbb{L}}

\newcommand{\NN}{\mathbb{N}}

\newcommand{\RR}{\mathbb{R}}

\newcommand{\Ba}{ {\cal B }}

\newcommand{\Da}{ {\cal D }}

\newcommand{\Qa}{ {\cal Q }}

\newcommand{\Ma}{ {\cal M }}

\newcommand{\Pa}{ {\cal P }}

\newcommand{\Ya}{ {\cal Y }}

\newcommand{\point}{\mbox{\LARGE .}}
\newcommand{\proof}{\noindent\mbox{\bf Proof:}\\}
\newcommand{\cqfd}{\hfill\blbx \\}
\def\blbx{\hbox{\vrule height 5pt width 5pt depth 0pt}\medskip}

\def \RR{\mathbb{R}}

\def \LL{\mathbb{L}}

\graphicspath{{./},{./logos/}}
\RRdate{\today}
\RRetitle{On the Stability and the Approximation of Branching Distribution Flows, with 
Applications to Nonlinear Multiple Target Filtering}
\RRtitle{Sur la stabilité et l'approximation de flots de distribution de processus de  branchements, avec applications au filtrage non linéaire multicibles}
\titlehead{Stability and Approximation of Branching Distribution Flows}
\RRauthor{ F. Caron\thanks{%
Centre INRIA Bordeaux et Sud-Ouest \& Institut de Mathématiques de Bordeaux
, Université Bordeaux, 351 cours de la Libération 33405 Talence cedex,
France, Francois.Caron@inria.fr}, P. Del Moral\thanks{%
Pierre.Del-Moral@inria.fr}, M. Pace\thanks{%
Michele.Pace@inria.fr}, B.-N. Vo\thanks{
Ba-Ngu Vo
School of Electrical Electronic \& Computer Engineering
M018 The University of Western Australia,
ba-ngu.vo@uwa.edu.au. The work of B.-N. Vo is supported by Australian Research Council under the
Future Fellowship FT0991854}}
\authorhead{F. Caron, P. Del Moral et M. Pace}
\RRprojet{ALEA}
\RRabstract{We analyse the exponential stability properties of a class of measure-valued equations arising in nonlinear multi-target filtering problems. We also prove the uniform convergence properties 
w.r.t. the time parameter
of a rather general class of stochastic filtering algorithms, including sequential Monte Carlo type models and mean field particle interpretation models. We illustrate these results in the context of the Bernoulli and the
Probability Hypothesis Density filter, yielding what seems to be the first results of this kind in this subject.
}
\RRresume{Nous analysons les propriétés de stabilité exponentielle d'une classe d'\'equations
\`a valeurs mesures que l'on rencontre dans des 
 probl\`emes de filtrage non-linéaire multicibles. 
Nous d\'emontrons ensuite les propri\'et\'es de convergence uniforme, par rapport \`a l'horizon temporel consid\'er\'e, d'une famille assez g\'en\'erale d'algorithmes de filtrage multicibles. Cette analyse s'applique notamment  aux m\'ethodes de type Monte Carlo séquentielles et aux algorithmes particulaires fond\'es sur l'\'evolution de syst\`emes de particules en interaction de type champ moyen. Nous illustrons ces résultats dans le cadre des filtres de Bernoulli et les filtres PHD (Probability Hypothesis Density). Ces r\'esultats semblent \^etre les premiers de ce type pour ces classes de mod\`eles  de filtrage stochastique multicibles.}

\RRkeyword{Measure-valued equations, nonlinear multi-target filtering, Bernoulli filter, Probability hypothesis density filter,
interacting particle systems, particle filters, sequential Monte Carlo methods, exponential concentration inequalities, semigroup stability, functional contraction inequalities.}

\RRmotcle{Processus \`a valeurs mesures, filtrage non-linéaire multicibles, filtre de Bernoulli, filtre PHD  (Probability hypothesis density), filtres particulaires, syst\`emes de particules en interaction, techniques de champ moyen, propri\'et\'es de concentration exponentielle, inégalités de contraction fonctionnelles.} 
\RCBordeaux
\begin{document}
\RRtheme{Mod\`eles et méthodes stochastiques}
\RRNo{7376}

\makeRR   

%
%
%


\tableofcontents
\newpage

\section{Introduction}

Let $(E_{n})_{n\geq 0}$ be a sequence of measurable spaces equipped with the
$\sigma $-fields $(\mathcal{E}_{n})_{n\geq 0}$, and for each with $n\geq 0$,
denote $\mathcal{M}(E_{n})$, $\mathcal{M}_{+}(E_{n})$ and $\mathcal{P}(E_{n})
$ the set of all finite signed measures, the subset of positive measures and
the subset of probability measures, respectively, over the space $E_{n}$.
The aim of this work is to present a stochastic interacting particle
interpretation for numerical solutions of the general measure-valued
dynamical systems $\gamma _{n}\in \mathcal{M}_{+}(E_{n})$ defined by the
following non-linear equation
\begin{equation}
\gamma _{n}(dx_{n})=\left( \gamma _{n-1}Q_{n,\gamma _{n-1}}\right)
(dx_{n}):=\int_{E_{n-1}}\gamma _{n-1}(dx_{n-1})Q_{n,\gamma
_{n-1}}(x_{n-1},dx_{n})  \label{defmod}
\end{equation}%
with initial measure $\gamma _{0}\in \mathcal{M}_{+}(E_{0})$, and positive
and bounded integral operators $Q_{n,\gamma}$ from $E_{n-1}$ into $E_{n}$,
indexed by the time parameter $n\geq 1$ and the set of measures $\gamma\in \mathcal{M}_{+}(E_{n})$.

This class of measure-valued equations arises in
a natural way in the analysis of the first moments evolution of nonlinear branching processes, as well as in signal processing and more particularly in multiple targets tracking models. A pair of filtering models is discussed in some details in section~\ref{bernsecintro} and in
section~\ref{phdmodel1}. In the context of multiple targets tracking problems these measure-valued equations represents the first-order statistical moments of the conditional distributions of the target
occupation measures given observation random measures obscured by clutter, detection uncertainty and data association uncertainty.

As in most of the filtering problems encountered in practice, the initial distribution of the targets
is usually unknown. It is therefore essential to check wether or not the filtering equation "forgets" any erroneous initial distribution. For a thorough discussion on the stability properties of traditional nonlinear filtering
problems with a detailed overview of theoretical developments on this subject, we
 refer to the book~\cite{fk} and to the more recent article by
M. L. Kleptsyna and A. Y. Veretennikov~\cite{vkstab}. Besides the fact that significant progress has been made in the recent years in the rigorous derivation of multiple target tracking nonlinear equations (see for instance~\cite{cddp1,maler,vothesis08,zuev}), up to our knowledge the stability and the robustness properties of these measure-valued models have never been addressed so far in  the literature on the subject. One aim of this paper is to study one such important property: the exponential stability properties of
multiple target filtering models. We present an original and general
perturbation type technique combining the continuity property and the
stability analysis of nonlinear semigroups of the form (\ref{defmod}). A more thorough presentation of these
results is provided in section~\ref{statements} dedicated to the statement of the main results of the present article. The detailed presentation of this perturbation technique can be found in section~\ref{secfctineq}.

On the other hand, while the integral equation (\ref%
{defmod}) appears to be simple at first glance, numerical solutions are
computationally intensive, often requiring integrations in high dimensional
spaces.
One natural way to solve the non-linear integral equation (\ref{defmod}) is
to use find a judicious probabilistic interpretation of the normalized distributions
flow given below
\begin{equation*}
\eta _{n}(dx_{n}):=\gamma _{n}(dx_{n})/\gamma _{n}(1)
\end{equation*}%
To describe with some conciseness these stochastic models, it is important to observe
that the pair process $(\gamma _{n}(1),\eta _{n})\in (\mathbb{R}_{+}\times
\mathcal{P}(E_{n}))$ satisfies an evolution equation of the following form
\begin{equation}
(\gamma _{n}(1),\eta _{n})=\Gamma _{n}(\gamma _{n-1}(1),\eta _{n-1})
\label{flotPhi}
\end{equation}%
Let the mappings $\Gamma _{n}^{1}:\mathbb{R}_{+}\times \mathcal{P}(E_{n})$ $%
\rightarrow $ $\mathbb{R}_{+}$ and $\Gamma _{n}^{2}:\mathbb{R}_{+}\times
\mathcal{P}(E_{n})\rightarrow \mathcal{P}(E_{n})$, denote the first and the
second components of $\Gamma _{n}$ respectively. By construction, we notice
that the total mass process can be computed using the recursive formula
\begin{equation}
\gamma _{n+1}(1)=\gamma _{n}(G_{n,\gamma _{n}})=\eta _{n}(G_{n,\gamma
_{n}})~\gamma _{n}(1)\quad \mbox{\rm with}\quad G_{n,\gamma
_{n}}:=Q_{n+1,\gamma _{n}}(1)  \label{massevol}
\end{equation}

Suppose that we are given an approximation $\left( \gamma _{n}^{N}(1),\eta
_{n}^{N}\right) $ of the pair $(\gamma _{n}(1),\eta _{n})$ at some time
horizon $n$, where $N$ stands for some precision parameter; that is $\left(
\gamma _{n}^{N}(1),\eta _{n}^{N}\right) $ converges (in some sense) to $%
\left( \gamma _{n}(1),\eta _{n}\right) $, as $N\rightarrow \infty$. Then,
the $N$-approximation of the measure $\gamma _{n}$ is given by $\gamma
_{n}^{N}=\gamma _{n}^{N}(1)\times \eta _{n}^{N}$. The central idea behind
any approximation model is to ensure that the total mass process at time $%
(n+1)$ defined by
\begin{equation}
\gamma _{n+1}^{N}(1)=\eta _{n}^{N}(G_{n,\gamma _{n}^{N}})~\gamma _{n}^{N}(1)
\label{massapprox}
\end{equation}%
can be "easily" computed in terms of the $N$-approximation measures
$\gamma _{n}^{N}$. Assuming that the initial mass $\gamma _{0}(1)=\gamma
_{0}^{N}(1)$ is known, the next step is to find some strategy to approximate
the quantities $\Gamma _{n+1}^{2}(\gamma _{n}^{N}(1),\eta _{n}^{N})$ by some
$N$-approximation measures $\eta _{n+1}^{N}$, and to set $\gamma
_{n+1}^{N}=\gamma _{n+1}^{N}(1)\times \eta _{n+1}^{N}$.

The local fluctuations of $\eta _{n}^{N}$ around the
measures $\Gamma _{n}^{2}(\gamma _{n-1}^{N}(1),\eta _{n-1}^{N})$ is
defined in terms of a collection of
random fields $W_{n}^{N}$ :%
\begin{equation}
W_{n}^{N}:=\sqrt{N}~\left[ \eta _{n}^{N}-\Gamma _{n}^{2}(\gamma
_{n-1}^{N}(1),\eta _{n-1}^{N})\right] \Longleftrightarrow \eta
_{n}^{N}=\Gamma _{n}^{2}\left( \gamma _{n-1}^{N}(1),\eta _{n-1}^{N}\right) +%
\frac{1}{\sqrt{N}}~W_{n}^{N}  \label{local}
\end{equation}%
which satisfies for any $r\geq 1$ and any test function $f$ with uniform
norm $\Vert f\Vert \leq 1$,
\begin{equation}
\mathbb{E}\left( W_{n}^{N}(f)~|~\mathcal{F}_{n-1}^{N}\right) =0\quad %
\mbox{and}\quad \mathbb{E}\left( \left\vert W_{n}^{N}(f)\right\vert ^{r}~|~%
\mathcal{F}_{n-1}^{N}\right) ^{\frac{1}{r}}\leq a_{r}  \label{localerr}
\end{equation}%
where $\mathcal{F}_{n-1}^{N}=\sigma \left( \eta _{p}^{N},0\leq p<n\right) $
is the $\sigma $-field generated by the random measures $\eta _{p}^{N}$, $%
0\leq p<n$,\ while $b$ and $a_{r}$ are universal constants whose values do
not depend on the precision parameter $N$. The stochastic analysis of the
resulting particle approximation model relies on the analysis of the
propagation of the local sampling errors defined in (\ref{local}). The main
objective is to control, at any time horizon $n$, the fluctuations of the
random measures $(\gamma _{n}^{N},\eta _{n}^{N})$ around their limiting
values $(\gamma _{n},\eta _{n})$ defined by the following random fields:
\begin{equation}\label{defchamps}
V_{n}^{\gamma ,N}:=\sqrt{N}~\left[ \gamma _{n}^{N}-\gamma _{n}\right] \quad %
\mbox{\rm with}\quad V_{n}^{\eta ,N}:=\sqrt{N}~\left[ \eta _{n}^{N}-\eta _{n}%
\right] .
\end{equation}

The construction of the $N$-approximation measures $\eta _{n}^{N}$
is far from being unique. In the present article, we devise
three different classes of stochastic particle approximation models.
These stochastic algorithms are discussed in section~\ref{secalgo}.
The first one is a mean field particle interpretation of the flow of
probability measures $\eta_n$, and it is presented in section~\ref{secmfield}. The second
model is an interacting particle association model while the third one is a combination of these
two approximation algorithms. These  pair of approximation models are respectively
discussed in section~\ref{secipsasso} and in section~\ref{secipsassomix}. In the context of multi-target
tracking models, the first two approximation models are closely related to the the sequential
Monte Carlo technique presented in the series of articles~\cite{sidenb,arnaud2,william1,william2,william3,zajic},
and respectively, the Gaussian mixture Probability
Hypothesis Density filter discussed in the article by B.-N. Vo, and W.-K. Ma~\cite{voc,vo06}, and the
the Rao-Blackwellized Particle multi-target filters
presented by S. Sarkka, A. Vehtari, and J. Lampinen in~\cite{sarkka1,sarkka2}. These modern stochastic algorithms are rather simple to implement and computationally tractable, and they exhibit excellent performance.

Nevertheless, despite advances in recent years~\cite{clark,arnaud3,arnaud2}, these Monte Carlo particle  type multi-target filters remain poorly understood theoretically. One aim of this article is to present a novel class of stochastic algorithms
with a refined analysis including uniform convergence results w.r.t. the time parameter. We also illustrate these results in the context of multi-target tracking models, yielding what seems to be the first uniform results of this type in this subject.

The rest of the article is organized as follows: In section~\ref{mvpmt} we illustrate
the abstract measure-valued equations (\ref{defmod}) with two recent multi-target filters models, namely
the Bernoulli filter and the
Probability Hypothesis Density filter ({\em abbreviate PHD filter}). Section~\ref{statements}
is devoted to the statement of our main results. In section~\ref{sgdescription}, we describe the semigroups and the continuity properties of the nonlinear equation~\ref{defmod}. We show that this semigroup analysis can be applied to analyse the convergence of the Bernoulli and the PHD approximation filters. Section~\ref{secfctineq} is devoted to the stability properties of nonlinear measure-valued processes
of the form (\ref{flotPhi}). We present a perturbation technique and a
series of functional contraction inequalities. In the next three sections, we illustrate these results in the context of Feynman-Kac models, as well as Bernoulli and PHD models. Section~\ref{secalgo} is concerned with the detailed presentation and the convergence analysis of
three different classes of particle type approximation models, including
 mean field type particle approximations and particle association stochastic algorithms.
Finally, the appendix of the article contains most of technical proofs in the text.

\subsection{Measure-valued systems in Multi-target tracking}\label{mvpmt}

The measure-valued process given by (\ref{defmod}) is a generalisation of
Feynman-Kac measures. Its continuous time version naturally arise in the
modeling and analysis of the first moments of spatial branching process~\cite%
{fk,dynkin}.

Our major motivation for studying this class of measure-valued system stems
from advanced signal processing, more specifically, multiple target
tracking. Driven primarily in the early 1970's by aerospace applications
such as radar, sonar, guidance, navigation, and air traffic control,  today
multi-target filtering has found applications in many diverse disciplines,
see for example the texts \cite{bsf88}, \cite{blackman86} \cite{mahlerbook07}
and references therein. These nonlinear filtering problems deal with jointly
estimating the number and states of several interacting targets given a
sequence of partial observations corrupted by noise, false measurements as
well as miss-detection. This rapidly developing subject is, arguably, one of
the most interesting contact points between the theory of spatial branching
processes, mean field particle systems and advanced signal processing.

The first connections between stochastic branching processes and
multi-target tracking seem to go back to the article by S. Mori, et. al.~%
\cite{mori} published in 1986. However it was Mahler's systematic treatment
of multi-sensor multi-target filtering using random finite sets theory \cite%
{maler4,maler5,maler0,maler} that lead to the development novel multi-target
filters and sparked world wide interests. To motivate the article, we
briefly outline two recent multi-target filters that do not fit the standard
Feynman-Kacs framework, but fall under the umbrella of the measure-valued
equation (\ref{defmod}). The first is the Bernoulli filter for joint
detection and tracking of a single target while the second is the
Probability Hypothesis Density filter.

\subsubsection{Bernoulli filtering}\label{bernsecintro}

A basic problem in target tracking is that the target of interest may not
always be present and exact knowledge of target existence/presence cannot be
determined from observations due to clutter and detection uncertainty \cite%
{mahlerbook07}. The \emph{Bernoulli} filter is a generalisation of the
standard Bayes filter, which accommodates presence and absence of the target
\cite{vothesis08}. In a Bernoulli model, the birth of the target at time $n+1
$ is modelled by a measure $\mu _{n+1}$ on $E_{n+1}$. The target enters the
scene with a probability $\mu _{n+1}(1)<1$ and its state is distributed
according to the normalised measure $\mu _{n+1}/\mu _{n+1}(1)$. At time $n$,
a target $X_{n}$ has a probability $s_{n}(X_{n})$ of surviving to the next
time and evolve to a new state according to a given elementary Markov
transition $M_{n+1}$ from $E_{n}$ into $E_{n+1}$. At time $n+1$, the target
(if it exists) generates with probability $d_{n+1}(X_{n+1})$ an observation $%
Y_{n+1}$ on some auxiliary state space, say $E_{n+1}^{Y}$ with likelihood
function $l_{n+1}(X_{n+1},y)$. This so-called Bernoulli observation point
process is superimposed with an additional and independent Poisson point
process with intensity function $h_{n}>0$ to form the occupation (or counting)
measure observation process $\mathcal{Y}_{n+1}=\sum_{1\leq i\leq
N_{n+1}^{Y}}\delta _{Y_{n+1}^{i}}$.

In its original form, the Bernoulli filter jointly propagates the
probability existence of the target and the distribution of the target state
\cite{vothesis08}. Combining the probability of existence and the state
distribution into a single measure, it can be shown that the Bernoulli
filter satisfies the integral equation (\ref{defmod}), with the probability
of existence of the target given by the mass $\gamma _{n}(1)$ and the
distribution of the target state given by the normalised measure $\eta_n=\gamma
_{n}/\gamma _{n}(1)$. The integral operator for the Bernoulli filter takes
the following form%
\begin{equation}
Q_{n+1,\gamma _{n}}(x_{n},dx_{n+1}):=\frac{
s_{n}(x_{n})g_{n}(x_{n})M_{n+1}(x_{n},dx_{n+1})+(\gamma _{n}(1)^{-1}-1)\mu _{n+1}(dx_{n+1})
}{(1-\gamma _{n}(1))+\gamma
_{n}(g_{n})}
 \label{Bernoullimodel}
\end{equation}%
where $g_{n}$ is a likelihood function given by
\begin{eqnarray}
g_{n}(x_{n}) &:&=(1-d_{n}(x_{n}))+d_{n}(x_{n})\mathcal{Y}_{n}\left({%
l_{n}(x_{n},\cdot )}/{h_{n}}\right)\label{likebe}
\end{eqnarray}

\subsubsection{PHD filtering}\label{phdmodel1}

A more challenging problem arises when the number of targets varies
randomly in time, obscured by clutter, detection uncertainty and
data association
uncertainty. Suppose that at a given time $n$ there are $N_{n}^{X}$ targets $%
(X_{n}^{i})_{1\leq i\leq N_{n}^{X}}$ each taking values in some measurable
state space $E_{n}$. A target $X_{n}^{i}$, at time $n$, survives to the next
time step with probability $s_{n}(X_{n}^{i})$ and evolves to a new state
according to a given elementary Markov transition $M^{\prime}_{n+1}$ from $E_{n}$
into $E_{n+1}$. In addition $X_{n}^{i}$ can spawn new targets at the next
time, usually modelled by a spatial Poisson process with intensity measure $%
B_{n+1}(X_{n}^{i},\cdot )$ on $E_{n+1}$. At the same time, an
independent collection of new targets is added to the current configuration.
This additional and spontaneous branching process is often modeled by a
spatial Poisson process with a prescribed intensity measure $\mu _{n+1}$ on $%
E_{n+1}$. Each target $X_{n+1}^{i}$ generates with probability $%
d_{n+1}(X_{n+1}^{i})$ an observation $Y_{n+1}^{i}$ on some auxiliary state
space, say $E_{n+1}^{Y}$, with probability density function $%
g_{n+1}(X_{n+1}^{i},y)$. In addition to this partial observation point
process we also observe an additional and independent Poisson point process
with intensity function $h_{n}$. Multi-target tracking concerns the
estimation of the random measures $\mathcal{X}_{n+1}=\sum_{1\leq i\leq
N_{n}^{X}}\delta _{X_{n}^{i}}$, given the observation occupation measures $%
\mathcal{Y}_{p}=\sum_{1\leq i\leq N_{p}^{Y}}\delta _{Y_{p}^{i}}.$

The multi-target tracking problem is computationally intractable in general
and the Probability Hypothesis Density PHD (filter), is an approximation
that propagates the first-order statistical moment, or intensity, of the
multi-target state forward in time \cite{maler}. The PHD filter satisfies
the integral equation (\ref{defmod}), with the integral operator given below
\begin{equation}
Q_{n+1,\gamma _{n}}(x_{n},dx_{n+1})=g_{n,\gamma
_{n}}(x_{n})M_{n+1}(x_{n},dx_{n+1})+\gamma _{n}(1)^{-1}~\mu _{n+1}(dx_{n+1})
\label{mttmodel}
\end{equation}%
where $M_{n+1}$ is a Markov kernel defined by%
\begin{equation}
M_{n+1}(x_{n},dx_{n+1}):=\frac{s_{n}(x_{n})M^{\prime}_{n+1}(x_{n},dx_{n+1})+
B_{n+1}(x_{n},dx_{n+1})}{s_{n}(x_{n})+b _{n}(x_{n})}
\end{equation}%
with the branching rate $b_n(x_n)=B_{n+1}(1)(x_n)$. The likelihood function $g_{n,\gamma _{n}}$ is given by
\begin{equation}
g_{n,\gamma _{n}}:=r_{n}\times
\widehat{g}_{n,\gamma _{n}}\quad\mbox{\rm with}\quad
r_n:=\left( s_{n}+b _{n}\right)
\label{decompo}
\end{equation}%
 and
\begin{equation}
\widehat{g}_{n,\gamma _{n}}(x_{n}):=
\left( 1-d_{n}(x_{n})\right)+d_{n}(x_{n})~\int \mathcal{Y}_{n}(dy)~
\frac{g_{n}(x_{n},y_{n})}{%
h_{n}(y_{n})+\gamma _{n}(d_{n}g_{n}(\mbox{\LARGE .},y_{n}))}
\end{equation}
Since its inception by Mahler \cite{maler} in 2003, the PHD filter has
attracted substantial interest to date. The development of numerical
solutions for the PHD filter \cite{arnaud2}, \cite{vo06} have opened the
door to numerous novel extensions and applications. More details on the
derivation of the PHD filter using random finite sets, Poisson techniques or
random measures theoretic approaches can be found in the series of articles~%
\cite{cddp1,maler,zuev}.

\subsection{Statement of the main results}\label{statements}

To describe with some conciseness the main result of this article, we need to
introduce some notation. We let $\mbox{Osc}_{1}(E_{n})$, be the set of $%
\mathcal{E}_n$-measurable functions $f$ on $E_{n}$ with oscillations $%
\mbox{osc}(f)=\sup_{x,x^{\prime }}{|f(x)-f(x^{\prime })|}\leq 1$. We denote
by $\mu (f)=\int ~\mu (dx)~f(x)$ the Lebesgue integral of $f$ w.r.t. some
measure $\mu \in \mathcal{M}(E_{n})$, and we let $\|\mu-\nu\|_{\rm tv}$ be the total variation distance between two probability measures $\nu$ and $\mu$ on $E_n$.

We assume that the following pair of regularity conditions are satisfied.\\

\emph{$(H_1)$ :
There exists a series of compact sets $I_n\subset (0,\infty)$ such that
the initial mass value $\gamma_0(1)\in I_0$, and for any $m\in I_n$ $\eta\in \mathcal{%
\ \mathcal{P }}(E_n)$,
we have
\begin{equation*}
\theta_{-,n}(m)\leq \eta \left(G_{n,m\eta}\right)\leq \theta_{+,n}(m)\quad\mbox{for some pair of positive functions $\theta_{+/-,n}$.}
\end{equation*}
}

The main implication of condition $(H_1)$ comes from the fact that the total
mass processes $\gamma_n(1)$ and their $N$-approximation models $\gamma_n^N(1)$ are finite and they evolves
at every time $n$ in a series of compact sets
$$I_n\subset[m^-_n,m^+_n]\subset (0,\infty)$$
with the sequence of parameters $m^{+/-}_n$ defined by the recursive equations
$m^-_{n+1}=m_{n}^-\theta_{-,n}(m_{n}^-)$ and $m^+_{n+1}=m_{n}^+\theta_{+,n}(m_{n}^+)$, with the initial conditions $m^-_0=m^+_0=\gamma_0(1)$.\\

\emph{$(H_2)$ : For any $n\geq 1$, $f\in\mbox{\rm Osc}_1(E_{n})$, and any $%
(m,\eta),(m^{\prime},\eta^{\prime})\in(I_{n}\times \mathcal{P }(E_{n}))$,
the one step mappings $\Gamma_n=\left(\Gamma_n^1,\Gamma_n^2\right)$ defined
in (\ref{flotPhi}) satisfy the following Lipschitz type inequalities:
\begin{eqnarray}
\left|\Gamma_{n}^1(m,\eta)-\Gamma_{n}^1(m^{\prime},\eta^{\prime})
\right|&\leq&
c(n)~|m-m^{\prime}|+\int~\left|[\eta-\eta^{\prime}](\varphi)\right|~%
\Sigma_{n,(m^{\prime},\eta^{\prime})}^1(d\varphi)  \label{lipGam1} \\
\left|\left[\Gamma_{n}^2(m,\eta)-\Gamma_{n}^2(m^{\prime},\eta^{\prime})%
\right](f) \right|&\leq&
c(n)~|m-m^{\prime}|+\int~\left|[\eta-\eta^{\prime}](\varphi)\right|~%
\Sigma_{n,(m^{\prime},\eta^{\prime})}^2(f,d\varphi)~~~~~~  \label{lipGam2}
\end{eqnarray}
for some finite constants $c(n)<\infty$, and some
 collection of bounded measures $\Sigma^1_{n,(m^{\prime},\eta^{%
\prime})}$ and $\Sigma^2_{n,(m^{\prime},\eta^{\prime})}(f,\mbox{\LARGE .})$
on $\mathcal{B }(E_{n})$ such that
\begin{equation*}
\int~\mbox{\rm osc}(\varphi)~\Sigma^1_{n,(m,\eta)}(d\varphi)\leq
\delta\left(\Sigma_{n}^1\right) \quad\mbox{and}\quad \int~\mbox{\rm osc}%
(\varphi)~\Sigma^2_{n,(m,\eta)}(f,d\varphi)\leq
\delta\left(\Sigma_{n}^2\right)
\end{equation*}
for some finite constant $\delta\left(\Sigma^i_{n}\right)<\infty$, $i=1,2$,
whose values do dot depend on the parameters $(m,\eta)\in(I_{n}\times
\mathcal{P }(E_{n})$ and $f\in\mbox{\rm Osc}_1(E_{n})$. }

Condition $(H_2)$ is a rather basic and weak continuity  type property. It states that the one step
transformations of the flow of measures (\ref{flotPhi}) are weakly Lipschitz, in the sense that the
mass variations and the integral differences w.r.t. some test function $f$ can be controlled by
the different initial masses  and measures w.r.t. a collection of integrals of a possibly infinite number of
test functions.
    It is satisfied for a large class of one step transformations $\Gamma_n$. In section~\ref{introregprop}, we will verify that it is satisfied for the general class of Bernoulli and the PHD filters discussed in section~\ref{bernsecintro} and section~\ref{phdmodel1}.

We are now in position to state the main results of this article. The first one is concerned with the exponential stability properties of the semigroup $
\Gamma_{p,n}=\left(\Gamma_{p,n}^1,\Gamma_{p,n}^2\right)$, with $0\leq p\leq
n $ associated with the one step transformations of the flow (\ref{flotPhi}). A more precise description and the complete proof of the next theorem is provided in section~\ref{secfctineq}.
\begin{theo}\label{theostabintro}
We let
$
\Phi^1_{p,n,\nu}$ and $
\Phi^2_{p,n,m}$
 be the semigroups associated with the one step transformations of the flow
of total masses $\Phi^1_{n,\nu_{n-1}}:=\Gamma^1_{n}\left(\point,\nu_{n-1}\right)$
and measures $\Phi^2_{n,m_{n-1}}:=\Gamma^2_{n}\left(m_{n-1},\point\right)$,
with a fixed collection of measures $\nu:=(\nu_n)_{n\geq 0}\in\prod_{n\geq 0}\Pa(E_n)$ and masses $m:=(m_n)_{n\geq 0}\in \prod_{n\geq 0}I_n$. When these semigroups are exponentially stable (in the sense that they forget exponentially fast  their initial conditions) and when the pair of mappings $\nu_{n-1}\mapsto \Phi^1_{n,\nu_{n-1}}$
and $m_{n-1}\mapsto\Phi^2_{n,m_{n-1}}$ are sufficiently regular then we have the following contraction inequalities
$$
\left|\Gamma_{p,n}^1(u^{\prime},\eta^{\prime})-\Gamma_{p,n}^1(u,\eta)\right|\vee \left\|\Gamma_{p,n}^2(u^{\prime},\eta^{\prime})-\Gamma^2_{p,n}(u,\eta)\right\|_{\rm tv}\leq c~e^{-\lambda (n-p)}
$$
for any $p\leq n$, $u,u^{\prime}\in I_p$, $\eta,\eta^{\prime}\in\Pa(E_p)$, and some finite constants $c<\infty$ and $\lambda>0$ whose values do not depend on the time parameters $p\leq n$.
 \end{theo}

The second theorem is concerned with estimating the approximation error associated with
a $N$-approximation model satisfying condition (\ref{localerr}).  The first part of the theorem is proved in
section~\ref{casgenintro}. The proof of the
uniform estimates is discussed in section~\ref {secfctineqst}
 (see for instance lemma~\ref{lemreff}).

\begin{theo}
\label{casgenintro} Under the assumptions $(H_1)$ and $(H_2)$, the semigroup $%
\Gamma_{p,n}$ satisfies the same Lipschitz type inequalities as those stated in (\ref%
{lipGam1}) and (\ref{lipGam2}) for some collection of measures $%
\Sigma^1_{p,n}$ and $\Sigma^2_{p,n}(f,\mbox{\LARGE .})$ on $\mathcal{B }%
(E_{p})$. In addition, for any $N$-approximation model satisfying condition (%
\ref{localerr}) we have the estimates:
\begin{equation}  \label{estimintro}
\mathbb{E}\left(\left|V^{\gamma,N}_n(1)\right|^r\right)^{\frac{1}{r}} \leq
a_r\sum_{p=0}^n \delta\left(\Sigma^1_{p,n}\right) \quad\mbox{and}\quad
\mathbb{E}\left(\left|V^{\eta,N}_n(f)\right|^r\right)^{\frac{1}{r}} \leq
a_r\sum_{p=0}^n \delta\left(\Sigma^2_{p,n}\right)
\end{equation}
for any $r\geq 1$, and $N\geq 1$, with some constants $a_r<\infty$ whose
values only depend on $r$. Furthermore, under the regularity conditions of
theorem~\ref{theostabintro} the couple of estimates
stated above are uniform w.r.t. the time horizon; that is, we have that $\sup_{n\geq 0}\sum_{p=0}^n \delta\left(\Sigma^i_{p,n}\right)<\infty$, for any $i=1,2$.

\end{theo}

These rather abstract theorems apply to a general class of discrete generation measure-valued equations of the form (\ref{defmod}). We illustrate the application of this pair of theorems in the
analysis of the stability properties and
the approximation convergence of the pair of multiple target filters presented in this introductory section.
These results can basically be stated as follows:

\begin{itemize}
\item The Bernoulli filter presented in section~\ref{bernsecintro} with a sufficiently mixing prediction and almost equal survival and spontaneous births rates $s_n\sim \mu_n(1)$ is exponentially stable.
\item The PHD filter presented in section~\ref{phdmodel1} is exponentially stable
for small clutter intensities and sufficiently high detection probability and spontaneous birth rates.
\item In both situations, the estimation error of any $N$-approximation model satisfying condition (%
\ref{localerr}) does not accumulate over time. Furthermore, the uniform rates of convergence provided in theorem~\ref{casgenintro} allows to design stochastic algorithms with prescribed performance index at any time horizon.
\end{itemize}

We end this section with some direct consequences of theorem~\ref{casgenintro}:

Firstly, we observe that the mean error estimates stated in the above theorem
clearly implies the almost sure convergence results
\begin{equation*}
\lim_{N\rightarrow\infty}\eta_{n}^N(f)=\eta_{n}(f) \quad \mbox{\rm and}
\quad \lim_{N\rightarrow\infty}\gamma_{n}^N(f)=\gamma_{n}(f)
\end{equation*}
for any bounded function $f$ on $E_n$. Furthermore, with some
information on the constants $a_{r}$, these $\mathbb{L}_{r}$-mean error
bounds can be turned to exponential concentration inequalities. To be more
precise, by lemma 7.3.3 in~\cite{fk}, the collection of constants $a_{r}$ in
theorem~\ref{casgenintro}, can be chosen so that
\begin{equation}
a_{2r}^{2r}\leq b^{2r}~(2r)!~2^{-r}/r!\quad \mbox{and}\quad
a_{2r+1}^{2r+1}\leq b^{2r+1} (2r+1)!~2^{-r}/r!  \label{amdef}
\end{equation}%
for some $b<\infty$, whose values do not depend on $r$. Using the above $%
\mathbb{L}_r$-mean error bounds we can establish the following non asymptotic
Gaussian tail estimates:
\begin{equation*}
\mathbb{P}\left(\left|\left[\eta_{n}^N-\eta_{n}\right](f)\right|\geq \frac{%
b_n}{\sqrt{N}}+\epsilon\right)\leq \exp{\left(-\frac{N\epsilon^2}{2b_n^2}%
\right)}\quad\mbox{\rm with}\quad b_n\leq b~\sum_{p=0}^n
\delta\left(\Sigma^2_{p,n}\right)
\end{equation*}
The above result is a direct consequence of the following observation
\begin{equation*}
\forall r\geq 1\qquad \mathbb{E}\left(U^r\right)^{\frac{1}{r}}\leq
a_r~b\Rightarrow \mathbb{P}\left(U\geq b+\epsilon\right)\leq \exp{\left(-{%
\epsilon^2}/{(2b)} \right)}
\end{equation*}
for any non negative random variable $U$. To check this claim, we use the following Laplace estimate
\begin{eqnarray*}
\forall t\geq 0\quad \mathbb{E}\left(e^{tU}\right) &\leq & \exp{\left(\frac{%
(bt)^2}{2}+bt\right)} \Rightarrow \mathbb{P}\left(U\geq
b+\epsilon\right)\leq \exp{\left(- \sup_{t\geq 0}{\left(\epsilon t-\frac{%
(bt)^2}{2}\right)} \right)}
\end{eqnarray*}

It is worth noting that the above constructions allows us to consider with
further work branching particle models in path spaces. These path space
models arise in the analysis of the historical process associated with a
branching models as well as the analysis of a filtering problem of the whole
signal path given a series of observations. For instance, let us suppose
that the Markov transitions $M_{n}$ defined in (\ref{mttmodel} are the
elementary transition of a Markov chain of the following form
\begin{equation*}
X_{n}:=\left( X_{p}^{\prime }\right) _{0\leq p\leq n}\in E_{n}:=\prod_{0\leq
p\leq n}E_{p}^{\prime }
\end{equation*}%
In other words $X_{n}$ represents the paths from the origin up to the
current time of an auxiliary Markov chain $X_{n}^{\prime }$ taking values in
some measurable state spaces $E_{n}^{\prime }$, with Markov transitions $%
M_{n}^{\prime }$. We assume that the potential functions $g_{n,\gamma _{n}}$
only depend on the terminal state of the path, in the sense that $%
g_{n,\gamma _{n}}(X_{n})=g_{n,\gamma _{n}}^{\prime }(X_{n}^{\prime })$, for
some potential function $g_{n,\gamma _{n}}^{\prime }$ on $E_{n}^{\prime }$.
In multiple target tracking problems, these path space models provide a way to
estimate the conditional intensity of the path of a given target in a
multi-target environment related to some likelihood function that only
depends on the terminal state of the signal path.

In practice, it is essential to observe that the mean field particle
interpretations of these path space models simply consist of keeping track
of the whole history of each particle. It can be shown that the resulting
particle model can be interpreted as the genealogical tree model associated
with a genetic type model (see for instance~\cite{fk}). In this situation, $%
\eta_n^N$ is the occupation measure of a random genealogical tree, each
particle represents the ancestral lines of the current individuals.

We end this section with some standard notation used in the paper:

We denote respectively by $\mathcal{M}(E)$, $\mathcal{P}(E)$, and $\mathcal{B%
}(E)$, the set of all finite positive measures $\mu$ on some measurable
space $(E,\mathcal{E})$, the convex subset of all probability measures, and
the Banach space of all bounded and measurable functions $f$ equipped with
the uniform norm $\Vert f\Vert$.
We denote by $f^-$ and $f^+$ the infimum and the supremum
of a function $f$.  For measurable subsets $A\in\mathcal{E}$,
in various instances we slightly abuse notation and we denote $\mu(A)$ instead of $%
\mu(1_{A})$; and we set $\delta_{a}$ the Dirac measure at $a\in E$. We
recall that a bounded and positive integral operator $Q$ from a measurable
space $(E_{1},\mathcal{E}_{1})$ into an auxiliary measurable space $(E_{2},%
\mathcal{E}_{2})$ is an operator $f\mapsto Q(f)$ from $\mathcal{B}(E_{2})$
into $\mathcal{B}(E_{1})$ such that the functions
\begin{equation*}
x\mapsto Q(f)(x):=\int_{E_{2}}Q(x,dy)f(y)
\end{equation*}
are $\mathcal{E}_{1}$-measurable and bounded for some measures $Q(x,%
\mbox{\LARGE .})\in\mathcal{M}(E_{2})$. These operators also generate a dual
operator $\mu\mapsto\mu Q$ from $\mathcal{M}(E_{1})$ into $\mathcal{M}%
(E_{2}) $ defined by $(\mu Q)(f):=\mu(Q(f))$. A Markov kernel is a positive
and bounded integral operator $M$ with $M(1)=1$. We denote by
$Q_{p,n}=Q_{p+1}Q_{p+2}\ldots Q_n$, with
$p\leq n$ the semigroup associated with
a given sequence of
 bounded and positive integral operator $Q_n$ from some measurable
spaces $(E_{n-1},\mathcal{E}_{n-1})$ into  $(E_{n},%
\mathcal{E}_{n})$. For $p=n$, we use the convention $Q_{n,n}=Id$, the identity operator.

We associate with a bounded positive potential function $G:x\in E\mapsto
G(x)\in\lbrack0,\infty)$, the Bayes-Boltzmann-Gibbs transformations
\begin{equation*}
\Psi_{G}~:~\eta\in\mathcal{M}(E)\mapsto\Psi_{G}(\eta)\in\mathcal{P}(E)\quad%
\mbox{\rm with}\quad\Psi_{G}(\eta)(dx):=\frac{1}{\eta(G)}~G(x)~\eta(dx)
\end{equation*}
provided $\eta(G)>0$. We recall that $\Psi _{G}(\eta )$ can be expressed
in terms of a Markov transport equation
\begin{equation}
\eta S_{\eta }=\Psi _{G}(\eta )  \label{TG}
\end{equation}%
for some selection type transition $S_{\eta }(x,dy)$. For instance, we can
take
\begin{equation}
S_{\eta }(x,dy):=\frac{\epsilon }{\eta (G)}~\delta _{x}(dy)+\left( 1-\frac{%
\epsilon }{\eta (G)}\right) ~\Psi _{(G-\epsilon )}(\eta )(dy)  \label{SGn1}
\end{equation}%
for any $\epsilon \geq 0$ s.t. $G(x)\geq \epsilon $. Notice that for $%
\epsilon =0$, we have $S_{\eta }(x,dy)=\Psi _{G}(\eta )(dy)$. We can also
choose
\begin{equation}
S_{\eta }(x,dy):=\epsilon G(x)~\delta _{x}(dy)+\left( 1-\epsilon G(x)\right)
~\Psi _{G}(\eta )(dy)  \label{SGn2}
\end{equation}%
for any $\epsilon \geq 0$ that may depend on the current measure $\eta $,
and s.t. $\epsilon G(x)\leq 1$. For instance, we can choose $1/\epsilon $ to
be the $\eta $-essential maximum of the potential function $G$.
Finally, in the context of Bernoulli and PHD filtering we set  $\overline{\mu}%
_{n+1}=\mu_{n+1}/\mu_{n+1}(1)$, for any $n\geq 0$, the normalized spontaneous birth measures.

\section{Semigroup description}\label{sgdescription}
\subsection{The Bernoulli filter semigroup}\label{bernoulliMcKean}
By construction, we notice that the mass process and the normalized measures are given by the rather simple recursive formulae
\begin{equation}\label{recmassbe}
\gamma_{n+1}(1)=
\frac{\gamma_n(1)\eta_n(g_n)}{(1-\gamma_n(1))+\gamma_n(1)\eta_n(g_n)}~\Psi_{g_n}(\eta_n)(s_n)+
\frac{(1-\gamma_n(1))}{(1-\gamma_n(1))+\gamma_n(1)\eta_n(g_n)}~\mu _{n+1}(1)
\end{equation}
and
$$
\eta_{n+1}:=\alpha_n(\gamma_n)~\Psi_{g_ns_n}(\eta_n)M_{n+1}
+(1-\alpha_n(\gamma_n))~\overline{\mu}_{n+1}
$$
with the mappings $\alpha_n~:~\gamma\in \Ma(E_n)\mapsto \alpha_n(\gamma)\in [0,1]$ defined by
$$
\alpha_n(\gamma)=\frac{\gamma(g_ns_n)}{\gamma(s_ng_n)+(1-\gamma(1))\mu_{n+1}(1)}
$$
By construction, if we set $\gamma=m\times \eta
$ then
\begin{eqnarray*}
\Gamma_{n+1}^1(m,\eta)&=&\frac{\gamma(g_n)}{(1-m)+\gamma(g_n)}~\Psi_{g_n}(\eta)(s_n)+
\frac{(1-m)}{(1-m)+\gamma(g_n)}~\mu _{n+1}(1)\\
\quad \Gamma_{n+1}^2(m,\eta)&=&\Psi_{g_ns_n}(\eta)M_{n+1,\gamma}
\end{eqnarray*}
 with the collection of Markov
transitions $M_{n+1,\gamma}$ defined below
\begin{equation}
M_{n+1,\gamma}(x,\mbox{\LARGE .}):=\alpha_{n}\left(\gamma\right) M_{n+1}(x,%
\mbox{\LARGE .})+\left(1-\alpha_{n}\left(\gamma\right)\right)~\overline{\mu}%
_{n+1}\label{bernoullimck}
\end{equation}
Next we provide an alternative interpretation of the mapping $\Gamma_{n+1}^2$. Firstly, observe that
\begin{equation}\label{bernm}
\Psi_{g_ns_n}(\eta)M_{n+1,\gamma}(f)=\frac{\eta\left(Q_{n+1,m}(f)\right)}{\eta\left(Q_{n+1,m}(1)\right)}
\end{equation}
with the integral operator
$$
Q_{n+1,m}(f)(x):=m g_n(x)s_n(x) M_{n+1}(f)(x)+ (1-m)~\mu_{n+1}(f)
$$
This implies that
$$
\Gamma_{n+1}^2(m,\eta)=\Psi_{\widehat{G}_{n,m}}(\eta)\widehat{M}_{n+1,m}
$$
with the potential function
\begin{equation}\label{bernmG}
\widehat{G}_{n,m}=m g_ns_n+ (1-m)~\mu_{n+1}(1)
\end{equation}
and the Markov transitions
\begin{equation}\label{bernmM}
\widehat{M}_{n+1,m}(f):=\frac{m g_ns_n}{m g_ns_n+ (1-m)~\mu_{n+1}(1)}~
M_{n+1}(f)+ \frac{(1-m)\mu_{n+1}(1)}{m g_ns_n+ (1-m)~\mu_{n+1}(1)}~\overline{\mu}_{n+1}(f)
\end{equation}
The condition $(H_1)$ is clearly not met for the Bernoulli filter (\ref{Bernoullimodel})
when $s_n=0$ and $\mu_{n+1}(1)=0$, since in this situation $\gamma_n=0$ for any $n\geq 1$.
Nevertheless, this condition is met with $I_n\subset (0,1]$ and $m\theta_{+,n}(m)=1$,
as long as $s_n$ and $\mu_{n+1}(1)$ are uniformly bounded from below. It is also met for
$s_n=0$, as long as $0<\mu_{n+1}(1)<1$ and the likelihood function given in (\ref{likebe})
is uniformly bounded. The condition is also met for
$\mu_{n+1}(1)=0$, as long as $\gamma_0(1)>0$, and the likelihood function given in (\ref{likebe})
and the function
$s_{n}$  are uniformly lower bounded.

We prove these assertions using the fact that
\begin{equation}\label{massber}
\gamma_{n+1}(1)=
\widehat{\gamma}_n(1)~\Psi_{g_n}(\eta_n)(s_n)+\left(1-\widehat{\gamma}_n(1)\right)
~\mu _{n+1}(1)
\end{equation}
with the updated mass parameters $\widehat{\gamma}_n(1)\in [0,1]$ given below
$$
\widehat{\gamma}_n(1):=\frac{\gamma_n(1)\eta_n(g_n)}{(1-\gamma_n(1))+\gamma_n(1)\eta_n(g_n)}
$$
If we set $s^-_n:=\inf_{E_{n}}s_{n}$ and $s^+_n=\sup_{E_{n}}s_{n}$
then
$$
\forall n\geq 1\qquad
\gamma_{n}(1)\in \left[m^-_n,m^+_n\right]
$$
with parameters
$$
m^-_n=\mu_{n}(1)\wedge s^-_{n-1}\quad\mbox{\rm and}\quad
m^+_n=\mu_{n}(1)\vee s^+_{n-1}~\left(\leq 1\right)
$$
If $s_n$ and $\mu_{n+1}(1)$ are uniformly bounded from below then we have $m^-_n>0$. In addition,
for the constant mapping $s_n=\mu_{n+1}(1)$, the total mass process is constant $$\gamma_{n+1}(1)=m^+_{n+1}=m^-_{n+1}=\mu_{n+1}(1)$$ for any $n\geq 0$. Furthermore, in this situation the flow of normalized measures is given by the updating-prediction transformation defined by
$$
\forall n\geq 0\qquad \eta_{n+1}=\Psi_{g^{(s)}_n}\left(\eta_n\right)M_{n+1}^{(s)}
$$
with the likelihood function $g^{(s)}_n$ and the Markov transitions $M^{(s)}_{n+1}$ defined by
\begin{equation}\label{gMs}
g^{(s)}_n:=s_ng_n+(1-s_n)\quad\mbox{\rm and}\quad
 M^{(s)}_{n+1}(f):=\frac{s_ng_nM_{n+1}(f)+(1-s_n)~\overline{\mu}_{n+1}(f)}{s_ng_n+(1-s_n)}
\end{equation}
When $\mu_{n+1}(1)=0$, the flow of normalized measures is again given by a simple updating-prediction
 equation
 \begin{equation}\label{ex1}
 \eta_{n+1}=\Psi_{g_ns_n}(\eta_n)M_{n+1}
 \quad\mbox{\rm and}
 \quad \gamma_{n+1}(1)=\Psi_{g_n}(\eta_n)(s_n)\times\theta_{\eta_n(g_n)}(\gamma_n(1))
 \end{equation}
 with the increasing mappings $\theta_a$ defined below
 \begin{equation}\label{theta}
 x\in [0,1]\mapsto \theta_a(x):=ax/[ax+(1-x)]
  \end{equation}
In addition, if $s^-_{n}>0$ then
$$
m^-_{n+1}\geq s^-_{n} \times \frac{g^-_{n}m^-_{n}}{g^-_{n}m^-_{n}+(1-m^-_{n})}>0
$$
as long as $g^-_{n}:=\inf_{E_{n}}g_{n}>0$, and $\gamma_0(1)>0$.
We prove this inequality using the fact that the mapping $(a,x)\in [0,\infty[\times [0,1]\mapsto \theta_a(x)$ is increasing in both coordinates. In the case where $s_n=1$, using the fact that and $\theta_a\circ\theta_b=\theta_{ab}$, we prove that
$$
\gamma_{n+1}(1)=\theta_{\eta_n(g_n)}\left(\gamma_n(1)\right)=\theta_{\prod_{p=0}^n\eta_p(g_p)}(\gamma_0(1))
$$

Conversely, when $\gamma_0(1)<1$ and $0<\mu_{n+1}(1)<1$ and $s_{n}=0$, for any $n\geq 0$, then we have a constant flow of normalized measures
$$
\forall n\geq 1\qquad \eta_{n}=\overline{\mu}_{n}
$$
and the total mass process is such that
$$
\gamma_{n}(1)\in ]0,1[~\Longrightarrow~
\gamma_{n+1}(1)=\mu_{n+1}(1)\times\left[1-\theta_{\overline{\mu}_{n}(g_n)}(\gamma_n(1))\right]\in ]0,1[
$$
with the convention $\overline{\mu}_{0}=\eta_0$, for $n=0$.
In addition, if $\mu_{n+1}(1)=1$  then we have
$$
\gamma_{2(n+1)}(1)=\theta_{\prod_{p=0}^n(b_{2p}/b_{2p+1})}(\gamma_0(1))
\quad\mbox{\rm and}\quad
\gamma_{2n+1}(1)=\theta_{b_{2n}^{-1}\prod_{p=0}^{n-1}(b_{2p+1}/b_{2p})}(\gamma_0(1))
$$
for any $n\geq 0$, with the parameters $b_n:=\overline{\mu}_{n}(g_n)$. We prove these formuae using the  the fact that $1-\theta_a(x)=\theta_{1/a}(1-x)$, and $\theta_a\circ\theta_b=\theta_{ab}$. This again implies that
$m^-_n>0$ as long as $\gamma_0(1)>0$ and the likelihood function
  are uniformly lower bounded.

\subsection{The PHD filter semigroup}\label{phdMcKean}

By construction,  if we set $\gamma=m\times \eta
$ then we find that
\begin{equation*}
\Gamma_{n+1}^1(m,\eta)=\gamma(g_{n,\gamma})+\mu_{n+1}(1) \quad \mbox{and}%
\quad \Gamma_{n+1}^2(m,\eta)=\Psi_{g_{n,\gamma} }(\eta)M_{n+1,\gamma} \end{equation*}
In the above display, $M_{n+1,\gamma}$ is the collection of Markov
transitions defined below
\begin{equation*}
M_{n+1,\gamma}(x,\mbox{\LARGE .}):=\alpha_{n}\left(\gamma\right) M_{n+1}(x,%
\mbox{\LARGE .})+\left(1-\alpha_{n}\left(\gamma\right)\right)~\overline{\mu}%
_{n+1}\quad\mbox{\rm with}\quad
\alpha_{n}\left(\gamma\right)= \frac{\gamma(g_{n,\gamma})}{%
\gamma(g_{n,\gamma})+\mu_{n+1}(1)}
\end{equation*}

The interpretation of the updating transformation $\Psi_{g_{n,\gamma} }(\eta)$
in terms of a Markov transport equation is non unique. For instance, using (%
\ref{decompo}) this Bolzmann-Gibbs transformation can be decomposed into two
parts. The first one relates to the undetectable targets and the second is
associated with non clutter observations. An alternative description is
provided below. We consider a virtual auxiliary observation point $c$
(corresponding to undetectable targets) and set $\mathcal{Y }_n^c=
\mathcal{Y }_n+\delta_c$. We also denote by $g^c_{n,\gamma}(\mbox{\LARGE .}%
,y)$ the function defined below
\begin{equation*}
g^{\gamma}_{n}(\mbox{\LARGE .},y)=\left\{
\begin{array}{lcl}
r_n(1-d_n) & \mbox{\rm if} & y=c \\
r_n~\displaystyle\frac{d_ng_n(\mbox{\LARGE .},y_n)}{h_n(y)+\gamma(d_ng_n(%
\mbox{\LARGE .},y))} & \mbox{\rm if} & y\not=c%
\end{array}%
\right.
\end{equation*}
In this notation, the updating transformation $\Psi_{g_{n,\gamma} }(\eta)$
can be rewritten in the following form
\begin{equation*}
\Psi_{g_{n,\gamma} }(\eta)=\Psi_{\overline{g}_{n,\gamma}}(\eta)\quad%
\mbox{\rm with}\quad \overline{g}_{n,\gamma}=\int \mathcal{Y }%
_n^c(dy)~g^{\gamma}_{n}(\mbox{\LARGE .},y)
\end{equation*}
The averaged potential function $\overline{g}_{n,\gamma}$ allows us to measure
the likelihood of signal states w.r.t. the current observation measure $%
\mathcal{Y }_n^c$. Using (\ref{TG}), the Bolzmann-Gibbs transformation $%
\Psi_{\overline{g}_{n,\gamma}}(\eta)$ can be interpreted as non linear
Markov transport equation of the following form
\begin{equation}
\Psi_{\overline{g}_{n,\gamma}}(\eta)=\eta S_{n,\gamma}\quad\mbox{\rm and}%
\quad \Gamma^2(m,\eta)=\eta K_{n+1,\gamma}\quad\mbox{\rm with}\quad K
_{n+1,\gamma}=S_{n,\gamma}M_{n+1,\gamma}\label{phdmck}
\end{equation}
for some Markov transitions $S_{n,\gamma}$ from $E_n$ into itself.

We also notice that condition $(H_1)$ holds as long as the
functions $s_{n},b_{n}$, and $g_{n}(\mbox{\LARGE .},y_{n})$ are uniformly bounded and $\mu_n(1)>0$.
It is also met when $\mu_n(1)=0$, as long as $r_n=(s_n+b_n)$ is uniformly lower bounded and
$\Ya_n\not=0$ or $d_n<1$.

\subsection{Lipschitz regularity properties}\label{introregprop}

Firstly, we mention that condition $(H_2)$ can be replaced by the following
regularity condition:

\emph{$(H_2^{\prime})$ : For any $n\geq 1$, $f\in\mbox{\rm Osc}_1(E_{n})$,
and any $(m,\eta),(m^{\prime},\eta^{\prime})\in(I_{n}\times \mathcal{P }%
(E_{n}))$, the integral operators $Q_{n,m\eta}$ satisfy the following
Lipschitz type inequalities:
\begin{equation}
\left\|Q_{n,m\eta}(f)-Q_{n,m^{\prime}\eta^{\prime}}(f)\right\|\leq
c(n)~|m-m^{\prime}|+\int~\left|[\eta-\eta^{\prime}](\varphi)\right|~%
\Sigma_{n,(m^{\prime},\eta^{\prime})}(f,d\varphi)~~~  \label{lipsgQ}
\end{equation}
for some collection of bounded measures $\Sigma_{n,(m^{\prime},\eta^{%
\prime})}(f,\mbox{\LARGE .})$ on $\mathcal{B }(E_{n})$ such that
\begin{equation*}
\int~\mbox{\rm osc}(\varphi)~\Sigma_{n,(m,\eta)}(f,d\varphi)\leq
\delta\left(\Sigma_{n}\right)
\end{equation*}
for some finite constant $\delta\left(\Sigma_{n}\right)<\infty$, whose
values do dot depend on the parameters $(m,\eta)\in(I_{n}\times \mathcal{P }%
(E_{n}))$ and $f\in\mbox{\rm Osc}_1(E_{n})$.}

We prove $(H_2^{\prime})\Rightarrow (\ref{lipGam1})$ using the
decompositions
\begin{equation*}
m\eta Q_{n,m\eta}-m^{\prime}\eta^{\prime}Q_{n,m^{\prime}\eta^{\prime}}=
m\eta \left[Q_{n,m\eta}-Q_{n,m^{\prime}\eta^{\prime}}\right]+\left[m\eta-
m^{\prime}\eta^{\prime}\right] Q_{n,m^{\prime}\eta^{\prime}}
\end{equation*}
and of course $\left[m\eta- m^{\prime}\eta^{\prime}\right]=\left[m-
m^{\prime}\right]\eta+m^{\prime}\left[\eta- \eta^{\prime}\right]$. To prove $%
(H_2^{\prime})\Rightarrow (\ref{lipGam2})$, we let $\gamma=m\eta$ and $%
\gamma^{\prime}=m^{\prime}\eta^{\prime}$ and we use the decomposition
\begin{equation*}
\left[\Gamma_{n}^2(m,\eta)-\Gamma_{n}^2(m^{\prime},\eta^{\prime})\right](f)=
\frac{1}{\gamma Q_{n,\gamma}(1)}~\left[\gamma
Q_{n,\gamma}-\gamma^{\prime}Q_{n,\gamma^{\prime}}\right]\left(
f-\Gamma_{n}^2(m^{\prime},\eta^{\prime})(f) \right)
\end{equation*}

The Bernoulli filter (\ref{Bernoullimodel}) satisfies $(H_{2}^{\prime })$, as long as the
 likelihood functions $g_n$ given in (\ref{likebe}) are uniformly bounded above. In this situation, (\ref{lipsgQ}) is met with
$$
\left\|Q_{n,m\eta}(f)-Q_{n,m^{\prime}\eta^{\prime}}(f)\right\|\leq
c(n)~|m-m^{\prime}|+c^{\prime}(n)\left|[\eta-\eta^{\prime}](g_n)\right|
$$
for some finite constant $c^{\prime}(n)<\infty$.

The PHD equation satisfies $(H_{2}^{\prime })$, as long as the functions $%
h_{n}(y)+g_{n,y}^{\prime }$ with $g_{n,y}^{\prime }:=d_{n}g_{n}(%
\mbox{\LARGE
.},y)$ are uniformly bounded above and below. To prove this claim, we simply
use the fact that
\begin{equation*}
\left\Vert \widehat{g}_{n,\gamma }-\widehat{g}_{n,\gamma ^{\prime
}}\right\Vert \leq c_{n}\left[ |m^{\prime }-m|+\displaystyle\int \mathcal{Y}%
_{n}(dy)~\left\vert [\eta ^{\prime }-\eta ](g_{n,y}^{\prime })\right\vert %
\right]
\end{equation*}%
This estimate is a direct consequence of the following one
\begin{equation*}
\widehat{g}_{n,\gamma }(x)-\widehat{g}_{n,\gamma ^{\prime }}(x)=\displaystyle%
\int \mathcal{Y}_{n}(dy)~\frac{g_{n,y}^{\prime }(x)}{h_{n}(y)+\gamma
(g_{n,y}^{\prime })}\frac{\left[ \gamma ^{\prime }-\gamma \right]
(g_{n,y}^{\prime })}{h_{n}(y)+\gamma ^{\prime }(g_{n,y}^{\prime })}
\end{equation*}

Next, we provide a pivotal regularity property of the semigroup $%
\left(\Gamma_{p,n}\right)_{0\leq p\leq n}$ associated with the one step
transformations of the flow (\ref{flotPhi}).

\begin{prop}
\label{proplojp} We assume that conditions $(H_1)$ and  $(H_2)$ are satisfied. Then, for any $0\leq p\leq n$, $f\in\mbox{\rm Osc}%
_1(E_{n})$, and any $(m,\eta),(m^{\prime},\eta^{\prime})\in(I_{p}\times
\mathcal{P }(E_{p}))$, we have the following Lipschitz type inequalities:
\begin{eqnarray*}
\left|\Gamma_{p,n}^1(m,\eta)-\Gamma_{p,n}^1(m^{\prime},\eta^{\prime})
\right|&\leq&
c_p(n)~|m-m^{\prime}|+\int~\left|[\eta-\eta^{\prime}](\varphi)\right|~%
\Sigma_{p,n,(m^{\prime},\eta^{\prime})}^1(d\varphi) \\
\left|\left[\Gamma_{p,n}^2(m,\eta)-\Gamma_{p,n}^2(m^{\prime},\eta^{\prime})%
\right](f) \right|&\leq&
c_p(n)~|m-m^{\prime}|+\int~\left|[\eta-\eta^{\prime}](\varphi)\right|~%
\Sigma_{p,n,(m^{\prime},\eta^{\prime})}^2(f,d\varphi)
\end{eqnarray*}
for some finite constants $c_p(n)<\infty$, and some collection of bounded
measures $\Sigma^1_{p,n,(m^{\prime},\eta^{\prime})}$ and $%
\Sigma^2_{p,n,(m^{\prime},\eta^{\prime})}(f,\mbox{\LARGE .})$ on $\mathcal{B
}(E_{p})$ such that
\begin{equation}
\int~\mbox{\rm osc}(\varphi)~\Sigma^1_{p,n,(m,\eta)}(d\varphi)\leq
\delta\left(\Sigma_{p,n}^1\right) \quad\mbox{and}\quad \int~\mbox{\rm osc}%
(\varphi)~\Sigma^2_{p,n,(m,\eta)}(f,d\varphi)\leq
\delta\left(\Sigma_{p,n}^2\right)\label{lescts}
\end{equation}
for some finite constant $\delta\left(\Sigma^i_{p,n}\right)<\infty$, $i=1,2$%
, whose values do dot depend on the parameters $(m,\eta)\in(I_{p}\times
\mathcal{P }(E_{p})$ and $f\in\mbox{\rm Osc}_1(E_{n})$.
\end{prop}

\noindent\mbox{\bf Proof:}\newline
To prove this proposition, we use a backward induction on the parameter $1\leq p\leq n$. For $p=(n-1)$,
we have $\Gamma_{n-1,n}^i=\Gamma_n^i$, with $i=1,2$, so that the desired
result is satisfied for $p=(n-1)$. We further assume that the estimates hold at a given rank $p<n$. To prove the estimates at rank $(p-1)$,
we recall that
\begin{equation*}
\Gamma_{p-1,n}(m,\eta)=\Gamma_{p,n}\left(\Gamma_p(m,\eta)\right)\Rightarrow%
\forall
i=1,2\qquad\Gamma_{p-1,n}^i(m,\eta)=\Gamma_{p,n}^i\left(\Gamma_p(m,\eta)%
\right)
\end{equation*}
Under the induction hypothesis
\begin{eqnarray*}
\left|\Gamma_{p-1,n}^1(m,\eta)-\Gamma_{p-1,n}^1(m^{\prime},\eta^{\prime})%
\right|&=&
\left|\Gamma_{p,n}^1\left(\Gamma_p(m,\eta)\right)-\Gamma_{p,n}^1\left(%
\Gamma_p(m^{\prime},\eta^{\prime})\right)\right| \\
&\leq &c_p(n)~|\Gamma_p^1(m,\eta)-\Gamma_p^1(m^{\prime},\eta^{\prime})| \\
&&\quad+\int~\left|\left[ \Gamma_p^2(m,\eta)-\Gamma_p^2(m^{\prime},\eta^{%
\prime}) \right](\varphi)\right|~\Sigma_{p,n,\Gamma_p(m^{\prime},\eta^{%
\prime})}^1(d\varphi)
\end{eqnarray*}
On the other hand
\begin{equation*}
|\Gamma_p^1(m,\eta)-\Gamma_p^1(m^{\prime},\eta^{\prime})|\leq
c(p)~|m-m^{\prime}|+\int~\left|[\eta-\eta^{\prime}](\varphi)\right|~%
\Sigma_{p,(m^{\prime},\eta^{\prime})}^1(d\varphi)
\end{equation*}
and
\begin{equation*}
\left|\left[ \Gamma_p^2(m,\eta)-\Gamma_p^2(m^{\prime},\eta^{\prime}) \right]%
(\varphi)\right| \leq
c(p)~|m-m^{\prime}|+\int~\left|[\eta-\eta^{\prime}](\psi)\right|~%
\Sigma_{p,(m^{\prime},\eta^{\prime})}^2(\varphi,d\psi)
\end{equation*}
The end of the proof is now clear. The analysis of $\Gamma_{p-1,n}^2$
follows the same line of arguments and is omitted. This ends the proof
of the proposition. \hfill\hbox{\vrule height 5pt width 5pt depth 0pt}%
\medskip \newline

\subsection{Proof of theorem~\protect\ref{casgenintro}}

This section is mainly concerned with the proof of  the couple of estimates (\ref{estimintro})
stated in theorem~\ref{casgenintro}.

We use the decomposition
\begin{eqnarray}
\left(\gamma_n^N(1),\eta_n^N\right)-\left(\gamma_n(1),\eta_n\right)&=&\left[%
\Gamma_{0,n}\left(\gamma_0^N(1),\eta^N_0\right)-\Gamma_{0,n}\left(%
\gamma_0(1),\eta_0\right)\right]  \notag \\
&&+\sum_{p=1}^n \left[\Gamma_{p,n}\left(\gamma^N_p(1),\eta^N_p\right)-%
\Gamma_{p-1,n}\left(\gamma^N_{p-1}(1),\eta_{p-1}^N\right)\right]
\label{keydecomp}
\end{eqnarray}
and the fact that
\begin{eqnarray*}
\Gamma_{p-1,p}\left(\gamma^N_{p-1}(1),\eta_{p-1}^N\right)&=&
\left(\gamma^N_{p}(1),\Gamma_{p-1,p}^2\left(\gamma^N_{p-1}(1),\eta_{p-1}^N%
\right)\right)
\end{eqnarray*}
to show that
\begin{eqnarray*}
\gamma_n^N(1)-\gamma_n(1)&=&\left[\Gamma^1_{0,n}\left(\gamma_0^N(1),\eta^N_0%
\right)-\Gamma^1_{0,n}\left(\gamma_0(1),\eta_0\right)\right] \\
&&\qquad\hskip1cm+\sum_{p=1}^n \left[\Gamma^1_{p,n}\left(\gamma^N_p(1),%
\eta^N_p\right)-\Gamma^1_{p,n}\left(\gamma^N_{p}(1),\Gamma_{p-1,p}^2\left(%
\gamma^N_{p-1}(1),\eta_{p-1}^N\right)\right)\right]
\end{eqnarray*}
Recalling that $\gamma_0^N(1)=\gamma_0(1)$, using proposition~\ref{proplojp}%
, we find that

\begin{equation*}
\sqrt{N}\left\vert \gamma _{n}^{N}(1)-\gamma _{n}(1)\right\vert \leq
\sum_{p=0}^{n}c_{p}(n)~\int ~\left\vert [W_{p}^{N}(\varphi )\right\vert
~\Sigma _{p,n}^{(N,1)}(d\varphi )
\end{equation*}%
with the predictable measure $\Sigma _{p,n}^{(N,1)}=\Sigma _{p,n,(m,\eta
)}^{1}$ associated with the parameters $(m,\eta )=(\gamma _{p}^{N}(1),\Gamma
_{p-1,p}^{2}\left( \gamma _{p-1}^{N}(1),\eta _{p-1}^{N}\right) )$, with $%
0<p\leq n$; and for $p=0$, we set $\Sigma _{0,n}^{(N,1)}=\Sigma
_{0,n,(\gamma _{0}(1),\eta _{0})}$. Combing the generalized Minkowski's
inequality with (\ref{localerr}) we have
\begin{equation*}
\mathbb{E}\left( \left\vert \int ~\left\vert W_{p}^{N}(\varphi )\right\vert
~\Sigma _{p,n}^{(N,1)}(d\varphi )\right\vert ^{r}\left\vert \mathcal{F}%
_{p-1}^{(N)}\right. \right) ^{\frac{1}{r}}\leq a_{r}~\delta \left( \Sigma
_{p,n}^{1}\right)
\end{equation*}%
for some constants $a_{r}$ whose values only depend on the time parameter.
This clearly implies that
\begin{equation*}
\mathbb{E}\left( \left\vert \gamma _{n}^{N}(1)-\gamma _{n}(1)\right\vert
^{r}\right) ^{\frac{1}{r}}\leq a_{r}~\sum_{p=0}^{n}\delta \left( \Sigma
_{p,n}^{1}\right)
\end{equation*}%
The normalized occupation measures can be analyzed in the same way using the
decomposition given below:
\begin{eqnarray*}
\eta _{n}^{N}-\eta _{n} &=&\left[ \Gamma _{0,n}^{2}\left( \gamma
_{0}^{N}(1),\eta _{0}^{N}\right) -\Gamma _{0,n}^{2}\left( \gamma
_{0}(1),\eta _{0}\right) \right] \\
&&\qquad \hskip2cm+\sum_{p=1}^{n}\left[ \Gamma _{p,n}^{2}\left( \gamma
_{p}^{N}(1),\eta _{p}^{N}\right) -\Gamma _{p,n}^{2}\left( \gamma
_{p}^{N}(1),\eta _{p-1}^{N}K_{p,(\gamma _{p-1}^{N}(1),\eta
_{p-1}^{N})}\right) \right]
\end{eqnarray*}%
This ends the proof of the theorem~\ref{casgenintro}. \hfill
\hbox{\vrule
height 5pt width 5pt depth 0pt}\medskip \newline

\section{Functional contraction inequalities}\label{secfctineq}
\subsection{Stability properties}\label{secfctineqst}
This section is concerned with the long time behavior of nonlinear measure-valued processes
of the form (\ref{flotPhi}). The complexity of these models depend in part on the
interaction function between the flow of masses $\gamma_n(1)$ and the flow of
probability measures $\eta_n=\gamma_n/\gamma_n(1)$. One natural way to start the
analysis of these models
is to study the stability properties of the measure-valued semigroup
associated with a fixed flow of masses, and vice versa. These two mathematical objects are defined below.
\begin{defi}
We associate with a flow of masses $m=(m_n)_{n\geq 0}\in \prod_{n\geq 0}I_n$ and probability
measures $\nu:=(\nu_n)_{n\geq 0}\in\prod_{n\geq 0}\Pa(E_n)$ the pair of semigroups
\begin{equation}
\Phi^1_{p,n,\nu}:=\Phi^1_{n,\nu_{n-1}}\circ\ldots\circ \Phi^1_{1,\nu_{0}}\quad\mbox{\rm and}\quad
\Phi^2_{p,n,m}:=\Phi^2_{n,m_{n-1}}\circ\ldots\circ\Phi^2_{1,m_{0}}\label{sgg}
\end{equation}
with $0\leq p\leq n$,
and the one step transformations
\begin{eqnarray*}
\Phi^1_{n,\nu_{n-1}}&:&u\in I_{n-1}\mapsto
\Phi^1_{n,\nu_{n-1}}(u):=\Gamma^1_{n}\left(u,\nu_{n-1}\right)\in I_{n}\\
\Phi^2_{n,m_{n-1}}&:&\eta\in \Pa(E_{n-1})\mapsto
\Phi^2_{n,m_{n-1}}(\eta):=\Gamma^2_{n}\left(m_{n-1},\eta\right)\in \Pa(E_{n})
\end{eqnarray*}
\end{defi}
By construction, using a simple induction on the time parameter $n$, we find that
$$
 (m_0,\nu_0)=(\gamma_0(1),\eta_0)
 \quad\mbox{\rm and}\quad
 \forall n\geq 1 \quad m_n=\Phi^1_{n,\nu_{n-1}}(m_{n-1})\quad\mbox{\rm and}\quad\nu_n=\Phi^2_{n,m_{n-1}}(\nu_{n-1})
$$
$$
\Updownarrow
$$
$$
\forall n\geq 0 \quad (m_n,\nu_n)=(\gamma_n(1),\eta_n)
$$
In the cases that are of particular interest, the semigroups $\Phi^1_{p,n,\nu}$ and $\Phi^2_{p,n,m}$
will have a Feynman-Kac representation. These models are rather well understood. A brief review on their contraction properties is provided in section~\ref{fksg}. Further details
can be found in the monograph~\cite{fk}.
The first basic
regularity property of these models which are needed is the following weak Lipschitz type property :

 {\em
$(\mbox{Lip$(\Phi)$})$ For any $p\leq n$, $u,u^{\prime}\in I_p$, $\eta,\eta^{\prime}\in\Pa(E_p)$ and $f\in \mbox{Osc}_{1}(E_{n})$ the following Lipschitz inequalities
\begin{eqnarray}
\left|\Phi_{p,n,\nu}^1(u)-\Phi_{p,n,\nu}^1(u^{\prime})
\right|&\leq& a^1_{p,n}~|u-u^{\prime}| \label{a1}\\
\left|\left[\Phi_{p,n,m}^2(\eta)-\Phi_{p,n,m}^2(\eta^{\prime})%
\right](f) \right|&\leq&a^2_{p,n} \int~\left|[\eta-\eta^{\prime}](\varphi)\right|~%
\Omega^2_{p,n,\eta^{\prime}}(f,d\varphi)\label{a2}
\end{eqnarray}
for some finite constants $a^i_{p,n}<\infty$, with $i=1,2$, and some collection of
Markov transitions $\Omega^2_{p,n,\eta^{\prime}}$ from $\mbox{Osc}_{1}(E_{n})$ into $\mbox{Osc}_{1}(E_{p})$, with $p\leq n$, whose values only depend on the parameters $p,n$, resp. $p,n$ and $\eta^{\prime}$.}\\

The semigroups  $\Phi^1_{p,n,\nu}$ and $\Phi^2_{p,n,m}$ may or may bot be asymptotically stable depending on
whether $a^i_{p,n}$ tends to $0$, as $(n-p)\rightarrow\infty$. In section~\ref{bersg}
we provide a set of easily checked regularity conditions under which the semigroups associated with the Bernoulli models discussed in~\ref{bernoulliMcKean} are asymptotically stable.

The second step in the study of the stability properties of the semigroups associated with the flow (\ref{flotPhi}) is the following continuity property:\\

{\em
$(\mbox{Cont$(\Phi)$})$
For any $n\geq 1$, $u,u^{\prime}\in I_{n-1}$, $\eta,\eta^{\prime}\in\Pa(E_{n-1})$ and any $f\in\mbox{Osc}_{1}(E_{n})$
\begin{eqnarray}
\left|\Phi_{n,\eta}^1(u)-\Phi_{n,\eta^{\prime}}^1(u)
\right|&\leq& \tau^1_{n}~ \int~\left|[\eta-\eta^{\prime}](\varphi)\right|~%
\Omega^1_{n,\eta^{\prime}}(d\varphi)\label{tau1}\\
\left|\left[\Phi_{n,u}^2(\eta)-\Phi_{n,u^{\prime}}^2(\eta)%
\right](f) \right|&\leq&\tau^2_{n}~ |u-u^{\prime}|\label{tau2}
\end{eqnarray}
for some finite constants $\tau^i_{n}<\infty$, with $i=1,2$, and some collection probability measures $\Omega^1_{n,\nu^{\prime}}$ on $\mbox{Osc}_{1}(E_{n-1})$, whose values only depend on the parameters $n$, resp. $n$ and $\nu^{\prime}$.

}

This elementary continuity condition allows us to enter the contraction properties of the semigroups  $\Phi^1_{p,n,\nu}$ and $\Phi^2_{p,n,m}$ in the stability analysis of the flow of measures (\ref{flotPhi}).
The resulting functional contraction inequalities will be described in terms of the following collection of parameters.

\begin{defi}
When the couple of conditions $(\mbox{Lip$(\Phi)$})$ and $(\mbox{Cont$(\Phi)$})$ stated above are satisfied, for any $i=1,2$ and $p\leq n$ we set
\begin{equation}\label{overa}
\overline{a}^i_{p,n}=\tau^i_{p+1}~a^i_{p+1,n}\qquad
b_{p,n}=\sum_{p<q<n}\overline{a}^{1}_{p,q}~\overline{a}^2_{q,n}
\quad\mbox{\rm and}\quad
b^{\prime}_{p,n}=\sum_{p\leq q<n}a^{1}_{p,q}~\overline{a}^2_{q,n}
\end{equation}
\end{defi}

The main result of this section is the following proposition.

\begin{prop}\label{propsgpn}
If  conditions $(\mbox{Lip$(\Phi)$})$ and $(\mbox{Cont$(\Phi)$})$ are satisfied, then for any $p\leq n$, $u,u^{\prime}\in I_p$, $\eta,\eta^{\prime}\in\Pa(E_p)$ and $f\in \mbox{Osc}_{1}(E_{n})$ we have the following Lipschitz inequalities
\begin{eqnarray*}
\left|\Gamma_{p,n}^1(u^{\prime},\eta^{\prime})-\Gamma_{p,n}^1(u,\eta)\right|
&\leq&
c^{1,1}_{p,n}
 ~|u-u^{\prime}|
 +
 c^{1,2}_{p,n}
\displaystyle\int~\left|[\eta-\eta^{\prime}](\varphi)\right|~%
\Sigma^1_{p,n,u^{\prime},\eta^{\prime}}(d\varphi)\\
\left|\Gamma_{p,n}^2(u^{\prime},\eta^{\prime})(f)-\Gamma^2_{p,n}(u,\eta)(f)\right]
&\leq& c^{2,1}_{p,n}~|u-u^{\prime}|+c^{2,2}_{p,n}
\displaystyle\int~\left|[\eta-\eta^{\prime}](\varphi)\right|~%
\Sigma^2_{p,n,u^{\prime}\eta^{\prime}}(f,d\varphi)
\end{eqnarray*}
 for some probability measures $\Sigma^1_{p,n,u^{\prime},\eta^{\prime}}(d\varphi)$ and
 Markov transitions  $\Sigma^2_{p,n,m^{\prime}\eta^{\prime}}$, with
 the collection of parameters
 \begin{eqnarray*}
c^{1,1}_{p,n}&=& a^1_{p,n}
+\displaystyle\sum_{p\leq q<n}~c^{2,1}_{p,q}~\overline{a}^{1}_{q,n}\qquad
\mbox{and}\qquad
c^{1,2}_{p,n}=\sum_{p\leq q<n}~
 c^{2,2}_{p,q}~\overline{a}^{1}_{q,n}\\
c^{2,1}_{p,n}&=&b^{\prime}_{p,n}+\sum_{l=1}^{n-p}~
 \sum_{p\leq r_1<\ldots r_{l}<n} b^{\prime}_{p,r_1}~\prod_{1\leq k\leq l} b_{r_k,r_{k+1}}\\
c^{2,2}_{p,n}&=&a^{2}_{p,n}+\sum_{l=1}^{n-p}~
 \sum_{p\leq r_1<\ldots r_{l}<n} a^{2}_{p,r_1}~\prod_{1\leq k\leq l} b_{r_k,r_{k+1}}~,\quad\mbox{with the convention $r_{l+1}=n$.}
\end{eqnarray*}
In particular, the collection of parameters $\delta\left(\Sigma^i_{p,n}\right)_{i=1,2}$, $p\leq n$ introduced in (\ref{estimintro}) and (\ref{lescts}) are such that
$$
\delta\left(\Sigma^1_{p,n}\right)\leq c^{1,2}_{p,n}\quad\mbox{\rm and}\quad
\delta\left(\Sigma^2_{p,n}\right)\leq c^{2,2}_{p,n}
$$
\end{prop}

The proof of this proposition is rather technical and it is postponed to section~\ref{propsgpnp}
in the appendix. Now we conclude this section with a direct application of the above estimates.
The proof of the theorem~\ref{theostabintro} stated in the introduction and the uniform estimates
discussed in theorem~\ref{casgenintro}  are a direct consequence of the
 following lemma.

\begin{lem}\label{lemreff}
Suppose that $\tau^i=\sup_{n\geq 1}\tau^i_n<\infty$, and
$a^{i}_{p,n}\leq c_i~e^{-\lambda_i(n-p)}$, for any $p\leq n$, and some finite parameters $c_i<\infty$
and $\lambda_i>0$,
 with $i=1,2$, satisfying the following condition
 $$\lambda_1\not=\lambda_2\quad\mbox{\rm and}
\quad c_1c_2~\tau^1\tau^2\leq\left(1-e^{-\left(\lambda_1\wedge\lambda_2\right)}\right)~\left(e^{-\left(\lambda_1\wedge\lambda_2\right)}-e^{-\left(\lambda_1\vee\lambda_2\right)}\right)$$
Then,  for any $i,j\in\{1,2\}$ we have
$$
c^{i,j}_{p,n}\leq
c^{i,j}~e^{-\lambda(n-p)}
\quad\mbox{\rm with}\quad\lambda=\left(\lambda_1\wedge\lambda_2\right)-\log{\left(1+c\tau^1\tau^2 \frac{e^{\left(\lambda_1\wedge\lambda_2\right)}}{e^{-\left(\lambda_1\wedge\lambda_2\right)}-e^{-\left(\lambda_1\vee\lambda_2\right)}}\right)}>0
$$
and the parameters $c^{i,j}$ defined below
$$
\begin{array}{rclcrll}
c^{2,2}&=&c_2&
 c^{2,1}&=&c_1c_2\tau^2/\left(e^{-(\lambda_1\wedge\lambda_2)}-e^{-(\lambda_1\vee\lambda_2)}\right)\\
 c^{1,1}&=& c_1~\left(1+c^{2,1}\tau^1/(e^{-\lambda}-e^{-\lambda_1})\right)&
 c^{1,2}&=&c_1c_2\tau^1/(e^{-\lambda}-e^{-\lambda_1})
\end{array}
$$
In particular, for any $N$-approximation models $(\gamma_{n}^{N}(1),\eta _{n}^{N})$ of the flow
$(\gamma_n(1),\eta_n)$  satisfying condition
 (\ref{localerr}), the $\LL_r$-mean error estimates presented in (\ref{estimintro}) are uniform w.r.t. the time parameter
 $$
 \sup_{n\geq 0}{\mathbb{E}\left(\left|V^{\gamma,N}_n(1)\right|^r\right)^{\frac{1}{r}}} \leq
a_r~c^{1,2}/(1-e^{-\lambda})\quad\mbox{and}\quad
 \sup_{n\geq 0}{\mathbb{E}\left(\left|V^{\eta,N}_n(f)\right|^r\right)^{\frac{1}{r}} }\leq
a_r~c^{2,2}/(1-e^{-\lambda})
 $$
 with some constants $a_r<\infty$ whose
values only depend on $r$.

\end{lem}

\proof
Under the premise of the lemma
$$
b_{p,n}\leq c\tau \sum_{p<q<n}~e^{-\lambda_1(q-(p+1))}
~e^{-\lambda_2(n-(q+1))}
\quad\mbox{\rm and}\quad
b^{\prime}_{p,n}\leq c\tau^2
\sum_{p\leq q<n}e^{-\lambda_1(q-p)}~e^{-\lambda_2(n-(q+1))}
$$
with $c=c_1c_2$ and $\tau=\tau^1\tau^2$.
We further assume that $\lambda_1>\lambda_2$ and we set $\Delta=|\lambda_1-\lambda_2|$.
$$
b_{p,n}\leq c\tau e^{-\lambda_2((n-1)-(p+1))}\sum_{p<q<n}~e^{-\Delta (q-(p+1))}\leq
c\tau e^{-\lambda_2((n-1)-(p+1))}/(1-e^{-\Delta})
$$
In the same way, if $\lambda_2>\lambda_1$ we have
$$
b_{p,n}\leq c\tau e^{-\lambda_1((n-1)-(p+1))}\sum_{p<q<n}~e^{-\Delta (n-(q+1))}\leq
c\tau e^{-\lambda_1((n-1)-(p+1))}/(1-e^{-\Delta})
$$
This implies that
$$
b_{p,n}\leq c
\tau e^{-(\lambda_1\wedge\lambda_2)((n-1)-(p+1))}/(1-e^{-\Delta})
$$
In much the same way, it can be shown that
\begin{equation}\label{bprime}
b^{\prime}_{p,n}=\leq c\tau^2e^{-(\lambda_1\wedge\lambda_2)((n-1)-p)}/(1-e^{-\Delta})
\end{equation}
We are now in a position to estimate the parameters $c^{i,j}_{p,n}$. Firstly, we observe that
$$
c^{2,2}_{p,n}\leq c_2~e^{-\lambda_2(n-p)}
+c_2\sum_{l=1}^{n-p}~\left(\frac{c\tau^1\tau^2 e^{2(\lambda_1\wedge\lambda_2)}}{1-e^{-\Delta}}\right)^l
 \sum_{p\leq r_1<\ldots r_{l}<n}
 e^{-\lambda_2(r_1-p)}
e^{-(\lambda_1\wedge\lambda_2)
(n-r_1)}
$$
When $\lambda_1>\lambda_2$, we find that
$$
c^{2,2}_{p,n}\leq
c_2~e^{-\lambda_2(n-p)}
\sum_{l=0}^{n-p}~
\left(\frac{c\tau e^{2\lambda_2}}{1-e^{-\Delta}}\right)^l
\left(
\begin{array}{c}
n-p\\
l
\end{array}
\right)
$$
and therefore
$$
c^{2,2}_{p,n}\leq
c_2~e^{-\lambda_2(n-p)}~\left(1+c\tau~\frac{e^{2\lambda_2}}{1-e^{-\Delta}}\right)^{n-p}
\Rightarrow
c^{2,2}_{p,n}=
c_2~e^{-\lambda(n-p)}
$$
with
$$
\lambda=\lambda_2-\log{\left(1+c\tau \frac{e^{\lambda_2}}{e^{-\lambda_2}-e^{-\lambda_1}}\right)}>0
$$
as long as
$$
c\tau \leq\left(1-e^{-\lambda_2}\right)~\left(e^{-\lambda_2}-e^{-\lambda_1}\right)
$$

When $\lambda_2>\lambda_1$ we have $\lambda_2=\lambda_1+\Delta$, we find that
$$
c^{2,2}_{p,n}\leq c_2~e^{-\lambda_2(n-p)}
+c_2e^{-\lambda_1
(n-p)}\sum_{l=1}^{n-p}~\left(\frac{c\tau  e^{2\lambda_1}}{1-e^{-\Delta}}\right)^l
 \sum_{p\leq r_1<\ldots r_{l}<n}
 e^{-\Delta(r_1-p)}
$$
from which it follows that
$$
c^{2,2}_{p,n}\leq c_2~e^{-\lambda_1(n-p)}
\left( 1+c\tau \frac{e^{2\lambda_1}}{1-e^{-\Delta}}\right)^{n-p}
$$
Using a similar line of argument as above, we have
$$
c^{2,2}_{p,n}\leq
c_2~e^{-\lambda(n-p)}
$$
with
$$
\lambda=\lambda_1-\log{\left(1+c\tau \frac{e^{\lambda_1}}{e^{-\lambda_1}-e^{-\lambda_2}}\right)}>0
$$
as long as
$$
c\tau \leq\left(1-e^{-\lambda_1}\right)~\left(e^{-\lambda_1}-e^{-\lambda_2}\right)
$$
We conclude that
$$
c^{2,2}_{p,n}\leq
c_2~e^{-\lambda(n-p)}
$$
with
$$
\lambda=\left(\lambda_1\wedge\lambda_2\right)-\log{\left(1+c\tau \frac{e^{\left(\lambda_1\wedge\lambda_2\right)}}{e^{-\left(\lambda_1\wedge\lambda_2\right)}-e^{-\left(\lambda_1\vee\lambda_2\right)}}\right)}>0
$$
as long as
$$
c\tau \leq\left(1-e^{-\left(\lambda_1\wedge\lambda_2\right)}\right)~\left(e^{-\left(\lambda_1\wedge\lambda_2\right)}-e^{-\left(\lambda_1\vee\lambda_2\right)}\right)
$$
Using (\ref{bprime}) we also show that
$$
c^{2,1}_{p,n}\leq
c^{2,1}~e^{-\lambda(n-p)}
 \quad\mbox{\rm
 with}
 \quad
 c^{2,1}=c\tau^2~\frac{1}{e^{ -(\lambda_1\wedge\lambda_2)}-e^{-(\lambda_1\vee\lambda_2)}}
 $$
 Using these estimates
 $$
 c^{1,1}_{p,n}= c_1~e^{-\lambda_1(n-p)}
+\displaystyle\sum_{p\leq q<n}~c^{2,1}_{p,q}~
c_1\tau^1~e^{-\lambda_1(n-(q+1))}
 $$
 and
 $$
 c^{1,1}_{p,n}= c_1~e^{-\lambda_1(n-p)}
+c^{2,1}c_1\tau^1\displaystyle\sum_{p\leq q<n}~e^{-\lambda(q-p)}
~e^{-\lambda_1(n-(q+1))}
 $$
Since $\lambda_1>\lambda$ we find that
 $$
 c^{1,1}_{p,n}\leq  c_1~e^{-\lambda_1(n-p)}
+c^{2,1}c_1\tau^1
~e^{-\lambda((n-1)-p)}/(1-e^{-\Delta^{\prime}})
 \quad\mbox{\rm
 with}\quad \Delta^{\prime}=\lambda_1-\lambda>0$$
 This yields
 $$
 c^{1,1}_{p,n}\leq c^{1,1}~e^{-\lambda(n-p)}
 \quad\mbox{\rm
 with}\quad
 c^{1,1}:= c_1~\left(1+c^{2,1}\tau^1/(e^{-\lambda}-e^{-\lambda_1})\right)
 $$
 Finally, we observe that
 $$
 c^{1,2}_{p,n}= c\tau^1
 \sum_{p\leq q<n}~e^{-\lambda(q-p)}
~e^{-\lambda_1(n-(q+1))} \leq c\tau^1~e^{-\lambda((n-1)-p)}/(1-e^{-\Delta^{\prime}})
 $$
which implies that
  $$
 c^{1,2}_{p,n}\leq c^{1,2}~e^{-\lambda(n-p)}
 \quad\mbox{\rm
 with}
 \quad
 c^{1,2}:=c\tau^1/(e^{-\lambda}-e^{-\lambda_1})
 $$
 This ends the proof of the lemma.\cqfd

 \subsection{Feynman-Kac models}\label{fksg}
We let $Q_{p,n}$, with $0\leq
p\leq n$, be the Feynman-Kac semi-group associated with a sequence of
 bounded and positive integral operator $Q_n$ from some measurable
spaces $(E_{n-1},\mathcal{E}_{n-1})$ into  $(E_{n},%
\mathcal{E}_{n})$. For any $n\geq 1$, we denote by $G_{n-1}$ and $M_n$ the potential function on $E_{n-1}$ and the Markov transition from $E_{n-1}$ into $E_n$ defined below
$$
G_{n-1}(x)=Q_n(1)(x)\quad\mbox{\rm and}\quad M_{n}(f)(x)=\frac{Q_{n}(f)(x)}{Q_n(1)(x)}
$$
We also denote by $\Phi_{p,n}$, $0\leq p\leq n$, the nonlinear semigroup from $\Pa(E_p)$ into $\Pa(E_n)$
defined below
\begin{equation}
\forall \eta\in\Pa(E_p),~\forall f\in\Ba(E_n)\quad
\Phi_{p,n}(\eta)(f)={\eta Q_{p,n}(f)}/{\eta Q_{p,n}(1)}
\end{equation}
As usual we use the convention $\Phi_{n,n}=Id$, for $p=n$. It is important to observe that this semigroup is
alternatively
defined by  the formulae
$$
\Phi_{p,n}(\eta)(f)=\frac{\eta(G_{p,n}~P_{p,n}(f))}{\eta(G_{p,n})}
\quad
\mbox{\rm with}\quad G_{p,n}=Q_{p,n}(1)\quad\mbox{\rm and}\quad P_{p,n}(f_n)=
{Q_{p,n}(f_n)}/{Q_{p,n}(1)}
$$
The next two parameters
\begin{equation}\label{thetwop}
r_{p,n}=\sup_{x,x^{\prime}\in E_p}{\frac{G_{p,n}(x)}{G_{p,n}(x^{\prime})}}
\quad\mbox{\rm
and}\quad
\beta(P_{p,n})=\sup_{x_p,y_p\in E_p}{
\|P_{p,n}(x_p,\point)-P_{p,n}(y_p,\point)\|_{\tiny\rm tv}}
\end{equation}
measure respectively
the relative oscillations of the potential
functions $G_{p,n}$
and the contraction properties of the Markov transition $P_{p,n}$.
Various  estimates  in the forthcoming sections
will be expressed in terms of these parameters. For instance and for further use in several places in this article, we have the  following Lipschitz regularity property.
\begin{prop}[\cite{dprmft}]\label{keyprop}
For any $f_n\in\mbox{Osc}_1(E_n)$ we have
\begin{equation}\label{lip}
\left|\left[\Phi_{p,n}(\eta_p)-\Phi_{p,n}(\mu_p)\right](f_n)\right|
\leq~2~r_{p,n}~\beta(P_{p,n})~
~\left|\left[\eta_p-\mu_p\right]\overline{P}_{p,n}^{\mu_p}(f_n)\right|
\end{equation}
for some function $\overline{P}_{p,n}^{\mu_p}(f_n)\in\mbox{Osc}_1(E_p)$ that doesn't depends on the measure $\eta_p$.
\end{prop}
Our next objective is to estimate the
 the contraction coefficients $r_{p,n}$ and $\beta(P_{p,n})$  in terms of the mixing type properties of the semigroup
 $$M_{p,n}(x_p,dx_n)
:=M_{p+1} M_{p+2}\ldots M_{n}(x_p,dx_n)$$ associated with the Markov operators $M_n$. We introduce the following regularity condition.\\
\\
{\em $(MG)_m$ \hskip.2cm There exists an integer $m\geq 1$
and a sequence $(\epsilon_p(M))_{p\geq 0}\in (0,1)^{\NN}$ and some
finite constant $r_p$ such that for any $p\geq 0$ and any $(x,x^{\prime})\in E_p^2$ we have
\begin{equation}\label{Hm}
  M_{p,p+m}(x_p,\point)
\geq \epsilon_p(m)~~
M_{p,p+m}(x^{\prime}_p,\point) \quad\mbox{\rm and}\quad
G_p(x)\leq r_{p}~
G_n(x^{\prime})
\end{equation}
}

It is well known that the above condition  is satisfied for any aperiodic and irreducible Markov chains on
finite spaces. Loosely speaking, for non compact spaces this condition is related to the tails
of the transition distributions on the boundaries of the state space. For instance, let us suppose that $E_n=\RR$ and $M_n$ is the bi-Laplace transition given by
$$
M_n(x,dy)=\frac{c(n)}{2}~e^{-c(n)\,|y-A_n(x)|}~dy
$$
for some $c(n)>0$ and some drift function $A_n$ with bounded oscillations $\mbox{\rm osc}(A_n)<\infty$.  In this case, it is readily checked that condition $(M)_m$
holds true for $m=1$ with the parameter
$
\epsilon_{n-1}(1)=\exp{(-c(n)~\mbox{\rm osc}(A_n))}
$.

Under  the mixing type condition $(M)_m$ we have for any $n\geq m\geq 1$, and $p\geq 1$
\begin{equation}\label{rpn}
r_{p,p+n}\leq\epsilon_p(m)^{-1}~\prod_{0\leq k<m}r_{p+k}
\end{equation}
and
\begin{equation}\label{betapn}
\beta(P_{p,p+n})\leq  \prod_{k=0}^{\lfloor n/m\rfloor-1}
\left(1-\epsilon_{p+km}^{(m)}\right)
\quad\mbox{\rm with}\quad
\epsilon_{p}^{(m)}:=\epsilon_{p}^2(m)~\prod_{0<k<m}r_{p+k}^{-1}
\end{equation}
Notice that these estimates are also valid for any $n\geq 0$.
Several contraction inequalities can be deduced from these estimates (see for instance
chapter 4 of the book~\cite{fk}). To give a flavor of these results, we further assume that $(M)_m$ is satisfied with $m=1$, and we have $\epsilon=\inf_n{\epsilon_{n}(1)}>0$. In this case, we can show that
$$r_{p,p+n}\leq r_p/\epsilon\quad \mbox{\rm and} \quad
\beta(P_{p,p+n})\leq  \left(1-\epsilon^2\right)^n
$$
We end this short section with a direct consequence of proposition~\ref{keyprop}.
\begin{cor}\label{propp}
Consider the Bernoulli semigroup presented in section~\ref{bernoulliMcKean}.
For constant mappings $s_n=\mu_{n+1}(1)$,  the first component mapping
is constant $ \Phi^1_{n+1,\nu_n}(u)=s_n$
and the second component mapping $
 \Phi^2_{n+1,m_n}(\eta)=
 \Psi_{g^{(s)}_n}(\eta)M_{n+1}^{(s)}
 $
 induces a Feynman-Kac semigroup
with the likelihood function $g^{(s)}_n$ and the Markov transitions $M^{(s)}_{n+1}$ defined in (\ref{gMs}).
In this situation, the condition  (\ref{a1}) is clearly met with $a^1_{p,n}=0$, for any $p<n$,.
We further assume that the semigroup of associated with the Markov transitions $M_{n}$ satisfies
the mixing property stated in the l.h.s. of (\ref{Hm})
for some integer $m\geq 1$ and some parameter $\epsilon_p(m)\in ]0,1]$.
In this situation,  the condition
(\ref{a2}) is also met with the collection of parameters $a^2_{p,n}$ given below
$$
a^2_{p,n}\leq 2~\rho_{p}(m)~\prod_{k=0}^{\lfloor (n-p)/m\rfloor-1}\left(
1-\epsilon^{(m,s)}_{p+km}
\right)
$$
with
$$
\rho_{p}(m):=\epsilon_p^{-1}(m)~\prod_{p\leq k<p+m}r_k^2(s_k)r_k(1)
~~
\mbox{and}~~
\epsilon^{(s,m)}_{p}= \epsilon_p^2(m)r_p(s_p)/\prod_{p\leq k<p+m}r_k(s_k)^{3}r_k(1)^{2}
$$
and the collection of parameters $r_n(s_n)$ defined below
$$
r_n(s_n):=
\frac{s_ng^+_n+(1-s_n)}{s_ng^-_n+(1-s_n)}\left(\leq r_n(1)\right)
$$
\end{cor}

\subsection{Bernoulli models}\label{bersg}
This section is concerned with the contraction properties of the semigroups
$\Phi_{p,n,\nu}^1$ and $\Phi_{p,n,m}^2$ associated with the Bernoulli filter discussed in section~\ref{bernoulliMcKean}. Before proceeding, we provide a brief discussion on the oscillations of the
 likelihood functions  $g_{n}$
given below
$$
g_{n}(x_{n})=(1-d_{n}(x_{n}))+d_{n}(x_{n})\mathcal{Y}_{n}\left({%
l_{n}(x_{n},\cdot )}/{h_{n}}\right)
$$
in terms of some $[0,1]$-valued detection probability functions $d_n$, some local likelihood functions $l_n$, and some
positive clutter intensity function $h_n$. The oscillations of these  likelihood functions
strongly depend on the nature of the functions $(d_n,h_n,l_n)$.

Assuming that $h^-_n>0$ we have
\begin{equation}\label{boundg}
(1-d^{\circ,-}_n)+d^{\circ,-}_n~\frac{l_n^-}{h_n^+}~\mathcal{Y}_{n}\left(1\right)\leq g_{n}^-\leq g^+_n\leq (1-d^{\circ,+}_n)+d^{\circ,+}_n~\frac{l_n^+}{h_n^-}~\mathcal{Y}_{n}\left(1\right)
\end{equation}
with the parameters
$$
d^{\circ,+}_n=d^+_n~1_{l_n^+\mathcal{Y}_{n}\left(1\right)\geq h^-_n}+
d^-_n~1_{l_n^+\mathcal{Y}_{n}\left(1\right)< h^-_n}
$$
$$
d^{\circ,-}_n=d^-_n~1_{l_n^-\mathcal{Y}_{n}\left(1\right)\geq h^+_n}+
d^+_n~1_{l_n^-\mathcal{Y}_{n}\left(1\right)< h^+_n}
$$
The semigroup contraction inequalities developed in this section will be expressed in terms
of the following parameters
$$
\delta_n(sg):=\frac{g^+_n s^+_n}{g^-_n s^-_n}~, \quad
\delta_n(g):=\frac{g^+_n}{g^-_n}\quad\mbox{\rm and}\quad \delta^{\prime}_n(g):=\frac{1}{g^-_n}\wedge g^+_n
$$
For time homogeneous models  $(d_n,h_n,l_n)=(d,h,l)$, with constant detection probability $d_n(x)=d$ and
uniformly bounded number of observations $\sup_n\Ya_n(1)\leq \Ya^+(1)<\infty$ we have the following estimates
$$
(1-d)\leq g_{n}^-\leq g^+_n\leq (1-d)+d~\frac{l^+}{h^-}~\mathcal{Y}^+\left(1\right)
$$
In this situation, we have
$$
\delta_n(g)\leq 1+\frac{d}{1-d}~\frac{l^+}{h^-}~\mathcal{Y}^+\left(1\right)
$$
For small clutter intensity function with $h^->0$ and $l^->0$ we also have
the observation free estimates $\frac{g^+_n}{g^-_n}\leq \frac{l^+h^+}{l^-h^-}$,
from which we find that the upper bound
\begin{equation}\label{bound1}
 \delta(g):=\sup_{n\geq 0}\delta_n(g)\leq\inf{\left\{1+\frac{d}{1-d}~\frac{l^+}{h^-}~\mathcal{Y}\left(1\right),\frac{l^+h^+}{l^-h^-}\right\}}
\end{equation}
and for $d<1$
\begin{equation}\label{bound2}
 \delta^{\prime}(g):=\sup_{n\geq 0}\delta^{\prime}_n(g)\leq\sup{\left\{(1-d)+d\frac{l^+}{h^-}~\mathcal{Y}\left(1\right),\frac{1}{1-d}\right\}}
\end{equation}
To be more precise, if we set $\inf_n\Ya_n(1)=\Ya^-(1)$ then
$$
1\leq \frac{l^-}{h^+}\Ya(1)^-\Rightarrow \delta^{\prime}(g)\leq (1-d)+d~\frac{l^+}{h^-}~\mathcal{Y}^+\left(1\right)
$$
In addition, if we have $d(1-d)\Ya(1)\leq h^-/l^+$ and $d<1$ then we find the the observation free estimates
$$
d\Ya(1)~l^+/h^-\leq 1/(1-d)\Rightarrow \delta^{\prime}(g)\leq(1-d)+ \frac{1}{1-d}
$$
Conversely, we have the observation free estimates
$$
\frac{l^+}{h^-}\Ya(1)^+\leq 1\Rightarrow \delta^{\prime}(g)\leq \frac{1}{(1-d)+d~\frac{l^-}{h^+}~\mathcal{Y}^-\left(1\right)}\leq \frac{1}{1-d}
$$

We are now in position to state the main result of this section.

\begin{theo}\label{theoBernoulli}
If $\mu_{n+1}(1)\in]0,1[$, $0<s^-_n\leq s_n^+<1$, and the semigroup $M_{p,n}$ satisfies the condition stated in the l.h.s. of (\ref{Hm})
for
some integer $m\geq 1$ and some positive constant $ \epsilon_p(m)$, then
 the condition $(\mbox{Lip$(\Phi)$})$ is met with
$$
a^1_{p,n}\leq 2~\epsilon_p^{-1} \delta^{\prime}_p(g)~\prod_{p\leq k<p+n}\left(1-\epsilon_{k}^{2} \right)
\quad
\mbox{and}
\quad
a^2_{p,n}\leq 2~\rho_p(m)
\prod_{k=0}^{\lfloor n/m\rfloor-1}
\left(1-\epsilon_{p+km}^{(m)}\right)
$$
with some parameters
$$
\epsilon_n\geq \inf{\left\{
\frac{s^-_n}{\mu_{n+1}(1)},\frac{\mu_{n+1}(1)}{s^+_n},
\frac{1-s^+_n}{1-\mu_{n+1}(1)},\frac{1-\mu_{n+1}(1)}{1-s^-_n}
\right\}}
$$
and
$$
\rho_p(m)
\leq \epsilon_p(m)^{-1}~\prod_{0\leq k<m}\delta_{p+k}(sg)^{3}
\quad
\mbox{and}
\quad
\epsilon_{p}^{(m)}\geq
\epsilon_p(m)^2~\delta_{p}(sg)^{-4}\prod_{0<k<m}\delta_{p+k}(sg)^{-5}
$$
In addition condition $(\mbox{Cont$(\Phi)$})$ is met with
$$
\tau_{n+1}^1\leq \delta_n(g)~\left[(s^+_n-s^-_n)+\left\|s_n-\mu_{n+1}(1)\right\|\right]
\quad
\mbox{and}
\quad
\tau_{n+1}^2\leq
\delta^{\prime}_n(g) ~
 \sup{\left\{\frac{\mu_{n+1}}{s^-_n},\frac{s^+_n}{\mu_{n+1}(1)}\right\}}
$$

\end{theo}
The proof of the theorem is postponed to section~\ref{ptheoBernoulli}.
To give a flavour of these estimates we examine time homogeneous models $$(d_n,h_n,l_n,s_n,\mu_n)=(d,h,l,s,\mu)$$
with constant detection and survival probabilities $d_n(x)=d$, $s_n(x)=s$, and
uniformly bounded number of observations $\sup_n\Ya_n(1)\leq \Ya(1)<\infty$. In this situation, we have $(\epsilon_p(m),\epsilon_p^{(s)}(m))=(\epsilon(m),\epsilon^{(s)}(m))$ and using the estimates
(\ref{bound1}) we prove
the following bounds
$$
\tau^1_{n+1}\leq \delta(g)~\left|s-\mu(1)\right|
\quad
\mbox{and}
\quad
\tau^2_{n+1}\leq
\delta^{\prime}(g) ~
\frac{\mu(1)\vee s}{\mu(1)\wedge s}
$$
and
$$
a^1_{0,n}\leq 2\epsilon^{-1} \delta^{\prime}(g)\left(1-\epsilon^{2} \right)^{n}
\quad\mbox{\rm and}\quad
a^2_{0,n}\leq 2\epsilon(m)^{-1}\delta(g)^{3m}~\left(
1-\epsilon(m)^2\delta(g)^{-5m+1}
\right)^{\lfloor n/m\rfloor}
$$
with
some parameter $\epsilon$ such that
$$
 \inf{\left\{
\frac{s}{\mu(1)},\frac{\mu(1)}{s},
\frac{1-s}{1-\mu(1)},\frac{1-\mu(1)}{1-s}
\right\}}\leq \epsilon\leq 1
$$
It is also readily verified that the assumptions of lemma~\ref{lemreff} are satisfied with
the parameters
$$
\begin{array}{rclrrll}
\tau^1&\leq& \delta(g)~\left|s-\mu(1)\right|\qquad
&
\tau^2&\leq&
\delta^{\prime}(g) ~
(({\mu(1)\vee s})/({\mu(1)\wedge s}))\\
c_1&=&2\epsilon^{-1} \delta^{\prime}(g)&
c_2&=&2\epsilon(m)^{-1}\left(
1-\epsilon(m)^2\delta(g)^{-5m+1}
\right)^{-1}\delta(g)^{3m}
\end{array}
$$
and the Lyapunov constants
$$
\lambda_1=-\log{(1-\epsilon^2)}\quad\mbox{\rm and}
\quad
\lambda_2=-\frac{1}{m}~{\log{\left(
1-\epsilon(m)^2\delta(g)^{-5m+1}
\right)}}
$$
We notice that $\epsilon$ tends to $1$ and $\tau^1$ tends to $0$, as $|s-\mu(1)|$ tends to $0$. Thus,
there exists some $\varsigma\geq 0$ such that
$$
\lambda_1>\lambda_2\quad\mbox{\rm and}\quad
c_1c_2\tau^1\tau^2< \left(1-e^{-\lambda_2}\right)~\left(e^{-\lambda_2}-e^{-\lambda_1}\right)
$$
as long as $
|s-\mu(1)|\leq \varsigma$.
We summarize this discussion with the following corollary.
\begin{cor}
Consider the time homogeneous model discussed above.
Under the assumptions of theorem~\ref{theoBernoulli}, for any $N$-approximation models $(\gamma_{n}^{N}(1),\eta _{n}^{N})$ of the Bernoulli model
$(\gamma_n(1),\eta_n)$  satisfying condition
 (\ref{localerr}), the $\LL_r$-mean error estimates presented in (\ref{estimintro}) are uniform w.r.t. the time parameter
$$
\sup_{n\geq 0}{\mathbb{E}\left(\left|V^{\gamma,N}_n(1)\right|^r\right)^{\frac{1}{r}}} \leq
a_r~c^{1,2}/(1-e^{-\lambda})\quad\mbox{and}\quad
 \sup_{n\geq 0}{\mathbb{E}\left(\left|V^{\eta,N}_n(f)\right|^r\right)^{\frac{1}{r}} }\leq
a_r~c^{2,2}/(1-e^{-\lambda})
$$
with the parameters $(c^{1,2},c^{2,2},\lambda)$ defined in lemma~\ref{lemreff}, and some finite constants $a_r<\infty$ whose
values only depend on $r$.
\end{cor}
\begin{rem}
When $\mu_{n+1}(1)=0$ we have seen in (\ref{ex1}) that
$$
 \Phi^1_{n+1,\nu_n}(u)=
 \Psi_{g_n}(\nu_n)(s_n)\times\theta_{\nu_n(g_n)}(u)
 \quad\mbox{\rm and}\quad
 \Phi^2_{n+1,m_n}(\eta)=
 \Psi_{g_ns_n}(\eta)M_{n+1}
 $$
with the collection of mappings $\theta_a$, with $a\in [0,\infty[$, defined in (\ref{theta}). Using the fact that
\begin{eqnarray*}
\left| \Phi^1_{n+1,\nu_n}(u)- \Phi^1_{n+1,\nu_n}(u^{\prime})\right|&=&\frac{\Psi_{g_n}(\nu_n)(s_n)~\nu_n(g_n)}{\left[\nu(g_n)u+(1-u)\right]\left[\nu(g_n)u^{\prime}+(1-u^{\prime})\right]}~|u-u^{\prime}|
\end{eqnarray*}
one proves that (\ref{a1}) is met with the rather crude upper bound
$$
a^1_{p,n}\leq \prod_{p\leq k<n}a_{k,k+1}\quad\mbox{\rm and}\quad
a^1_{k,k+1}\leq {(s^+_kg^+_k)}/(1\wedge g^-_k)^2
$$
We also notice that the second component mapping $\Phi^2_{n+1,m_n}$ doesn't  depends on the parameter $m_n$, and it induces a
 Feynman-Kac semigroup of the same form as the one studied in section~\ref{fksg}. Assuming that
the mixing condition  stated in the l.h.s. of (\ref{Hm}) is satisfied some integer $m\geq 1$ and some parameter $\epsilon_p(m)>0$,
one can prove that (\ref{a2}) is met with the collection of parameters $a^2_{p,n}$ given below
$$
a^2_{p,n}\leq 2~\rho_{p}(m)~\prod_{k=0}^{\lfloor (n-p)/m\rfloor-1}\left(
1-\epsilon^{(m)}_{p+km}
\right)
\quad\mbox{with}
\quad
\rho_p(m)=\epsilon_p^{-1}(m)\prod_{p\leq q<p+m}\delta_{q}(sg)
$$
and the collection of parameters $\epsilon^{(m)}_p=\epsilon^{(m)}_p=\epsilon^2_{p}(m)/\prod_{p<q<p+m}\delta_{q}(sg)$.
\end{rem}
\subsection{PHD Models}
This section is concerned with the contraction properties of the semigroups
$\Phi_{p,n,\nu}^1$ and $\Phi_{p,n,m}^2$ associated with the PHD filter discussed in section~\ref{phdmodel1} and in section~\ref{phdMcKean}.

The analysis of these nonlinear models is much more involved than the one of the Bernoulli models. We simplify the analysis and we further assume that the clutter intensity function, the detectability rate as well as the survival and the spawning
rates introduced in section~\ref{phdMcKean}
are time homogeneous and constants functions, and we set
$$(b_n(x),h_n(x),s_n(x),r_n(x))=(b,h,s,r)
$$
To simplify the presentation, we also assume that the state spaces, the Markov transitions of the targets, the likelihood functions and the spontaneous birth measures
are time homogeneous, that is we have that $E_n=E$, $E^Y_n=E^Y$, $M_{n}=M$,
$g_n(x,y)=g(x,y)$ and $\mu_{n+1}=\mu$. Without further mention, we suppose that $r(1-d)<1$, $\mu(1)>0$, $r>0$, and for any
$y\in E^Y$ we have
$$
0\leq g^-(y):=\inf_{x\in E}g(x,y)\leq g^+(y):=\sup_{x\in E}g(x,y)<\infty
$$
Given a mapping $\theta$ from $E^Y$ into  $\RR$, we set $\Ya^-(\theta):=\inf_n\Ya_n(\theta)$ and  $\Ya^+(\theta):=\sup_n\Ya_n(\theta)$.

We recall from (\ref{mttmodel}) that the PDH filter is defined by the measure-valued equation
$$\gamma_{n+1}=\gamma_nQ_{n+1,\gamma _{n}}$$ with the integral operator
$$
Q_{n+1,\gamma _{n}}(x_{n},dx_{n+1})=g_{n,\gamma
_{n}}(x_{n})M_{n+1}(x_{n},dx_{n+1})+\gamma _{n}(1)^{-1}~\mu _{n+1}(dx_{n+1})
$$
with the function $g_{n,\gamma_n}$ defined below
$$
g_{n,\gamma_n}(x)=r(1-d)+rd~\int \mathcal{Y}_{n}(dy)~
\frac{g(x,y)}{%
h+d\gamma _{n}(g(\mbox{\LARGE .},y))}
$$
We also notice that the total mass process and the normalized distribution flow are given by the following equations
\begin{eqnarray*}
\gamma_{n+1}(1)&=&\Phi^1_{n+1,\eta_n}(\gamma_n(1))\\
&=&\gamma_{n}(1)~r(1-d)+\int \mathcal{Y}_{n}(dy)~
w_{\gamma_n(1)}(\eta_n,y)+\mu(1)\\
\eta_{n+1}(1)&=&\Phi^2_{n+1,\gamma_n(1)}(\eta_n)\\
&\propto& \gamma_n(1)~r(1-d)~\eta_n M+\int \Ya_n(dy)~w_{\gamma_n(1)}(\eta_n,y)~\Psi_{g(\point,y)}(\eta_n)M + \mu(1)~\overline{\mu}
\end{eqnarray*}
with the probability measure $\overline{\mu}$ and weight functions $w$ defined below
$$
\overline{\mu}(dx)=\mu(dx)/\mu(1)\quad\mbox{\rm and}\quad w_u(\eta,y):=r\left(1-\frac{h}{h+du\eta(g(\point,y))}\right)
$$

For null clutter parameter $h=0$, we already observe that the total mass transformation $\Phi^1_{n+1,\eta_n}$
doesn't depend on the flow of probability measures $\eta_n$ and it
is simply given by
$$
\Phi^1_{n+1,\eta_n}(\gamma_n(1))=
\gamma_{n}(1)~r(1-d)+r~\mathcal{Y}_{n}(1)+\mu(1)
$$
In this particular situation, we have $$
\gamma^N_n(1)=\gamma_n(1)=(r(1-d))^{n}\gamma_0(1)+\sum_{0\leq k<n}(r(1-d))^{n-1-k} (r~\mathcal{Y}_{k}(1)+\mu(1)
)
$$ Now, we easily show that the pair of conditions (\ref{a1}) and (\ref{tau1}) are satisfied
with the parameters $a^1_{p,n}=(r(1-d))^{n-p}$ and $\tau^1_n=0$. In more general situations, the total mass process is not explicitly known. Some useful estimates are provided by the following lemma.

\begin{lem}
We assume that the number of observations is uniformly bounded; that is, we have that $\Ya^+(1)<\infty$. In this situation,
the total mass process $\gamma_n(1)$ and any approximation model
$\gamma^N_n(1)$ given by the recursion (\ref{massapprox}) (with the initial condition $\gamma^N_0(1)=\gamma_0(1)$) take values in a sequence of compact sets $I_n\subset [m^-,m^+]$
with
$$
m^-:=\frac{\mu(1)}{1-r(1-d)}~\left(1+rd~\Ya^-\left(\frac{g^-}{h+d\mu(1)g^-}\right)\right)
\quad\mbox{
and}
\quad
m^+:= \gamma_0(1)+\frac{r\mathcal{Y}^+(1)+\mu(1)}{1-r(1-d)}
$$
\end{lem}
\proof
Using the fact that $\gamma_n(1)\geq \mu(1)$ we prove that
$$
r\left(1-\frac{h}{h+d\mu(1)~g^-(y)}\right)
\leq w_{\gamma_n(1)}(\eta_n,y)
\leq r
$$
from which we conclude that
$$
\gamma_{n}(1)~r(1-d)+r~\mathcal{Y}_{h,n}(1)+\mu(1)\leq
\Phi^1_{n+1,\eta_n}(\gamma_n(1))\leq \gamma_{n}(1)~r(1-d)+r~\mathcal{Y}_{n}(1)+\mu(1)
$$
with the random measures
$$
\mathcal{Y}_{h,n}(dy):=\mathcal{Y}_{n}(dy)~\frac{d\mu(1)~g^-(y)}{h+d\mu(1)~g^-(y)}
$$
For any sequence of probability
measures $\nu:=(\nu_n)_{n\geq 0}\in\Pa(E)^{\NN}$, and any starting mass $u\in [0,\infty[$ one conclude that
$$
\left(r(1-d)\right)^n u+\frac{r\mathcal{Y}_{h}^-(1)+\mu(1)}{1-r(1-d)}
\leq \Phi_{0,n,\nu}^1(u)\leq
\left(r(1-d)\right)^n u+\frac{r\mathcal{Y}^+(1)+\mu(1)}{1-r(1-d)}
$$
 This implies that  $\gamma_n(1),\gamma^N_n(1)\in I_n\subset [m^-,m^+]$
with
$$
m^-:=\frac{r\mathcal{Y}_{h}^-(1)+\mu(1)}{1-r(1-d)}=\frac{\mu(1)}{1-r(1-d)}~\left(1+rd~\Ya^-\left(\frac{g^-}{h+d\mu(1)g^-}\right)\right)
$$
The end of the proof of the lemma is now completed.
\cqfd
We are now in position to state the main result of this section.
\begin{theo}\label{theophdunif}
We assume that the number of observations is uniformly bounded; that is, we have that $\Ya^+(1)<\infty$. In this situation,  the condition $(\mbox{Lip$(\Phi)$})$ is met with the Lipschitz constants $a^i_{p,n}\leq \prod_{p\leq k<n}a^i_{k,k+1}$, with $i=1,2$, and the sequence of parameters $\left(a^i_{n,n+1}\right)_{n\geq 0}$, $i=1,2$, defined below
$$
a^1_{n,n+1}\leq
r(1-d)+rd h~\mathcal{Y}_{n}\left(
\frac{g^+}{%
[h+dm^- g^-]^2}\right)
$$
and
$$
a^2_{n,n+1}\leq m^+
\frac{\beta(M)\left[
(1-d)+d~\Ya_n\left(\frac{g^+}{h+dm^+g^+}~\frac{g^+}{g^-}\right)\right]+hd\Ya_n\left(\frac{g^+-g^-}{(h+dm^-g^-)^2}\right)}{(1-d)~m^-+dm^-\Ya_n\left(\frac{g^-}{h+dm^-g^-}\right)+\mu(1)/r}
$$
In addition, condition $(\mbox{Cont$(\Phi)$})$ is met with the sequence of parameters
$$
\tau^1_{n+1}\leq rdh m^+~
\mathcal{Y}_{n}\left(
\frac{g^+-g^-}{%
[h+dm^-g^-]^2}\right)
\quad
\tau^{2}_{n+1}\leq \frac{
(1-d)+hd~
\Ya_n\left(\frac{g^+}{(h+dm^-g^-)^2}\right)}{(1-d)~m^-+dm^-\Ya_n\left(\frac{g^-}{h+dm^-g^-}\right)+\mu(1)/r}
$$
\end{theo}
The proof of theorem~\ref{theophdunif} is postponed to section~\ref{ptheophdunif}.
\begin{cor}
We assume that $\Ya^+\left({g^+}/{g^-}\right)$ and $\Ya^+\left({g^+}/{(g^-)^2}\right)<\infty$. In this situation, there exists some parameters $0<\kappa_0\leq 1$, $\kappa_1<\infty$, and $\kappa_2>0$
such that for any $d\geq \kappa_0$, $\mu(1)\geq \kappa_1$, and $h\leq \kappa_2$,
 the semigroups  $\Phi^1_{p,n,\nu}$ and $\Phi^2_{p,n,m}$ satisfy the  pair of conditions $(\mbox{Lip$(\Phi)$})$ and $(\mbox{Cont$(\Phi)$})$ with some parameters $(a^i_{p,n},\tau^i_n)_{i=1,2,p\leq n}$, satisfying the assumptions of lemma~\ref{lemreff}. In particular, for any $N$-approximation models $(\gamma_{n}^{N}(1),\eta _{n}^{N})$ of the PHD equation
$(\gamma_n(1),\eta_n)$  satisfying condition
 (\ref{localerr}), the $\LL_r$-mean error estimates presented in (\ref{estimintro}) are uniform w.r.t. the time parameter
 $$
 \sup_{n\geq 0}{\mathbb{E}\left(\left|V^{\gamma,N}_n(1)\right|^r\right)^{\frac{1}{r}}} \leq
a_r~c^{1,2}/(1-e^{-\lambda})\quad\mbox{and}\quad
 \sup_{n\geq 0}{\mathbb{E}\left(\left|V^{\eta,N}_n(f)\right|^r\right)^{\frac{1}{r}} }\leq
a_r~c^{2,2}/(1-e^{-\lambda})
 $$
with the parameters $(c^{1,2},c^{2,2},\lambda)$ defined in lemma~\ref{lemreff}, and some finite constants $a_r<\infty$ whose
values only depend on $r$.
\end{cor}
\proof
There is no loss of generality to assume that $r(1-d)<1/2\leq d$ and $\mu(1)\geq 1\geq h$. Recalling that $
m^-\geq \mu(1)
$, one readily  proves that
$$
\frac{m^+}{\mu(1)}=\frac{\gamma_0(1)}{\mu(1)}+\frac{1}{1-r(1-d)}
\left(1+\frac{r}{\mu(1)}~\mathcal{Y}^+(1)\right)\leq 2+
\gamma_0(1)+2r\mathcal{Y}^+(1):=\rho
$$
If we set
$
\delta(g):=\rho\vee\Ya^+\left(\frac{g^+}{g^-}\right)\vee\Ya^+\left(\frac{g^+}{(g^-)^2}\right)
$,
then we find the rather crude estimates
$$
a^1_{n,n+1}/r\leq
(1-d)+\frac{2h}{\mu(1)^2}~\delta(g)
\quad
\mbox{\rm and}
\quad
a^2_{n,n+1}/r
\leq \left[\beta(M)
(1-d)+\frac{2h+\beta(M)}{\mu(1)}~\right]\delta(g)
$$
as well as
$$
\tau^1_{n+1}/r\leq \frac{2h}{\mu(1)}
~\delta(g)^2
\quad
\mbox{\rm and}
\quad
\tau^{2}_{n+1}/r\leq \frac{1}{\mu(1)}\left[
(1-d)+\frac{2h}{\mu(1)^2}~
\delta(g)\right]
$$
from which we find that
\begin{equation}\label{tauphd}
\tau^1\tau^2\leq
\frac{2hr^2}{\mu(1)^2}\left[
(1-d)+\frac{2h}{\mu(1)^2}~
\delta(g)\right]~\delta(g)^2
\end{equation}
Thus, there exists some $0<\kappa_0\leq 1$ and some $\kappa_1<\infty$ so that
for any $d\geq \kappa_0$ and any $\mu(1)\geq \kappa_1$ we have
\begin{eqnarray*}
a^1_{n,n+1}&\leq& r\left[
(1-d)+\frac{2}{\mu(1)^2}\right]~\delta(g):=e^{-\lambda_1}<1\\
a^2_{n,n+1}
&\leq &r \left[
(1-d)+\frac{3}{\mu(1)}~\right]\delta(g):=e^{-\lambda_2}<1\quad\mbox{\rm with}\quad
0<\lambda_2<\lambda_1
\end{eqnarray*}
Finally, using (\ref{tauphd}) we find some $\kappa_2>0$ such that for any $h\leq \kappa_2$, we have that
$
\tau^1\tau^2\leq  \left(1-e^{-\lambda_2}\right) \left(e^{-\lambda_2}-e^{-\lambda_1}\right)
$.
The end of the proof is now a direct consequence of lemma~\ref{lemreff}. This ends the proof of the corollary.
\cqfd
\section{Stochastic particle approximations}\label{secalgo}

\subsection{Mean field interacting particle systems}\label{secmfield}

\subsubsection{Description of the models}

The mean field type interacting particle system associated with the equation
(\ref{flotPhi}) relies on the fact that the one step mappings $%
\Gamma^2_{n+1} $ can be rewritten in the following form
\begin{equation}  \label{flotPhimc}
\Gamma_{n+1}^2(\gamma_{n}(1),\eta_{n})=\eta_{n}K_{n+1,\gamma_{n}}\quad%
\mbox{\rm with}\quad \gamma_n=\gamma_n(1)\times\eta_n
\end{equation}
for some collection of Markov kernels $K_{n+1,\gamma}$ indexed by the time
parameter $n$ and the set of measures $\gamma\in \mathcal{M}_{+}(E_n)$. We
mention that the choice of the Markov transitions $K_{n,\gamma}$ is not
unique. In the literature on mean field particle models, $K_{n,\gamma}$ are
called a choice of McKean transitions. Some McKean interpretation models of the
Bernoulli and the PHD
filter models (\ref{Bernoullimodel}
) and  (\ref{mttmodel}) are discussed in section~\ref{phdMcKean}
(see for instance (\ref{phdmck})) and in section~\ref{bernoulliMcKean} (see for instance
~\ref{bernoullimck})

These models provide a natural interpretation of the distribution laws $%
\eta_n$ as the laws of a non linear Markov chain $\overline{X}_n$ whose
elementary transitions $\overline{X}_n\leadsto \overline{X}_{n+1}$ depends
on the distribution $\eta_n=\mbox{\rm Law}(\overline{X}_n)$, as well as on
the current mass process $\gamma_n(1)$. In contrast to traditional McKean
model, the dependency on the mass process induce a dependency of all the
flow of measures $\eta_p$, for $0\leq p\leq n$. For a thorough description
of these discrete generation and non linear McKean type models, we refer the
reader to~\cite{fk}.

In further developments of the article, we always assume that the
mappings
\begin{equation*}
\left(m,x_n,\left(x^i\right)_{1\leq i\leq N}\right)\mapsto
K_{n+1,m\sum_{j=1}^N\delta_{x^j}}\left(x_n,A_{n+1}\right) ~\mbox{\rm
and} ~ G_{n+1,m\sum_{j=1}^N\delta_{x^j}}(x_n)
\end{equation*}
are pointwise known, and of course measurable w.r.t. the corresponding
product sigma fields, for any $n\geq 0$, $N\geq 1$, $A_{n+1}\in \mathcal{E }%
_{n+1}$, and any $x_n\in E_n$. In this situation, the mean field particle
interpretation of this nonlinear measure-valued model is an $E^N_n$-valued
Markov chain
\index{$\xi^{(N)}_n$} $\xi^{(N)}_n=\left(\xi^{(N,i)}_n\right)_{1\leq i\leq
N} $, with elementary transitions defined as
\begin{eqnarray}
\gamma_{n+1}^N(1)&=&\gamma_n^N(1)~\eta_{n}^N(G_{n,\gamma^N_n})
\label{meanfieldmass} \\
\mathbb{P}\left(\xi^{(N)}_{n+1}\in dx ~\left|~ \mathcal{F }%
^{(N)}_n\right.\right) &=&\prod_{i=1}^N~K_{n+1,\gamma^N_n}(
\xi^{(N,i)}_n,dx^i)  \label{meanfieldeta}
\end{eqnarray}
with the pair of occupation measures $\left(\gamma^N_n,\eta^N_n\right)$
defined below
\begin{equation*}
\eta^N_n:=%
\frac{1}{N} \sum_{j=1}^N\delta_{\xi_n^{(N,j)}} \quad\mbox{\rm and}\quad
\gamma_n^N(dx):=\gamma_n^N(1)~\eta^N_n(dx)
\end{equation*}
In the above displayed formula, $\mathcal{F }^{N}_n$ stands for the $\sigma$%
-field generated by the random sequence $(\xi^{(N)}_p)_{0\leq p\leq n}$, and
$dx=dx^1\times\ldots\times dx^N$ stands for an infinitesimal neighborhood of
a point $x=(x^1,\ldots,x^N)\in E_n^N$. The initial system $\xi^{(N)}_0$
consists of $N$ independent and identically distributed random variables
with common law $\eta_0$. As usual, to simplify the presentation, when there
is no possible confusion we suppress the parameter $N$, so that we write $%
\xi_n$ and $\xi^i_n$ instead of $\xi^{(N)}_n$ and $\xi^{(N,i)}_n$.

\subsubsection{Convergence analysis}

The rationale behind the mean field particle model described in (\ref%
{meanfieldeta}) is that $\eta^N_{n+1}$ is the empirical measure associated
with $N$ independent variables with distributions $K_{n+1,\gamma^N_n}\left(%
\xi^{i}_n,dx\right)$, so as long as $\gamma^N_{n}$ is a good approximation
of $\gamma_{n}$ then $\eta^N_{n+1}$ should be a good approximation of $%
\eta_{n+1}$. Roughly speaking, this induction argument shows that $\eta^N_n$
tends to $\eta_n$, as the population size $N$ tends to infinity.

These stochastic particle algorithms can be thought of in various ways: From
the physical view point, they can be seen as microscopic particle
interpretations of physical nonlinear measure-valued equations. From the
pure mathematical point of view, they can also be interpreted as natural
stochastic linearizations of nonlinear evolution semigroups. From the
probabilistic point of view, they can be interpreted as a interacting
recycling acceptance-rejection sampling techniques. In this case, they can
be seen as a sequential and interacting importance sampling technique.

By construction, the local fluctuation random fields $(W^N_n)_{n\geq 0}$
defined in (\ref{local}) can be rewritten as follows
\begin{equation*}
\eta^N_n=\eta^N_{n-1}K_{n,\gamma^N_{n-1}}+\frac{1}{\sqrt{N}}~W^N_n
\end{equation*}
Using Khintchine's inequality, we can check that (\ref{localerr}) is met for
any $r\geq 1$ and any $f_n\in \mbox{Osc}_{1}(E_n)$, with the collection of
universal constants given below
\begin{equation*}
a_{2r}^{2r}\leq (2r)!~2^{-r}/r!\quad \mbox{and}\quad a_{2r+1}^{2r+1}\leq
(2r+1)!~2^{-r}/r!
\end{equation*}

We end this section with a brief discussion on the PHD equation presented in
(\ref{mttmodel}). This model combines in a single step the traditional
updating and a prediction filtering transition. This combination allows us to
reduce the fluctuations of the local sampling errors and their propagations
w.r.t. the time parameter. Since these updating-prediction models are
often used in the literature of multiple target tracking, we provide below a
short summary. If we set
\begin{equation*}
\widehat{g}^c_{n,\gamma}(\mbox{\LARGE .},y)=\left\{
\begin{array}{lcl}
(1-d_n) & \mbox{\rm if} & y=c \\
\displaystyle\frac{d_ng_n(\mbox{\LARGE .},y_n)}{h_n(y)+\gamma(d_ng_n(%
\mbox{\LARGE .},y))} & \mbox{\rm if} & y\not=c%
\end{array}%
\right.
\end{equation*}
then
\begin{equation*}
\gamma_{n+1}=\widehat{\gamma}_n Q_{n+1}+\mu_{n+1}\quad\mbox{\rm with}\quad
Q_{n+1}(f):=r_n~M_{n+1}(f)
\end{equation*}
with the updated measures defined below
\begin{equation*}
\widehat{\gamma}_n(f):=\gamma_n(\widetilde{g}^c_{n,\gamma_n}f) \quad%
\mbox{\rm with}\quad \widetilde{g}^c_{n,\gamma_n}=\int \mathcal{Y }_n^c(dy)~%
\widehat{g}^c_{n,\gamma_n}(\mbox{\LARGE .},y)
\end{equation*}
Notice that
\begin{equation*}
\widehat{\gamma}_n(1)=\gamma_n(\widetilde{g}^c_{n,\gamma_n}f)\quad%
\mbox{\rm
and}\quad \widehat{\eta}_n(dx):={\widehat{\gamma}_n(dx)}/{\widehat{\gamma}%
_n(1)}=\Psi_{\widetilde{g}^c_{\gamma_n,n}}(\eta_n)(dx)
\end{equation*}
from which we find the recursive formulae
\begin{equation*}
\left(
\begin{array}{c}
\gamma_n(1) \\
\eta_n%
\end{array}
\right)\overset{\mbox{\small updating}}{-\!\!\!-\!\!\!-\!\!\!-\!\!\!-\!\!\!-%
\!\!\!-\!\!\!-\!\!\!-\!\!\!-\!\!\!-\!\!\!-\!\!\!\longrightarrow} \left(
\begin{array}{c}
\widehat{\gamma}_n(1) \\
\widehat{\eta}_n%
\end{array}
\right) \overset{\mbox{\small prediction}}{-\!\!\!-\!\!\!-\!\!\!-\!\!\!-\!\!%
\!-\!\!\!-\!\!\!-\!\!\!-\!\!\!-\!\!\!-\!\!\!-\!\!\!\longrightarrow} \left(
\begin{array}{c}
\gamma_{n+1}(1) \\
\eta_{n+1}%
\end{array}
\right)
\end{equation*}
with the prediction transition described below
\begin{equation*}
\gamma_{n+1}(1)=\widehat{\gamma}_n(r_n)+\mu_{n+1}(1) \quad \mbox{\rm and}
\quad \eta_{n+1}=\Psi_{r_n}\left(\widehat{\eta}_n\right)M_{n+1,\widehat{%
\gamma}_n}^{\prime}
\end{equation*}
In the above displayed formula, $M_{n+1,\widehat{\gamma}_n}^{\prime}$ is the
Markov transition defined by
\begin{equation*}
M_{n+1,\widehat{\gamma}_n}^{\prime}(x,\mbox{\LARGE .})=\alpha^{\prime}_n(%
\widehat{\gamma}_n)~M_{n+1}(x,\mbox{\LARGE .})+ \left(1-\alpha^{\prime}_n(%
\widehat{\gamma}_n)\right)~\overline{\mu}_{n+1}
\end{equation*}
with the collection of $[0,1]$-valued parameters $\alpha^{\prime}_n(\widehat{%
\gamma}_n)={\widehat{\gamma}_n(r_n)}/{\left(\widehat{\gamma}%
_n(r_n)+\mu_{n+1}(1)\right)}$. It should be clear that the updating and
the prediction transitions can be approximated using a genetic type
selection and mutation transition. Each of these sampling transitions
introduces a separate local sampling fluctuation error. The stochastic
analysis of the corresponding mean field particle interpretations can be
developed using the same line of arguments as those used for the particle
model discussed above.

\subsection{Interacting particle association systems}\label{secipsasso}

\subsubsection{Description of the models}

We let $( \mathcal{A }_{n})_{n\geq 0}$ be a sequence of finite sets equipped
with some finite positive measures $(\nu_n)_{n\geq 0}$. We further assume
that the initial distribution $\gamma_0$ and the integral operators $%
Q_{n+1,\gamma_n}$ in (\ref{defmod}) have the following form
\begin{equation*}
\gamma_0=\int \nu_{0}(da)~\eta^{(a)}_0 \quad \mbox{\rm and} \quad
Q_{n+1,\gamma_n}=\int \nu_{n+1}(da)~Q^{(a)}_{n+1,\gamma_n}
\end{equation*}
In the above display $\eta^{(a)}_0$ stands for a collection of measures on $%
E_0$, indexed by the parameter $a\in \mathcal{A }_0$, and $%
Q^{(a)}_{n+1,\gamma_n}$ is a collection of integral operators indexed by the
parameter $a\in \mathcal{A }_{n+1}$. In this situation, we observe that
\begin{equation*}
\gamma_0(1)=\nu_0(1)\quad\mbox{\rm and}\quad\eta_0=\int
A_0(da)~\eta^{(a)}_0\quad\mbox{\rm with}\quad A_0(da):={\nu_{0}(da)}/{%
\nu_{0}(1)}
\end{equation*}

We also assume that the following property is met
\begin{equation}  \label{defGM}
G^{(a)}_{n,\gamma}:=Q^{(a)}_{n+1,\gamma}(1)\propto G^{(a)}_n \quad%
\mbox{\rm
and}\quad {Q^{(a)}_{n+1,\gamma}(f)}/{Q^{(a)}_{n+1,\gamma}(1)}:=
M^{(a)}_{n+1}(f)
\end{equation}
for some function $G^{(a)}_n$ on $E_n$, and some Markov transitions $%
M^{(a)}_{n+1}$ from $E_n$ into $E_{n+1}$ whose values do not depend on the
measures $\gamma$. For clarity of presentation, sometimes we write $%
\Psi^{(a)}_{G_{n}}$ instead of $\Psi_{G_{n}^{(a)}}$.
\begin{defi}
We consider the collection of probability measures $\eta_n^{(\mathrm{a}%
_n)}\in \mathcal{P }(E_n)$, indexed by sequences of parameters
\begin{equation*}
\mathrm{a}_n=(a_0,\ldots,a_n)\in \mathcal{A }_{[0,n]}:=\left( \mathcal{A }%
_0\times\ldots\times \mathcal{A }_n\right)
\end{equation*}
and defined by the following equations
\begin{eqnarray}
\eta_{n}^{(\mathrm{a}_{n})}&=&\left(\Phi_{n}^{(a_{n})}\circ\ldots\circ%
\Phi_{1}^{(a_{1})}\right)\left(\eta^{(a_0)}_0\right)  \label{reca}
\end{eqnarray}
with the mappings $\Phi_{n}^{(a)}: \mathcal{P }(E_{n-1})\rightarrow \mathcal{%
P }(E_n)$ indexed by $a\in \mathcal{A }_{n}$ and defined by the updating-prediction
transformation
\begin{equation*}
\Phi_{n}^{(a)}\left(\eta\right)=\Psi^{(a)}_{G_{n-1}}\left(\eta%
\right)M^{(a)}_{n}
\end{equation*}
\end{defi}
We illustrate these abstract conditions in the context of the multiple
target tracking equation presented in (\ref{mttmodel}). In this situation,
it is convenient to add a pair of virtual observation states $c,c^{\prime}$
to $E_{n}^Y$. Using this notation, the above conditions are satisfied with
the finite sets $\mathcal{A }_{n+1}$ and their counting measures $\nu_{n+1}$
defined below
\begin{equation*}
\mathcal{A }_{n+1}=\left\{Y^i_n, 1\leq i\leq
N^Y_n\}\cup\{c,c^{\prime}\right\}\qquad \nu_{n+1}= \mathcal{Y }%
_{n}+\delta_{c}+\delta_{c^{\prime}}\in \mathcal{M }( \mathcal{A }_{n+1})
\end{equation*}
Using (\ref{mttmodel}) and (\ref{decompo}), we check that (\ref{defGM}) is
met with the couple of potential functions and Markov transitions defined by
\begin{equation*}
(G^{(y_n)}_{n},M^{(y_n)}_{n+1})=\left\{
\begin{array}{lcl}
\left(r_{n}d_{n}g_{n}(\mbox{\LARGE .},y_n),M_{n+1}\right) &
\mbox{\rm
for} & y_n\not\in\{c,c^{\prime}\} \\
\left(r_{n}(1-d_{n}),M_{n+1}\right) & \mbox{\rm for} & y_n=c \\
\left(1,\overline{\mu}_{n+1}\right) & \mbox{\rm for} & y_n=c^{\prime}%
\end{array}
\right.
\end{equation*}
In this case, we observe that
\begin{equation*}
Q^{(y_n)}_{n+1,\gamma_n}(x_{n},\mbox{\LARGE .})=G^{(y_n)}_{n,%
\gamma_n}(x_n)~M^{(y_n)}_{n+1}(x_{n},\mbox{\LARGE .})
\end{equation*}
with the potential function $G^{(y_n)}_{n,\gamma_n}$ defined below
\begin{equation}  \label{phdex}
G^{(y_n)}_{n,\gamma_n}/G^{(y_n)}_{n}=\left\{
\begin{array}{ccl}
\left[{h_n(y_n)+\gamma_n(d_ng_n(\mbox{\LARGE .},y_n))}\right]^{-1}~ & %
\mbox{\rm for} & y_n\not\in\{c,c^{\prime}\} \\
1 & \mbox{\rm for} & y_n=c \\
{\mu_{n+1}(1)}/{\gamma_n(1)} & \mbox{\rm for} & y_n=c^{\prime}%
\end{array}
\right.
\end{equation}

Under our assumptions, using (\ref{flotPhi}), we have the following result.

\begin{prop}
The solution the equation (\ref{flotPhi}) has the following form
\begin{equation*}
\eta_n=\int A_n(da)~\eta_n^{(a)}
\end{equation*}
with a total mass process $\gamma_n(1)$ and the association measures $A_n\in
\mathcal{P }( \mathcal{A }_{[0,n]})$ defined by the following recursive
equations
\begin{equation*}
\gamma_{n+1}(1)=\gamma_n(1)~\eta_n(G_{n,\gamma_n})\quad\mbox{\rm and}\quad
A_{n+1}=\Omega_{n+1}\left(\gamma_n(1),A_n\right)
\end{equation*}
With the mapping
\begin{equation*}
\Omega_{n+1}~:~(m,A)\in \left(]0,\infty[\times \mathcal{P }( \mathcal{A }%
_{[0,n]})\right)\mapsto \Omega_{n+1}(m,A)\in \mathcal{P }( \mathcal{A }%
_{[0,n+1]})
\end{equation*}
defined by the following formula
\begin{equation}  \label{defOmega}
\begin{array}{l}
\Omega_{n+1}\left(m,A\right)(d(a,b)) \propto
A(da)~\nu_{n+1}(db)~\eta_n^{(a)}\left( G^{(b)}_{n,m \int A(da)~\eta^{(a)}_n}
\right)%
\end{array}%
\end{equation}
\end{prop}

\noindent\mbox{\bf Proof:}\newline
The proof of the above assertion is simply based on the fact that
\begin{eqnarray*}
\eta_{n+1}\propto \int \nu_{n+1}(db)~\eta_nQ^{(b)}_{n+1,\gamma_n} &=&\int
A_n(da)~\nu_{n+1}(db)~\eta_n^{(a)}Q^{(b)}_{n+1,\gamma_n} \\
&=& \int A_n(da)~\nu_{n+1}(db)~\eta_n^{(a)}\left( G^{(b)}_{n,\gamma_n}
\right)~\eta_{n+1}^{(a,b)}
\end{eqnarray*}
This clearly implies that
\begin{equation*}
\Gamma^2_{n}\left(m, \int A(da)~\eta^{(a)}_{n-1}) \right)=\int
\Omega_{n}\left(m,A\right)(d(a,b))~\eta^{(a,b)}_n
\end{equation*}
This ends the proof of the proposition.
\hfill\hbox{\vrule height 5pt width 5pt depth 0pt}\medskip \newline

By construction, we notice that for any discrete measure $A\in \mathcal{P }(
\mathcal{A }_{[0,n-1]})$, and any collection of measures $\eta^{(a)}\in
\mathcal{P }(E_{n-1})$, with $a\in \mathcal{A }_{[0,n-1]}$ we have the
formula
\begin{equation*}
\Gamma^2_{n}\left(m, \int A(da)~\eta^{(a)}) \right)=\int
\Omega_{n}\left(m,A\right)(d(a,b))~~\Phi_{n}^{(b)}\left(\eta^{(a)}\right)
\end{equation*}

\subsubsection{Particle approximation models}

To get some feasible solution, we further assume that $\eta_n^{(a)}\left(
G^{(b)}_{n,\gamma_n}\right)$ are explicitly known for any sequence of
parameters $(a,b) \in\left( \mathcal{A }_{[0,n]}\times \mathcal{A }%
_{n+1}\right)$. This rather strong condition is satisfied for the multiple
target tracking model discussed above as long as the quantities
\begin{equation*}
\eta_n^{(a_0,y_0,\ldots,y_{n-1})}(r_{n}d_{n}g_{n}(\mbox{\LARGE .}%
,y_n))\quad\eta_n^{(a_0,y_0,\ldots,y_{n-1})}\left(r_{n}(1-d_{n})\right)%
\quad \eta_n^{(a_0,y_0,\ldots,y_{n-1})}(d_n g_n(\mbox{\LARGE .},y_n))
\end{equation*}
are explicitly known. This condition is clearly met for linear gaussian
target evolution and observation sensors as long as the survival and
detection probabilities $s_n$ and $d_n$ are state independent, and
spontaneous birth $\overline{\mu}_n$ and spawned targets branching rates $%
b_n $ are Gaussian mixtures. In this situation, the collection of measures $%
\eta_n^{(a_0,y_0,\ldots,y_{n-1})}$ are gaussian distributions and the
equation (\ref{reca}) coincides with the traditional updating-prediction
transitions of the discrete generation Kalman-Bucy filter.

We let $A^N_0=\frac{1}{N}\sum_{i=1}^N\delta_{\mathrm{a}^i_0}$, be the
empirical measure associated with $N$ independent and identically
distributed random variables $(\mathrm{a}^i_0)_{1\leq i\leq N}$ with common
distribution $A_0$. By construction, we have
\begin{equation*}
\eta^N_0:=\int A^N_0(da)~\eta^{(a)}_0=\eta_0+\frac{1}{\sqrt{N}}~W^N_0
\end{equation*}
with some local sampling random fields satisfying (\ref{localerr}). We
further assume that $\gamma_0(1)$ is known and we set $\gamma_0^N=%
\gamma_0(1)~\eta^N_0$.
\begin{equation*}
\gamma_1^N(1)=\gamma_0^N(1)~\eta^N_0(G_{0,\gamma^N_0})\quad\mbox{\rm and}%
\quad \eta^N_1:=\int A^N_1(da)~\eta^{(a)}_1
\end{equation*}
with the occupation measure $A^N_1=\frac{1}{N}\sum_{i=1}^N\delta_{\mathrm{a}%
^i_{1}}$ associated with $N$ conditionally independent and identically
distributed random variables $\mathrm{a}^i_1:=(a^i_{0,1},a^i_{1,1})$ with
common law $\Omega_{1}\left(\gamma_0^N(1),A_0^N\right)$. By construction, we
also have
\begin{equation*}
\eta^N_1:=\int \Omega_{1}\left(\gamma_0^N(1),A_0^N\right)(da)~\eta^{(a)}_1 +%
\frac{1}{\sqrt{N}}~W^N_1=\Gamma^2_{1}\left(\gamma_{0}^N(1),\eta^N_{0}\right)+%
\frac{1}{\sqrt{N}}~W^N_1
\end{equation*}
with some local sampling random fields satisfying (\ref{localerr}).
Iterating this procedure, we define by induction a sequence of $N$-particle
approximation measures
\begin{equation*}
\gamma_n^N(1)=\gamma_{n-1}^N(1)~\eta^N_{n-1}(G_{n-1,\gamma^N_{n-1}})\quad%
\mbox{\rm and}\quad \eta^N_n:=\int A^N_n(da)~\eta^{(a)}_n
\end{equation*}
with the occupation measure $A^N_n=\frac{1}{N}\sum_{i=1}^N\delta_{\mathrm{a}%
^i_{n}}$ associated with $N$ conditionally independent and identically
distributed random variables $\mathrm{a}^i_n:=(a^i_{0,n},a^i_{1,n},%
\ldots,a_{n,n}^i)$ with common law $\Omega_{n}\left(%
\gamma_{n-1}^N(1),A_{n-1}^N\right)$. Arguing as above, we find that
\begin{equation*}
\eta^N_n=\int
\Omega_{n}\left(\gamma_{n-1}^N(1),A_{n-1}^N\right)(da)~\eta^{(a)}_n+\frac{1}{%
\sqrt{N}}~W^N_n =\Gamma^2_{n}\left(\gamma_{n-1}^N(1),\eta^N_{n-1}\right)+%
\frac{1}{\sqrt{N}}~W^N_n
\end{equation*}
with some local sampling random fields satisfying (\ref{localerr}).

\subsubsection{Convergence analysis}

The main objective of this section is to show that $N$-particle occupation
measures $A^N_n$ converge in a sense to be given, as $N$ tends to $\infty$,
to the association probability measures $A_n$. To this end we observe that
the one step mapping $\Omega_{n+1}$ introduced in (\ref{defOmega}) can be
rewritten in the following form
\begin{equation*}
\Omega_{n+1}\left(m,A\right)(F)=\frac{A \mathcal{Q }_{n+1,mA}(F)}{A \mathcal{%
Q }_{n+1,mA}(1)}
\end{equation*}
with the collection of integral operators $\mathcal{Q }_{n+1,m A}$ from $%
\mathcal{A }_{[0,n]}$ into $\mathcal{A }_{[0,n+1]}$ defined below
\begin{equation*}
\mathcal{Q }_{n+1,B}(a,d(a^{\prime},b)):= \delta_{a}(da^{\prime})~
\nu_{n+1}(db)~\eta_n^{ (a^{\prime}) }\left( \mathcal{G }_{n,B}^{(b)} \right)
\quad\mbox{\rm where}\quad \mathcal{G }_{n,B}^{(b)}:=G^{(b)}_{n,\int
B(da)~\eta^{(a)}_n}
\end{equation*}
with $B=mA$. In the above display $d(a^{\prime},b)=da^{\prime}\times db$
stands for an infinitesimal neighborhood of the point $(a^{\prime},b)\in
\mathcal{A }_{[0,n+1]}$, with $a=(a_0^{\prime},\ldots,a_n^{\prime})\in
\mathcal{A }_{[0,n]}$ and $b\in \mathcal{A }_{n+1}$, and $%
a=(a_0,\ldots,a_n)\in \mathcal{A }_{[0,n]}$. It is important to point out
that
\begin{equation*}
B_n:=\gamma_n(1)\times A_n\Longrightarrow B_{n+1}=B_n \mathcal{Q }_{n+1,B_n}
\end{equation*}
Notice that the flow of measures $(B_n)_{n\geq 0}$ satisfies the same type
of equation as in (\ref{defmod}), with the a total mass evolution of the
same form as (\ref{massevol}):
\begin{equation*}
B_{n+1}(1)=B_n(1)~A_n\left( \mathcal{G }_{n,B_n}\right)\quad\mbox{\rm with}%
\quad \mathcal{G }_{n,mA}:=\int \nu_{n+1}(db)~ \mathcal{G }_{n,mA}^{(b)}
\end{equation*}

\begin{equation*}
\mathcal{Q }_{n+1,B_n}(F)(a)=\int \nu_{n+1}(db)~\eta_n^{ (a) }\left(
\mathcal{G }_{n,B_n}^{(b)} \right)~F(a,b)
\end{equation*}

\begin{equation*}
\left[ \mathcal{Q }_{n+1,B}(F)- \mathcal{Q }_{n+1,B^{\prime}}(F)\right]%
(a)=\int \nu_{n+1}(db)~ ~\left[ \eta_n^{ (a) }\left( \mathcal{G }%
_{n,B}^{(b)} \right)- \eta_n^{ (a) }\left( \mathcal{G }_{n,B^{\prime}}^{(b)}
\right) \right] ~F(a,b)
\end{equation*}
If we set $B=mA$ and $B^{\prime}=m^{\prime}A^{\prime}$ then condition $%
(H^{\prime}_2)$ is met as long as
\begin{equation*}
\left| \eta_n^{ (a) }\left( \mathcal{G }_{n,B}^{(b)} \right)- \eta_n^{ (a)
}\left( \mathcal{G }_{n,B^{\prime}}^{(b)} \right) \right|\leq
c(n)~|m-m^{\prime}|+\int
|[A-A^{\prime}](\varphi)|~\Sigma^{(b)}_{n,B^{\prime}}(d\varphi)
\end{equation*}
for some collection of bounded measures $\Sigma^{(b)}_{n,B^{\prime}}$ on $%
\mathcal{B }( \mathcal{A }_{n})$ such that $\int~\mbox{\rm osc}%
(\varphi)~\Sigma^{(b)}_{n,B^{\prime}}\leq
\delta\left(\Sigma^{(b)}_{n}\right) $, for some finite constant $%
\delta\left(\Sigma_{n}^{(b)}\right)<\infty$, whose values do dot depend on
the parameters $(m,A)\in(I_{n}\times \mathcal{P }( \mathcal{A }_{n}))$.
Under the assumptions (\ref{defGM}), we have
\begin{equation*}
\mathcal{G }_{n,B}^{(b)}(x)=\alpha_n^{(b)}(B)~G^{(b)}_n(x)
\end{equation*}
for some collection of parameters $\alpha_n^{(b)}(B)$ satisfying
\begin{equation*}
\left| \alpha_n^{(b)}(B)- \alpha_n^{(b)}(B^{\prime}) \right|\leq
c(n)~|m-m^{\prime}|+\int
|[A-A^{\prime}](\varphi)|~\Sigma^{(b)}_{n,B^{\prime}}(d\varphi)
\end{equation*}
This condition is clearly satisfied for the PHD model discussed in (\ref%
{phdex}), as long as the functions $h_n(y_n)+d_ng_{n}(\mbox{\LARGE .},y_n)$
are uniformly bounded from above and below.

For instance, for $b=y_n\not\in\{c,c^{\prime}\}$ we have
\begin{equation*}
\alpha_n^{(b)}(B)=\left[{h_n(b)+\int B(da)~\eta^{(a)}_n(d_ng_n(%
\mbox{\LARGE
.},b))}\right]^{-1}
\end{equation*}
In this case, we can check that
\begin{equation*}
\left| \alpha_n^{(b)}(B)- \alpha_n^{(b)}(B^{\prime}) \right|\leq
c(n)\left|[B-B^{\prime}](\varphi^{(b)}_n)\right| \quad \mbox{\rm with}
\quad\varphi^{(b)}_n(a):=\eta^{(a)}_n(d_ng_n(\mbox{\LARGE .},b))
\end{equation*}
In the same way, we show that the condition $(H_1)$ is also met for the PHD model.
This, by construction of $A^N_n$ we find that
\begin{equation*}
A^N_n=\Omega_{n}\left(\gamma_{n-1}^N(1),A_{n-1}^N\right)+\frac{1}{\sqrt{N}}~
\mathcal{W }^N_n
\end{equation*}
with some local sampling random fields satisfying (\ref{localerr}). Notice that
\begin{equation*}
\Omega_{n+1}\left(m,A\right)=\Psi_{ \mathcal{H }_{n,mA}}\left(A\right)
\mathcal{M }_{n+1,mA}(a,d(a^{\prime},b))
\end{equation*}
with the collection of potential functions
\begin{eqnarray*}
\mathcal{H }_{n,mA}(a)&:=& \mathcal{Q }_{n+1,m A}(1)(a)=\eta_n^{ (a) }\left( \mathcal{G }_{n,mA} \right)
\end{eqnarray*}
and the Markov transitions
\begin{eqnarray*}
\mathcal{M }_{n+1,mA}(a,d(a^{\prime},b))&:=&\frac{ \mathcal{Q }_{n+1,m
A}(a,d(a^{\prime},b))}{ \mathcal{Q }_{n+1,m A}(1)(a)}=\delta_{a}(da^{%
\prime})~ \frac{\nu_{n+1}(db)~\eta_n^{ (a^{\prime}) }\left( \mathcal{G }%
_{n,mA}^{(b)} \right)}{ \int \nu_{n+1}(db^{\prime})~\eta_n^{ (a^{\prime})
}\left( \mathcal{G }_{n,mA}^{(b^{\prime})} \right) }
\end{eqnarray*}

\subsection{Mixed particle association models}\label{secipsassomix}

We consider the association mapping
\begin{equation*}
\Omega_{n+1}~:~(m,A,\eta)\in \left(]0,\infty[\times \mathcal{A }%
_{[0,n]}\times \mathcal{P }(E_{n})^{ \mathcal{A }_{[0,n]}}\right)\mapsto
\Omega_{n+1}(m,A,\eta)\in \mathcal{P }( \mathcal{A }_{[0,n+1]})
\end{equation*}
defined for any $(m,A)\in \left(]0,\infty[\times \mathcal{A }_{[0,n]}\right)$
and any mapping $\eta : a\in \mbox{\rm Supp}(A)\mapsto \eta^{(a)}\in Pa(E_n)$
by
\begin{equation*}
\begin{array}{l}
\Omega_{n+1}\left(m,A,\eta\right)(d(a,b)) \propto
A(da)~\nu_{n+1}(db)~\eta^{(a)}\left( G^{(b)}_{n,m \int A(da)~\eta^{(a)}}
\right)%
\end{array}%
\end{equation*}
By construction, for any discrete measure $A\in \mathcal{P }( \mathcal{A }%
_{[0,n-1]})$, and any mapping $a\in \mbox{\rm Supp}(A) \mapsto \eta^{(a)}\in
\mathcal{P }(E_{n-1})$, we have the formula
\begin{equation*}
\Gamma^2_{n}\left(m, \int A(da)~\eta^{(a)}) \right)=\int
\Omega_{n}\left(m,A,\eta^{(\mbox{\LARGE .})}\right)(d(a,b))~~\Phi_{n}^{(b)}%
\left(\eta^{(a)}\right)
\end{equation*}
We also mention that the updating-prediction transformation defined in (\ref%
{reca})
\begin{equation}  \label{reca2}
\Phi_{n}^{(a)}\left(\eta\right)=\Psi^{(a)}_{G_{n-1}}\left(\eta%
\right)M^{(a)}_{n}=\eta \mathcal{K}^{(a)}_{n,\eta}\quad \mbox{\rm with}\quad
\mathcal{K}^{(a)}_{n,\eta}=\mathcal{S}^{(a)}_{n-1,\eta}M^{(a)}_n
\end{equation}
In the above displayed formula $\mathcal{S}^{(a)}_{n,\eta}$ stands for some
updating Markov transition from $E_{n-1}$ into itself satisfying the
compatibility condition $\eta\mathcal{S}^{(a)}_{n-1,\eta}=%
\Psi^{(a)}_{G_{n-1}}\left(\eta\right)$.

We let $A^N_0=\frac{1}{N}\sum_{i=1}^N\delta_{\mathrm{a}^i_0}$, be the
empirical measure associated with $N$ independent and identically
distributed random variables $(\mathrm{a}^i_0)_{1\leq i\leq N}$ with common
distribution $A_0$. For any $\mathrm{a}\in \mathcal{A }_0$, we let
\begin{equation*}
\eta^N_0:=\int A^N_0(d\mathrm{a})~\eta^{(\mathrm{a},N^{\prime})}_0\quad%
\mbox{\rm and}\quad \eta^{(\mathrm{a},N^{\prime})}_0=\frac{1}{N^{\prime}}%
\sum_{i=1}^{N^{\prime}}~\delta_{\xi^{[\mathrm{a},j]}_0}
\end{equation*}
with the empirical measure $\eta^{(\mathrm{a},N^{\prime})}_0$ associated
with $N^{\prime}$ random variables $\xi^{[\mathrm{a}]}_0=\left(\xi^{[\mathrm{%
a},j]}_0\right)_{1\leq j\leq N^{\prime}}$ with common law $\eta^{(\mathrm{a}%
)}_0$. We further assume that $\gamma_0(1)$ is known and set
\begin{equation*}
\gamma_0^N:=\gamma_0(1)~\eta^N_0\quad\mbox{\rm and} \quad
\gamma_1^N(1):=\gamma_0^N(1)~\eta^N_0(G_{0,\gamma^N_0})
\end{equation*}

It is readily checked that the fluctuation random fields given below
\begin{equation*}
\mathcal{W }_0^{(a,N^{\prime})}=\sqrt{N^{\prime}}~\left(
\eta^{(a,N^{\prime})}_0-\eta^{(a)}_0\right)
\end{equation*}
satisfies (\ref{localerr}), with $N=N^{\prime}$, for any given $a\in
\mathcal{A }_0$. Using the fact that
\begin{equation*}
\int A^N_0(d\mathrm{a})~\eta^{(\mathrm{a},N^{\prime})}_0=\int A^N_0(d\mathrm{%
a})~\eta^{(\mathrm{a})}_0+\frac{1}{\sqrt{N^{\prime}}}~ \int A^N_0(d\mathrm{a}%
)~ \mathcal{W }_0^{(a,N^{\prime})}
\end{equation*}
we conclude that
\begin{equation*}
\eta^N_0:=\eta_0+\frac{1}{\sqrt{N}}~W^N_0
\end{equation*}
with some local sampling random fields $W^N_0$ satisfying the same estimates
as in (\ref{localerr}) by replacing $1/\sqrt{N}$ by the sum $\left(1/\sqrt{N}%
+1/\sqrt{N^{\prime}}\right)$.

Using (\ref{reca2}), for any $\mathrm{a}_1=(a_0,a_1)$ we find that
\begin{equation*}
\Phi^{(a_1)}_1\left(\eta^{(a_0,N^{\prime})}_0\right)=\eta^{(a_0,N^{%
\prime})}_0 \mathcal{K}^{(a_1)}_{n,\eta^{(a_0,N^{\prime})}_0}
\end{equation*}
We let $A^N_1=\frac{1}{N}\sum_{i=1}^N\delta_{\mathrm{a}^i_{1}}$ be the
occupation measure associated with $N$ conditionally independent and
identically distributed random variables $\mathrm{a}%
^i_1:=(a^i_{0,1},a^i_{1,1})$ with common law
\begin{equation*}
\Omega_{1}\left(\gamma_0^N(1),A_0^N,\eta^{(\mbox{\LARGE .}%
,N^{\prime})}_0\right)
\end{equation*}
In the above displayed formula $\eta^{(\mbox{\LARGE .},N^{\prime})}_0$
stands for the mapping $a_0\in \mathcal{A }_0\mapsto\eta^{(a_0,N^{%
\prime})}_0\in \mathcal{P }(E_0)$.

We consider a sequence of conditionally independent random variables $%
\xi^{[a_0,a_1,j]}_1$ with distribution $\mathcal{K}^{(a_1)}_{n,%
\eta^{(a_0,N^{\prime})}_0}\left(\xi^{[a_0,j]}_0,\mbox{\LARGE .}\right)$,
with $1\leq j\leq N^{\prime}$, and we set
\begin{equation*}
\eta^{((a_0,a_1),N^{\prime})}_1=\frac{1}{N^{\prime}}\sum_{i=1}^{N^{\prime}}~%
\delta_{\xi^{[(a_0,a_1),j]}_1}\quad\mbox{\rm and}\quad \eta^N_1:=\int A^N_1(d%
\mathrm{a})~\eta^{(\mathrm{a},N^{\prime})}_1
\end{equation*}
Arguing as before, for any given $\mathrm{a}_1:=(a_0,a_1)\in \mathrm{Supp}%
(A^N_1)$, the sequence of random fields
\begin{equation*}
\mathcal{W }_1^{(\mathrm{a}_1N^{\prime})}:=\sqrt{N}~\left(
\eta^{((a_0,a_1),N^{\prime})}_1-\Phi^{(a_1)}_1\left(\eta^{(a_0,N^{%
\prime})}_0\right) \right)
\end{equation*}
satisfies (\ref{localerr}), with $N=N^{\prime}$. Thus, we conclude that
\begin{eqnarray*}
\eta^N_1&=&\int \Omega_{1}\left(\gamma_0^N(1),A_0^N,\eta^{(\mbox{\LARGE .}%
,N^{\prime})}_0\right)(d(a_0,a_1))~\Phi^{(a_1)}_1\left(\eta^{(a_0,N^{%
\prime})}_0\right) +\frac{1}{\sqrt{N}}~W^N_1 \\
&=&\Gamma^2_{1}\left(\gamma_{0}^N(1),\eta^N_{0}\right)+\frac{1}{\sqrt{N}}%
~W^N_1
\end{eqnarray*}
with some local sampling random fields $W^N_1$ satisfying the same estimates
as in (\ref{localerr}) by replacing $1/\sqrt{N}$ by the sum $\left(1/\sqrt{N}%
+1/\sqrt{N^{\prime}}\right)$. Iterating this procedure, we define by
induction a sequence of $N$-particle approximation measures
\begin{equation*}
\gamma_n^N(1)=\gamma_{n-1}^N(1)~\eta^N_{n-1}(G_{n-1,\gamma^N_{n-1}})\quad%
\mbox{\rm and}\quad \eta^N_n:=\int A^N_n(da)~\eta^{(a,N^{\prime})}_n
\end{equation*}
with the occupation measure $A^N_n=\frac{1}{N}\sum_{i=1}^N\delta_{\mathrm{a}%
^i_{n}}$ associated with $N$ conditionally independent and identically
distributed random variables $\mathrm{a}^i_n:=(a^i_{0,n},a^i_{1,n},%
\ldots,a_{n,n}^i)$ with common law $\Omega_{n}\left(%
\gamma_{n-1}^N(1),A_{n-1}^N,\eta^{(\mbox{\LARGE .},N^{\prime})}_{n-1}\right)$%
. Arguing as above, we find that
\begin{equation*}
\eta^N_n=\int \Omega_{n}\left(\gamma_{n-1}^N(1),A_{n-1}^N,\eta^{(%
\mbox{\LARGE .},N^{\prime})}_{n-1}\right)(d(a,b))~\Phi_{n}^{(b)}\left(%
\eta^{(a,N^{\prime})}_{n-1}\right)
=\Gamma^2_{n}\left(\gamma_{n-1}^N(1),\eta^N_{n-1}\right)+\frac{1}{\sqrt{N}}%
~W^N_n
\end{equation*}
with some local sampling random fields satisfying the same estimates as in (%
\ref{localerr}) by replacing $1/\sqrt{N}$ by the sum $\left(1/\sqrt{N}+1/%
\sqrt{N^{\prime}}\right)$. As before, the $N$-particle occupation measures $%
A^N_n$ converge as $N$ tends to $\infty$ to the association probability
measures $A_n$.
\section{Appendix}
\subsection{Proof of corollary~\ref{propp}}
For constant mappings $s_n=\mu_{n+1}(1)$, the mappings $ \Phi^1_{n+1,\nu_n}$ and $ \Phi^2_{n+1,m_n}$ are given by
$$
 \Phi^1_{n+1,\nu_n}(u)=s_n
 \quad\mbox{\rm and}\quad
 \Phi^2_{n+1,m_n}(\eta)=
 \Psi_{g^{(s)}_n}(\eta)M_{n+1}^{(s)}
 $$
with the likelihood function $g^{(s)}_n$ and the Markov transitions $M^{(s)}_{n+1}$ defined in (\ref{gMs}).   Firstly, we observe that
$
r_n(s_n):=\sup_{x,x^{\prime}\in E_n}{g^{(s)}_n(x)}/{g^{(s)}_n(x^{\prime})}$.
We also notice that the second component mapping $\Phi^2_{n+1,m_n}$ does not depends on the parameter $m_n$, and it induces a
 Feynman-Kac semigroup of the same form as the one discussed in section~\ref{fksg}.

 Under the premise of the proposition, the semigroup of associated with the Markov transitions $M_{n}$ satisfies
the mixing property stated in the l.h.s. of (\ref{Hm})
for some integer $m\geq 1$ and some parameter $\epsilon_p(m)\in ]0,1]$.
In this situation, we also have that
$$
M_{p,p+m}^{(s)}(x,\point)\geq \epsilon_p^{(s)}(m)~M^{(s)}_{p,p+m}(x^{\prime},\point)
$$
with some positive parameter
$$
\epsilon_p^{(s)}(m)\geq \epsilon_p(m)/\prod_{p\leq k<p+m}r_k(s_k)r_k(1)\quad\mbox{\rm and}
\quad
r_n(s_n):=
\frac{s_ng^+_n+(1-s_n)}{s_ng^-_n+(1-s_n)}\left(\leq r_n(1)\right)
$$
To prove this claim, firstly we observe that
$
M_{p,p+m}^{(s)}(x,\point)\ll M_{p,p+m}^{(s)-}(x,\point)
$
and
$$
\prod_{p\leq k<p+m} r_k(s_k)^{-1}
\leq
{dM_{p,p+m}^{(s)}(x,\point)}/{dM_{p,p+m}^{(s)-}(x,\point)}
\leq \prod_{p\leq k<p+m}
r_k(1)
$$
with the semigroup $M_{p,n}^{(s)-}$ associated with the Markov transition
$$
M_{p,p+1}^{(s)-}(x,\point)=\alpha_{p+1}~M_{p+1}(x,\point)+(1-\alpha_{k+1})~\overline{\mu}_{k+1}
\quad\mbox{\rm with}\quad
\alpha_{p+1}:=\frac{s_kg_k^-}{s_kg^-_k+(1-s_k)}
$$
Using the geometric representation
$$
M_{p,n}^{(s)-}(x,\point)=\left(\prod_{p<k\leq n}\alpha_k\right)M_{p,n}(x,\point)+\sum_{p<k\leq n}
(1-\alpha_k)~\left(\prod_{k<l\leq n}\alpha_l\right)~\overline{\mu}_kM_{k,n}
$$
it can be verified that
$$
M_{p,p+m}^{(s)-}(x,\point)\geq \epsilon_p(m)~M_{p,p+m}^{(s)-}(x^{\prime},\point)\geq
\epsilon_p(m)~\left(\prod_{p\leq k<p+m}
{g^-_k}/{g^+_k}\right)~M_{p,p+m}^{(s)}(x^{\prime},\point)
$$
from which we conclude that
$$
M_{p,p+m}^{(s)}(x,\point)\geq \epsilon_p^{(s)}(m)~M^{(s)}_{p,p+m}(x^{\prime},\point)
\quad\mbox{\rm with}\quad
\epsilon_p^{(s)}(m)\geq \epsilon_p(m)/\prod_{p\leq k<p+m}r_k(s_k)r_k(1)
$$
We end the proof of the proposition combing the proposition~\ref{keyprop} with the couple of estimates presented in (\ref{rpn}) and
(\ref{betapn}).
This ends the proof of the corollary.
\cqfd

\subsection{Proof of theorem~\ref{theoBernoulli}}\label{ptheoBernoulli}

The formulae presented  in (\ref{massber}) can be rewritten in terms of matrix operations as follows
$$
\left[\gamma_{n+1}(1)~,~1-\gamma_{n+1}(1)\right]=\left[\widehat{\gamma}_n(1)~,~1-\widehat{\gamma}_n(1)\right]~\left[
\begin{array}{cc}
\Psi_{g_n}(\eta_n)(s_n)&1-\Psi_{g_n}(\eta_n)(s_n)\\
\mu_{n+1}(1)&1-\mu_{n+1}(1)
\end{array}
\right]
$$
and
$$
\left[\widehat{\gamma}_n(1)~,~1-\widehat{\gamma}_n(1)\right]=\frac{\left[\gamma_n(1)~,~1-\gamma_n(1)\right]~\left[
\begin{array}{cc}
\eta_n(g_n)&0\\
0&1
\end{array}
\right]}{\left[\gamma_n(1)~,~1-\gamma_n(1)\right]~\left[
\begin{array}{cc}
\eta_n(g_n)&0\\
0&1
\end{array}
\right]~\left[
\begin{array}{c}
1\\
1
\end{array}
\right]}
$$
With a slight abuse of notation, we set
$$ \vartheta_n:=\left[\gamma_{n}(1)~,~1-\gamma_{n}(1)\right]
\quad
 \widehat{\vartheta}_n:=\left[\widehat{\gamma}_n(1)~,~1-\widehat{\gamma}_n(1)\right]\quad
 \mbox{\rm and}\quad
 1=\left[
\begin{array}{c}
1\\
1
\end{array}
\right]\quad
$$
We also denote by $\Ma_{n+1,\eta_n}$ and $\Da_{n,\eta_n}$ the stochastic and the diagonal
matrices defined by
\begin{equation}\label{defMG}
\Ma_{n+1,\eta_n}:=~\left[
\begin{array}{cc}
\Psi_{g_n}(\eta_n)(s_n)&1-\Psi_{g_n}(\eta_n)(s_n)\\
\mu_{n+1}(1)&1-\mu_{n+1}(1)
\end{array}
\right]
\quad\mbox{\rm and}\quad
\Da_{n,\eta_n}:=\left[
\begin{array}{cc}
\eta_n(g_n)&0\\
0&1
\end{array}
\right]\end{equation}
In this notation, the above recursion can be rewritten in a more compact form
$$
\vartheta_{n+1}=\widehat{\vartheta}_n~\Ma_{n+1,\eta_n}
\quad\mbox{\rm and}\quad
 \widehat{\vartheta}_n=\frac{\vartheta_n~\Da_{n,\eta_n}
 }{
\vartheta_n\Da_{n,\eta_n} 1
}~\Longrightarrow~
\vartheta_{n+1}=\frac{\vartheta_n~\Qa_{n+1,\eta_n}
 }{
\vartheta_n\Qa_{n+1,\eta_n} 1
}
$$
with the product of matrices $\Qa_{n+1,\eta_n}=\Da_{n,\eta_n}\Ma_{n+1,\eta_n}$.
$$
\forall u\in I_{p}(\subset[0,1])\qquad
\left[\Phi^1_{p,n,\nu}(u),1-\Phi^1_{p,n,\nu}(u)\right]=\frac{[u,1-u]~\Qa_{p,n,\nu}}{[u,1-u]~\Qa_{p,n,\nu}(1)}
$$
with the matrix semigroup
$$
\Qa_{p,n,\nu}=\Qa_{p+1,\nu_p}\Qa_{p+2,\nu_{p+1}}\ldots \Qa_{n,\nu_{n-1}}
$$
These semigroups are again of the same form as the Feynman-Kac models discussed in section~\ref{fksg} with a two point state space. When $\mu_{n+1}(1)\in]0,1[$ and $0<s^-_n\leq s_n^+<1$, we have
for any $n\geq 0$ and any $i,i^{\prime},j\in\{1,2\}$
$$
\Ma_{n+1,\nu_n}(i,j)\geq \epsilon_n~\Ma_{n+1,\nu_n}(i^{\prime},j)\quad\mbox{\rm and}\quad
\sup_{i,i^{\prime}\in\{1,2\}}\frac{\Qa_{n+1,\nu_n}(1)(i)}{\Qa_{n+1,\nu_n}(1)(i^{\prime})}\leq \delta^{\prime}_{n}(g)
$$
The first assertion is a direct consequence of the proposition~\ref{keyprop} with the couple of estimates presented in (\ref{rpn}) and
(\ref{betapn}).

Using (\ref{bernm}), we find that $\Phi^2_{n+1,m_n}$ induces
a Feynman-Kac models of the same form as the one discussed in section~\ref{fksg}. More precisely, we have that
$$
\Phi^2_{n+1,m_n}(\eta)=\Psi_{\widehat{G}_{n,m_n}}(\eta)\widehat{M}_{n+1,m_n}
$$
with the potential functions $G_{n,m_n}$ and the Markov transitions $\widehat{M}_{n+1,m_n}$ defined in (\ref{bernmG}) and (\ref{bernmM}). Notice that
$$
\sup_{x,x^{\prime}\in E_n}\frac{\widehat{G}_{n,m_n}(x)}{\widehat{G}_{n,m_n}(x^{\prime})}\leq \delta_n(sg)
$$
and for any $x\in E_n$ and any $n\geq 0$
$$
\delta_n(sg)^{-1}~\widehat{M}^-_{n+1,m_n}(x,\point) \leq \widehat{M}_{n+1,m_n}(x,\point) \leq \delta_n(sg)~\widehat{M}^-_{n+1,m_n}(x,\point)
$$
with the Markov transitions $\widehat{M}^-_{n+1,m_n}$ defined as $\widehat{M}^-_{n+1,m_n}$
by replacing the functions $(s_n,g_n)$ by their lower bounds $(s^-_n,g_n^-)$.
To prove this claim, we use the fact that for any positive function $f$ we have
$$
\frac{d\widehat{M}_{n+1,m_n}(f)}{d\widehat{M}^-_{n+1,m_n}(f)}=
\frac{m_ng^-_ns^-_n+(1-m_n)\mu_{n+1}(1)}{m_ng_ns_n+(1-m_n)\mu_{n+1}(1)}\times
\frac{m_ng_ns_nM_{n+1}(f)+(1-m_n)\mu_{n+1}(1)\overline{\mu}_{n+1}(f)}{m_ng^-_ns^-_n+(1-m_n)\mu_{n+1}(1)\overline{\mu}_{n+1}(f)}
$$
and the two series of inequalities
$$
\delta_n(sg)^{-1}\leq \frac{m_ng^-_ns^-_n+(1-m_n)\mu_{n+1}(1)}{m_ng_ns_n+(1-m_n)\mu_{n+1}(1)}\leq 1
$$
and
$$
1\leq \frac{m_ng_ns_nM_{n+1}(f)+(1-m_n)\mu_{n+1}(1)\overline{\mu}_{n+1}(f)}{m_ng^-_ns^-_n+(1-m_n)\mu_{n+1}(1)\overline{\mu}_{n+1}(f)}\leq\delta_n(sg)
$$

With a slight abuse of notation, we write $\widehat{M}_{p,n}$, and respectively $\widehat{M}^-_{p,n}$, the semigroup associated with the Markov transitions $ \widehat{M}_{n+1,m_n}$, and resp. $ \widehat{M}_{n+1,m_n}^-$. Using the same argument as in the proof of corollary~\ref{propp} it follows that
 $$
 \widehat{M}_{p,p+m}^-(x,\point)\geq\epsilon_p(m)~\widehat{M}_{p,p+m}^-(x^{\prime},\point)
 $$
 from which we conclude that
$$
\widehat{M}_{p,p+m}(x,\point)\geq \widehat{\epsilon}_p(m)~\widehat{M}_{p,p+m}(x^{\prime},\point)
\quad\mbox{\rm
with}\quad \widehat{\epsilon}_p(m)\geq \epsilon_p(m)~\prod_{0\leq k<m}\delta_{p+k}(sg)^{-2}
$$
using proposition~\ref{keyprop} with the couple of estimates presented in (\ref{rpn}) and
(\ref{betapn}), we check that (\ref{a2}) is satisfied with
$$
a^2_{p,n}\leq 2~\rho_p(m)
\prod_{k=0}^{\lfloor n/m\rfloor-1}
\left(1-\epsilon_{p+km}^{(m)}\right)
$$
and some parameters
$$
\epsilon_{p}^{(m)}:=\widehat{\epsilon}_{p}(m)^2~\prod_{0<k<m}\delta_{p+k}(sg)^{-1}\geq
\epsilon_p(m)^2~\delta_{p}(sg)^{-4}\prod_{0<k<m}\delta_{p+k}(sg)^{-5}
$$
and
$$
\rho_p(m):=\widehat{\epsilon}_p(m)^{-1}~\prod_{0\leq k<m}\delta_{p+k}(sg)
\leq \epsilon_p(m)^{-1}~\prod_{0\leq k<m}\delta_{p+k}(sg)^{3}
$$
This ends the proof of the first assertion of the  theorem. Next, we discuss condition $(\mbox{Cont$(\Phi)$})$.
We observe that
$$
\Phi_{n+1,\nu}^1(u)=\frac{u~\nu(g_ns_n)+(1-u)\mu_{n+1}(1)}{u~\nu(g_n)+(1-u)}
$$
After some manipulations
$$
\begin{array}{l}
\Phi_{n+1,\nu}^1(u)-\Phi_{n+1,\nu^{\prime}}^1(u)\\
\\
=
\frac{u\nu^{\prime}(g_n)}{u\nu^{\prime}(g_n)+(1-u)}~\left[\Psi_{g_n}(\nu)-\Psi_{g_n}(\nu^{\prime})\right](s_n)\\
\\\qquad+
\frac{u}{u\nu(g_n)+(1-u)}\frac{(1-u)}{u\nu^{\prime}(g_n)+(1-u)}~\left[\Psi_{g_n}(\nu)(s_n)-\mu_{n+1}(1)\right]~
\left[\nu-\nu^{\prime}\right](g_n)
\end{array}
$$
Recalling that the mapping $\theta_a(x)=ax/(ax+(1-x))$ in increasing on $[0,1]$ and using the fact that
$$
\Psi_{g_n}(\nu)=\nu S_{n,\nu}\Longrightarrow \Psi_{g_n}(\nu)-\Psi_{g_n}(\nu^{\prime})=\frac{g^+_n}{\nu(g_n)}
~(\nu-\nu^{\prime}) S_{n,\nu^{\prime}}
$$
with the Markov transition
$$
S_{n,\nu^{\prime}}(x,dx^{\prime})=\frac{g_n(x)}{g_n^+(x)}~\delta_x(dx^{\prime})+\left(1-\frac{g_n(x)}{g_n^+(x)}\right)~\Psi_{g_n}(\nu^{\prime})(dx^{\prime})
$$
we prove
\begin{equation}\label{refboltzman}
\left|\Psi_{g_n}(\nu)(s_n)-\Psi_{g_n}(\nu^{\prime})(s_n)\right|\leq \frac{g^+_n}{g^-_n}~\left|(\nu-\nu^{\prime})S_{n,\nu^{\prime}}(s_n)\right|
\end{equation}
and for any $u\in I_n=[m^-_n,m^+_n]$
$$
\begin{array}{l}
\left|\Phi_{n+1,\nu}^1(u)-\Phi_{n+1,\nu^{\prime}}^1(u)\right|\\
\\
=
\frac{m^+_ng_n^+}{m^+_ng_n^++(1-m^+_n)}~\frac{g^+_n}{g^-_n}~\left|(\nu-\nu^{\prime})S_{n,\nu^{\prime}}(s_n)\right|
\\
\\\qquad+
\frac{m^+_ng_n^+}{m^+_ng_n^++(1-m^+_n)}
\frac{(1-m^-_n)}{m^-_ng_n^-+(1-m_n^-)}~\left\|s_n-\mu_{n+1}(1)\right\|~
\left|\left[\nu-\nu^{\prime}\right](g_n/g^-_n)\right|
\end{array}
$$
This implies that
\begin{eqnarray*}
\tau^1_{n+1}&\leq&
\frac{m^+_ng_n^+}{m^+_ng_n^++(1-m^+_n)}~\frac{g^+_n}{g^-_n}~(s^+_n-s^-_n)\\
&&\qquad+
\frac{m^+_ng_n^+}{m^+_ng_n^++(1-m^+_n)}~
\frac{(1-m^-_n)}{m^-_ng_n^-+(1-m_n^-)}~\left\|s_n-\mu_{n+1}\right\|~\left(\frac{g^+_n}{g^-_n}-1\right)\\
&\leq&\frac{g^+_n}{g^-_n}~\left[(s^+_n-s^-_n)+\left\|s_n-\mu_{n+1}(1)\right\|\right]
\end{eqnarray*}

Using (\ref{bernm}) we also find that
$$
\Phi_{n+1,m}^2(\eta)(f)=\frac{m\eta(s_ng_nM_{n+1}(f))+(1-m)~\mu_{n+1}(f)}{m\eta(s_ng_n)+(1-m)\mu_{n+1}(1)}
$$
It is also readily check that
$$
\left[\Phi_{n+1,m}^2(\eta)-\Phi_{n+1,m^{\prime}}^2(\eta)\right](f)
=
\frac{\mu_{n+1}(1)~\eta(g_ns_n)~\left[
\Psi_{g_ns_n}(\eta)M_{n+1}-\overline{\mu}_{n+1}
\right](f)~(m-m^{\prime})
}{\left[m\eta(s_ng_n)+(1-m)\mu_{n+1}(1)\right]\left[m^{\prime}\eta(s_ng_n)+(1-m^{\prime})\mu_{n+1}(1)\right]}$$
from which we conclude that
$$
\tau^2_{n+1}\leq \sup{\left\{\frac{\mu_{n+1}}{s^-_ng^-_n},\frac{s^+_ng^+_n}{\mu_{n+1}(1)}\right\}}\leq
\delta^{\prime}_n(g) ~
 \sup{\left\{\frac{\mu_{n+1}}{s^-_n},\frac{s^+_n}{\mu_{n+1}(1)}\right\}}
$$
This ends the proof of the theorem.
\cqfd
\subsection{Proof of proposition~\ref{propsgpn}}\label{propsgpnp}

The proof of proposition~\ref{propsgpn} is based on the following technical lemma.

\begin{lem}\label{lemtech}
We assume that the regularity conditions $(\mbox{Lip$(\Phi)$})$ and $(\mbox{Cont$(\Phi)$})$  are satisfied. In this situation, for any $p\leq n$, $u,u^{\prime}\in I_p$, $\eta,\eta^{\prime}\in\Pa(E_p)$ and $f\in \mbox{Osc}_{1}(E_{n})$ and any flow of masses
and probability
measures
 $m=(m_n)_{n\geq 0}\in \prod_{n\geq 0}I_n$ and $\nu:=(\nu_n)_{n\geq 0}\in\prod_{n\geq 0}\Pa(E_n)$ we have the following estimates
$$
\left|\Phi_{p,n,\nu^{\prime}}^1(u^{\prime})-\Phi_{p,n,\nu}^1(u)\right|\leq
 a^1_{p,n}~|u-u^{\prime}|+\sum_{p\leq q<n}~\overline{a}^{1}_{q,n}~\int~\left|[\nu_{q}-\nu^{\prime}_{q}](\varphi)\right|~\Omega^1_{q+1,\nu^{\prime}_{q}}(d\varphi)
$$
$$
\begin{array}{l}
\left|\Phi_{p,n,m^{\prime}}^2(\eta^{\prime})(f)-\Phi_{p,n,m}^2(\eta)(f)\right|
\leq a^2_{p,n} \displaystyle\int\left|[\eta-\eta^{\prime}](\varphi)\right|
\Omega^2_{p,n,\eta^{\prime}}(f,d\varphi)
+\sum_{p\leq q< n}~\overline{a}^2_{q,n}~|m_{q}-m^{\prime}_{q}|
\end{array}$$
with the collection of parameters $\overline{a}^i_{p,n}$, $i=1,2$, defined in (\ref{overa}).
\end{lem}
\proof
We use the decomposition
\begin{eqnarray*}
\Phi_{p,n,\nu^{\prime}}^1(u^{\prime})-\Phi_{p,n,\nu}^1(u)
&=&\Phi_{p,n,\nu}^1(u^{\prime})-\Phi_{p,n,\nu}^1(u)\\
&&\qquad\qquad+\sum_{p<q\leq n}
\left[\Phi_{q,n,\nu}^1(\Phi_{p,q,\nu^{\prime}}^1(u^{\prime}))
-
\Phi_{q-1,n,\nu}^1(\Phi_{p,q-1,\nu^{\prime}}^1(u^{\prime}))
\right]
\end{eqnarray*}
and the fact that
\begin{eqnarray*}
\Phi_{q-1,n,\nu}^1(\Phi_{p,q-1,\nu^{\prime}}^1(u^{\prime}))&=&\Phi_{q,n,\nu}^1\left(
\Phi_{q-1,q,\nu}^1\left[\Phi_{p,q-1,\nu^{\prime}}^1(u^{\prime})\right]
\right)\\
\Phi_{q,n,\nu}^1(\Phi_{p,q,\nu^{\prime}}^1(u^{\prime}))&=&
\Phi_{q,n,\nu}^1\left(\Phi_{q-1,q,\nu^{\prime}}^1\left[\Phi_{p,q-1,\nu^{\prime}}^1(u^{\prime})\right]
\right)
\end{eqnarray*}
and
$$
\left|\Phi_{p,n,\nu}^1(u^{\prime})-\Phi_{p,n,\nu}^1(u)\right|\leq a^1_{p,n}~|u-u^{\prime}|
$$
and
$$
\begin{array}{l}
\left|\Phi_{q,n,\nu}^1(\Phi_{p,q,\nu^{\prime}}^1(u^{\prime}))
-
\Phi_{q-1,n,\nu}^1(\Phi_{p,q-1,\nu^{\prime}}^1(u^{\prime}))
\right|\\
\\
\leq a^{1}_{q,n}~\left|
\Phi_{q,\nu_{q-1}}^1\left[\Phi_{p,q-1,\nu^{\prime}}^1(u^{\prime})\right]
-\Phi_{q,\nu^{\prime}_{q-1}}^1\left[\Phi_{p,q-1,\nu^{\prime}}^1(u^{\prime})\right]
\right|
\\
\\
\leq ~\overline{a}^{1}_{q-1,n}~\int~\left|[\nu_{q-1}-\nu^{\prime}_{q-1}](\varphi)\right|~\Omega_{q,\nu^{\prime}_{q-1}}(d\varphi)
\end{array}$$
to show that
$$
\left|\Phi_{p,n,\nu^{\prime}}^1(u^{\prime})-\Phi_{p,n,\nu}^1(u)\right|\leq
 a^1_{p,n}~|u-u^{\prime}|+\sum_{p<q\leq n}~\overline{a}^{1}_{q-1,n}~\int~\left|[\nu_{q-1}-\nu^{\prime}_{q-1}](\varphi)\right|~\Omega^1_{q,\nu^{\prime}_{q-1}}(d\varphi)
$$
In the same way, we use the decomposition
\begin{eqnarray*}
\left[\Phi_{p,n,m^{\prime}}^2(\eta^{\prime})-\Phi_{p,n,m}^2(\eta)\right]
&=&\left[\Phi_{p,n,m}^2(\eta^{\prime})-\Phi_{p,n,m}^2(\eta)\right]\\
&&\qquad+\sum_{p<q\leq n}
\left[\Phi_{q,n,m}^2(\Phi_{p,q,m^{\prime}}^2(\eta^{\prime}))
-
\Phi_{q-1,n,m}^2(\Phi_{p,q-1,m^{\prime}}^2(\eta^{\prime}))
\right]
\end{eqnarray*}
and the fact that
\begin{eqnarray*}
\Phi_{q-1,n,m}^2(\Phi_{p,q-1,m^{\prime}}^2(\eta^{\prime}))&=&\Phi_{q,n,m}^2\left(
\Phi_{q-1,q,m}^2\left[\Phi_{p,q-1,m^{\prime}}^2(\eta^{\prime})\right]
\right)\\
\Phi_{q,n,m}^2(\Phi_{p,q,m^{\prime}}^2(\eta^{\prime}))&=&
\Phi_{q,n,m}^2\left(\Phi_{q-1,q,m^{\prime}}^2\left[\Phi_{p,q-1,m^{\prime}}^2(\eta^{\prime})\right]
\right)
\end{eqnarray*}
and
$$
\left|\Phi_{p,n,m}^2(\eta^{\prime})(f)-\Phi_{p,n,m}^2(\eta)(f)\right|\leq a^2_{p,n} \int~\left|[\eta-\eta^{\prime}](\varphi)\right|~%
\Omega^2_{p,n,\eta^{\prime}}(f,d\varphi)
$$
to show that
$$
\begin{array}{l}
\left|\Phi_{q,n,m}^2(\Phi_{p,q,m^{\prime}}^2(\eta^{\prime}))
-
\Phi_{q-1,n,m}^2(\Phi_{p,q-1,m^{\prime}}^2(\eta^{\prime}))
\right|\\
\\
\leq a^2_{q,n} \displaystyle\int~\left|[\Phi_{q,m_{q-1}}^2\left[\Phi_{p,q-1,m^{\prime}}^2(\eta^{\prime})\right]-\Phi_{q,m^{\prime}_{q-1}}^2\left[\Phi_{p,q-1,m^{\prime}}^2(\eta^{\prime})\right]
](\varphi)\right|~%
\Omega^2_{q,n,\Phi_{p,q,m^{\prime}}^2(\eta^{\prime})}(f,d\varphi)

\\
\\
\leq  \overline{a}^2_{q-1,n}~|m_{q-1}-m^{\prime}_{q-1}|
\end{array}$$
Using these estimates we conclude that
$$
\left|\left[\Phi_{p,n,m^{\prime}}^2(\eta^{\prime})-\Phi_{p,n,m}^2(\eta)\right](f)\right]
\leq a^2_{p,n} \int~\left|[\eta-\eta^{\prime}](\varphi)\right|~%
\Omega^2_{p,n,\eta}(f,d\varphi)
+\sum_{p<q\leq n}~\overline{a}^2_{q-1,n}~|m_{q-1}-m^{\prime}_{q-1}|
$$
This ends the proof of the lemma.
\cqfd
Now we come to the proof of  proposition~\ref{propsgpn}.

{\bf Proof of proposition~\ref{propsgpn}:}

We fix a parameter $p\geq 0$, and we let $(m_n)_{n\geq p},(m^{\prime}_n)_{n\geq p}\in \prod_{n\geq p}I_n$ and $(\nu_n)_{n\geq p}$, and $(\nu^{\prime}_n)_{n\geq p}\in\prod_{n\geq p}\Pa(E_n)$ be defined by the following recursive formulae
$$
\forall q>p\qquad m_q^{\prime}=\Phi_{q,\nu_{q-1}^{\prime}}^1(m_{q-1}^{\prime})\quad
\mbox{\rm and}\quad
\nu_q^{\prime}=\Phi_{q,m_{q-1}^{\prime}}^2(\nu_{q-1}^{\prime})
$$
$$
\forall q>p\qquad m_q=\Phi_{q,\nu_{q-1}}^1(m_{q-1})\quad
\mbox{\rm and}\quad
\nu_q=\Phi_{q,m_{q-1}}^2(\nu_{q-1})
$$
with the initial condition for $q=p$
$$
(\nu_p,\nu_p^{\prime})=(\eta,\eta^{\prime})
\quad\mbox{\rm and}\quad
(m_p,m_p^{\prime})=(u,u^{\prime})
$$
By construction, we have
$$
\nu_q^{\prime}=\Phi_{p,q,m^{\prime}}^2(\eta^{\prime})\quad\mbox{\rm and}\quad
\nu_q=\Phi_{p,q,m}^2(\eta)
$$
as well as
$$
m_q^{\prime}=\Phi_{p,q,\nu^{\prime}}^1(u^{\prime})\quad\mbox{\rm and}\quad
m_q=\Phi_{p,q,\nu}^1(u)
$$
In this case, using lemma~\ref{lemtech} it follows that
$$
\begin{array}{l}
\left|\left[\Gamma_{p,n}^2(m^{\prime},\eta^{\prime})-\Gamma^2_{p,n}(m,\eta)\right](f)\right]
\\
\\
\leq a^2_{p,n} \displaystyle\int~\left|[\eta-\eta^{\prime}](\varphi)\right|~%
\Omega^2_{p,n,\eta^{\prime}}(f,d\varphi)
+\sum_{p\leq q< n}~\overline{a}^2_{q,n}~|\Gamma_{p,q}^1(m^{\prime},\eta^{\prime})-\Gamma^1_{p,q}(m,\eta)|
\end{array}$$
and
$$
\begin{array}{l}
\left|\Gamma_{p,n}^1(m^{\prime},\eta^{\prime})-\Gamma_{p,n}^1(m,\eta)\right|\\
\\
\leq
 a^1_{p,n}~|m-m^{\prime}|+\displaystyle\sum_{p\leq q<n}~\overline{a}^{1}_{q,n}~\displaystyle\int~\left|[\Gamma_{p,q}^2(m^{\prime},\eta^{\prime})-\Gamma^2_{p,q}(m,\eta)](\varphi)\right|~\overline{\Omega}^1_{p,q,m^{\prime},\eta^{\prime}}(d\varphi)
\end{array}$$
with the probability measure $~\overline{\Omega}^1_{p,q,m^{\prime},\eta^{\prime}}=\Omega^1_{q+1,\Gamma^2_{p,q}(m^{\prime},\eta^{\prime})}$.

Combining these two estimates, we arrive at the following inequality
$$
\begin{array}{l}
\left|\left[\Gamma_{p,n}^2(m^{\prime},\eta^{\prime})-\Gamma^2_{p,n}(m,\eta)\right](f)\right]
\\
\\
\leq a^2_{p,n} \displaystyle\int~\left|[\eta-\eta^{\prime}](\varphi)\right|~%
\Omega^2_{p,n,\eta^{\prime}}(f,d\varphi)
+\left[\sum_{p\leq q< n}~a^1_{p,q}~\overline{a}^2_{q,n}\right]~|m-m^{\prime}|\\
\\
+
\displaystyle\sum_{p\leq r<q<n}~\overline{a}^{1}_{r,q}~\overline{a}^2_{q,n}~\displaystyle\int~\left|[\Gamma_{p,r}^2(m^{\prime},\eta^{\prime})-\Gamma^2_{p,r}(m,\eta)](\varphi)\right|~\overline{\Omega}^1_{p,r,m^{\prime},\eta^{\prime}}(d\varphi)
\end{array}$$
This implies that
$$
\begin{array}{l}
\left|\left[\Gamma_{p,n}^2(m^{\prime},\eta^{\prime})-\Gamma^2_{p,n}(m,\eta)\right](f)\right]
\\
\\
\leq b^{\prime}_{p,n}~|m-m^{\prime}|+a^2_{p,n} \displaystyle\int~\left|[\eta-\eta^{\prime}](\varphi)\right|~%
\Omega^2_{p,n,\eta^{\prime}}(f,d\varphi)
\\
\\+
\displaystyle\sum_{p\leq r_1<n}~b_{r_1,n}
~\displaystyle\int~\left|[\Gamma_{p,r_1}^2(m^{\prime},\eta^{\prime})-\Gamma^2_{p,r_1}(m,\eta)](\varphi)\right|~\overline{\Omega}^1_{p,r_1,m^{\prime},\eta^{\prime}}(d\varphi)
\end{array}$$

Our next objective is to show that
$$
\begin{array}{l}
\left|\left[\Gamma_{p,n}^2(m^{\prime},\eta^{\prime})-\Gamma^2_{p,n}(m,\eta)\right](f)\right]
\\
\\
\leq \alpha^{k}_{p,n}~|m-m^{\prime}|+\beta^{k}_{p,n}
\displaystyle\int~\left|[\eta-\eta^{\prime}](\varphi)\right|~%
\Theta^{k}_{p,n,\eta^{\prime}}(f,d\varphi)
\\
\\+
\displaystyle
\sum_{p\leq r_1<r_2<\ldots<r_k<n}~b_{r_1,r_2}~\ldots~b_{r_k,n}
~\displaystyle\int~\left|[\Gamma_{p,r_1}^2(m^{\prime},\eta^{\prime})-\Gamma^2_{p,r_1}(m,\eta)](\varphi)\right|~\overline{\Omega}^1_{p,r_1,m^{\prime},\eta^{\prime}}(d\varphi)
\end{array}$$
for any $k\leq (n-p)$ for some Markov transitions $\Theta^{k}_{p,n,m^{\prime}\eta^{\prime}}(f,d\varphi)$ and the parameters
\begin{eqnarray*}
\alpha^{k}_{p,n}&=&b^{\prime}_{p,n}+\sum_{l=1}^{k-1}
 \sum_{p\leq r_1<\ldots r_{l}<n} b^{\prime}_{p,r_1}~b_{r_1,r_2}\ldots b_{r_{l},n}\\
\beta^k_{p,n}&=&a^{2}_{p,n}+\sum_{l=1}^{k-1}
 \sum_{p\leq r_1<\ldots r_{l}<n} a^{2}_{p,r_1}~b_{r_1,r_2}\ldots b_{r_{l},n}
\end{eqnarray*}
We proceed by induction on the parameter $k$. Firstly, we observe that the result is satisfied for $k=1$
with
$$
\left(\alpha^{1}_{p,n},\beta^1_{p,n}\right)=\left(b^{\prime}_{p,n},a^{2}_{p,n}\right)\quad\mbox{\rm and}
\quad \Theta^{1}_{p,n,\eta^{\prime}}=\Omega^2_{p,n,\eta^{\prime}}
$$
We further assume that the result is satisfied at rank $k$. In this situation, using the fact that
$$
\begin{array}{l}
\left|\left[\Gamma_{p,r_1}^2(m^{\prime},\eta^{\prime})-\Gamma^2_{p,r_1}(m,\eta)\right](\varphi)\right]
\\
\\
\leq b^{\prime}_{p,r_1}~|m-m^{\prime}|+a^2_{p,r_1} \displaystyle\int~\left|[\eta-\eta^{\prime}](\varphi^{\prime})\right|~%
\Omega^2_{p,r_1,\eta^{\prime}}(\varphi,d\varphi^{\prime})
\\
\\+
\displaystyle\sum_{p\leq r_0<r_1}~b_{r_0,r_1}
~\displaystyle\int~\left|[\Gamma_{p,r_0}^2(m^{\prime},\eta^{\prime})-\Gamma^2_{p,r_0}(m,\eta)](\varphi)\right|~\overline{\Omega}^1_{p,r_0,m^{\prime},\eta^{\prime}}(d\varphi)
\end{array}$$
we conclude that
$$
\begin{array}{l}
\left|\left[\Gamma_{p,n}^2(m^{\prime},\eta^{\prime})-\Gamma^2_{p,n}(m,\eta)\right](f)\right]
\\
\\
\leq \alpha^{k+1}_{p,n}~|m-m^{\prime}|+\beta^{k+1}_{p,n}
\displaystyle\int~\left|[\eta-\eta^{\prime}](\varphi)\right|~%
\Theta^{k+1}_{p,n,m^{\prime}\eta^{\prime}}(f,d\varphi)
\\
\\+
\displaystyle
\sum_{p\leq r_0<r_1<r_2<\ldots<r_k<n}~b_{r_0,r_1}~b_{r_1,r_2}~\ldots~b_{r_k,n}
~\displaystyle\int~\left|[\Gamma_{p,r_0}^2(m^{\prime},\eta^{\prime})-\Gamma^2_{p,r_0}(m,\eta)](\varphi)\right|~\overline{\Omega}^1_{p,r_0,m^{\prime},\eta^{\prime}}(d\varphi)
\end{array}$$
with
\begin{eqnarray*}
\alpha^{k+1}_{p,n}&=&\alpha^{k}_{p,n}+\sum_{p\leq r_1<r_2<\ldots<r_k<n}~b^{\prime}_{p,r_1}~b_{r_1,r_2}~\ldots~b_{r_k,n}\\
\beta^{k+1}_{p,n}&=&\beta^{k}_{p,n}+\sum_{p\leq r_1<r_2<\ldots<r_k<n}~a^{2}_{p,r_1}~b_{r_1,r_2}~\ldots~b_{r_k,n}
\end{eqnarray*}
and the Markov transition
\begin{eqnarray*}
\beta^{k+1}_{p,n}~\Theta^{k+1}_{p,n,m^{\prime}\eta^{\prime}}(f,d\varphi)&=&\beta^{k}_{p,n}~\Theta^{k}_{p,n,\eta^{\prime}}(f,d\varphi)\\
&&+\sum_{p\leq r_1<r_2<\ldots<r_k<n}~a^2_{p,r_1}~b_{r_1,r_2}~\ldots~b_{r_k,n}~\left(\overline{\Omega}^1_{p,r_1,m^{\prime},\eta^{\prime}}\Omega^2_{p,r_1,\eta^{\prime}}\right)(d\varphi)
\end{eqnarray*}
We end the proof of the proposition using the fact that
\begin{eqnarray*}
\left|\Gamma_{p,n}^1(m^{\prime},\eta^{\prime})-\Gamma_{p,n}^1(m,\eta)\right|
&\leq&
\left[ a^1_{p,n}
+\displaystyle\sum_{p\leq q<n}~c^{2,1}_{p,q}~\overline{a}^{1}_{q,n}
\right]
 ~|m-m^{\prime}|\\
 &&
 +\displaystyle\sum_{p\leq q<n}~\overline{a}^{1}_{q,n}~
 c^{2,2}_{p,q}
\displaystyle\int~\left|[\eta-\eta^{\prime}](\varphi^{\prime})\right|~%
\left[\overline{\Omega}^1_{p,q,m^{\prime},\eta^{\prime}}\Theta_{p,q,\eta^{\prime}}\right](d\varphi^{\prime})
\end{eqnarray*}
This proof of the proposition is now completed.
\cqfd
\subsection{Proof of theorem~\ref{theophdunif}}\label{ptheophdunif}

For any $\eta\in \Pa(E)$ and any $u,u^{\prime}\in I_n$, we have
$$
\begin{array}{l}
\left|\Phi^1_{n+1,\eta}(u)-\Phi^1_{n+1,\eta}(u^{\prime})\right|\\
\\
=\left|u-u^{\prime}\right|\left[
r(1-d)+rd h~\int \mathcal{Y}_{n}(dy)~
\frac{\eta(g(\point,y))}{%
[h+du\eta(g(\mbox{\LARGE .},y))][h+du^{\prime}\eta(g(\mbox{\LARGE .},y))]}
\right]\\
\\
\leq
\left|u-u^{\prime}\right|\left[
r(1-d)+rd h~\mathcal{Y}_{n}\left(
\frac{g^+}{%
[h+dm^- g^-)]^2}\right)
\right]
\end{array}
$$
This implies that condition (\ref{a1}) is satisfied with
$$
a^1_{n,n+1}\leq
r(1-d)+rd h~\mathcal{Y}_{n}\left(
\frac{g^+}{%
[h+dm^- g^-)]^2}\right)
$$
In the same way, for any $\eta,\eta^{\prime}\in \Pa(E)$ and any $u\in I_n$, we have
$$
\begin{array}{l}
\Phi^1_{n+1,\eta}(u)-\Phi^1_{n+1,\eta^{\prime}}(u)
=rdh u~\int \mathcal{Y}_{n}(dy)~
\frac{1}{%
[h+du\eta(g(\mbox{\LARGE .},y))][h+du\eta^{\prime}(g(\mbox{\LARGE .},y))]}~(\eta-\eta^{\prime})\left(g(\mbox{\LARGE .},y)\right)
\end{array}
$$
$$
\tau^1_{n+1}\leq rdh m^+~
\mathcal{Y}_{n}\left(
\frac{g^+-g^-}{%
[h+dm^-g^-]^2}\right)
$$
and the probability measure
$$
\Omega^1_{n,\eta^{\prime}}(d\varphi)\propto
\int \mathcal{Y}_{n}(dy)~
\frac{g^+(y)-g^-(y)}{%
[h+dm^-g^-(y)]^2}~\delta_{\frac{g(\mbox{\LARGE .},y)}{g^+(y)-g^-(y)}}(d\varphi)
$$
Now, we come to the analysis of the mappings
$$
\Phi^{2}_{n+1,u}(\eta)\propto r(1-d)u~\eta M +\int \Ya_n(dy)~w_u(\eta,y)~\Psi_{g(\point,y)}(\eta)M+ \mu(1)~\overline{\mu}
$$
with the weight functions
$$
w_u(\eta,y):=\frac{rdu\eta(g(\point,y))}{h+du\eta(g(\point,y))}=r\left(1-\frac{h}{h+du\eta(g(\point,y))}\right)
$$
Notice that
$$
w^-(y):=\frac{rdm^-g^-(y)}{h+dm^-g^-(y)}\leq w_u(\eta,y)\leq w^+(y):=\frac{rdm^+g^+(y)}{h+dm^+g^+(y)}
$$
To have a more synthetic formula, we extend the observation state space with two auxiliary points $c_1,c_2$ and we set $$
\Ya_n^c=\Ya_n+\delta_{c_1}+\delta_{c_2}
$$
we extend the likelihood and the weight functions  by setting
$$
g(x,c_1)=g(x,c_2)=1$$
and
\begin{eqnarray*}
w^-(c_1)&:=&r(1-d) m^-\leq w_{u}(\eta,c_1):=r(1-d)u \leq w^+(c_1):=r(1-d) m^+
\\
w_{u}(\eta,c_2)&=&w^+(c_2)=w^-(c_2):=\mu(1)
\end{eqnarray*}
In this notation, we find that
$$
\Phi^{2}_{n+1,u}(\eta)\propto \int \Ya_n^c(dy)~w_u(\eta,y)~\Psi_{g(\point,y)}(\eta)M_y
$$
with the collection of Markov transitions $M_y$ defined below
$$
\forall y\not\in \{c_2\}\qquad M_y=M\quad\mbox{\rm and}\quad
M_{c_2}=\overline{\mu}
$$
Notice that the normalizing constants $\Ya_n^c(w_u(\eta,\point))$  satisfy the following lower bounds
\begin{eqnarray*}
\Ya_n^c(w_u(\eta,\point))&\geq &\Ya_n^c(w^-)=r(1-d)~m^-+\Ya_n\left(w^-\right)+\mu(1)
\end{eqnarray*}

We analyze the Lipschitz properties of the mappings $\Phi^{2}_{n+1,u}$ using the following decomposition
$$
\Phi^{2}_{n+1,u}(\eta)-\Phi^{2}_{n+1,u}(\eta^{\prime})
=\Delta_{n+1,u}(\eta,\eta^{\prime})+\Delta_{n+1,u}^{\prime}(\eta,\eta^{\prime})
$$
with the signed measures
$$
\Delta_{n+1,u}(\eta,\eta^{\prime})=\int
\Ya_n^c(dy)~\frac{w_u(\eta,y)}{\Ya_n^c(w_u(\eta,\point))}~\left[
\Psi_{g(\point,y)}(\eta)M_y-\Psi_{g(\point,y)}(\eta^{\prime})M_y\right]
$$
and
$$
\Delta_{n+1,u}^{\prime}(\eta,\eta^{\prime})
=\frac{1}{\Ya_n^c(w_u(\eta,\point))}\int
\Ya_n^c(dy)~\left[w_u(\eta,y)-w_u(\eta^{\prime},y)\right]~\left(\Psi_{g(\point,y)}(\eta^{\prime})M_y-\Phi^{2}_{n+1,u}(\eta^{\prime})\right)
$$
Arguing as in the proof of theorem~\ref{theoBernoulli} given in the appendix (see for instance (\ref{refboltzman})), one checks that
$$
\begin{array}{l}
\left|\Delta_{n+1,u}(\eta,\eta^{\prime})(f)\right|\\
\\
\leq \frac{1}{\Ya_n^c(w^-)}~\left(
r(1-d)m^+~|(\eta-\eta^{\prime})(M(f))|+\int \Ya_n(dy)~w^+(y)~\frac{g^+(y)}{g^-(y)}~\left|(\eta-\eta^{\prime})(S^{y}_{\eta^{\prime}}M(f))\right|
\right)
\end{array}
$$
for some collection of Markov transitions $S^{y}_{\eta^{\prime}}$ from $E$ into itself. It is also readily checked that
$$
\left|\Delta_{n+1,u}^{\prime}(\eta,\eta^{\prime})(f)\right|\leq
\frac{hrdm^+}{\Ya_n^c(w^-)}~\int\Ya_n(dy)\frac{1}{(h+m^-dg^-(y))^2}
\left|(\eta-\eta^{\prime})(g(\point,y))\right|
$$
This clearly implies that condition (\ref{a2}) is satisfied with
$$
a^2_{n,n+1}\leq
\frac{1}{\Ya_n^c(w^-)}~\left(\beta(M)\left[
r(1-d)m^++\Ya_n\left(\frac{w^+g^+}{g^-}\right)\right]+hrdm^+\Ya_n\left(\frac{g^+-g^-}{(h+m^-dg^-)^2}\right)\right)
$$

We analyze the continuity properties of the mappings $u\mapsto \Phi^{2}_{n+1,u}(\eta)$ using the
following decomposition
$$
\begin{array}{l}
\Phi^{2}_{n+1,u}(\eta)-\Phi^{2}_{n+1,u^{\prime}}(\eta)
\\
\\=\frac{1}{\Ya_n^c(w_u(\eta,\point))}\int
\Ya_n^c(dy)~\left[w_u(\eta,y)-w_{u^{\prime}}(\eta,y)\right]~\left(\Psi_{g(\point,y)}(\eta)M_y-\Phi^{2}_{n+1,u^{\prime}}(\eta)\right)
\end{array}
$$
This implies that
$$
\begin{array}{l}
\left|\left[\Phi^{2}_{n+1,u}(\eta)-\Phi^{2}_{n+1,u^{\prime}}(\eta)\right](f)\right|
\leq \frac{1}{\Ya_n^c(w^-)}
\left[ r(1-d)+hrd~
\Ya_n\left(\frac{g^+}{(h+dm^-g^-)^2}\right)\right]
~|u-u^{\prime}|
\end{array}
$$
This shows that condition (\ref{tau2}) is satisfied with
$$
\tau^{2}_{n+1}\leq \frac{1}{\Ya_n^c(w^-)}
\left[ r(1-d)+hrd~
\Ya_n\left(\frac{g^+}{(h+dm^-g^-)^2}\right)\right]
$$
This ends the proof of the theorem.
\cqfd

\end{document}